\documentclass[11pt]{article}
\usepackage{amssymb,amsmath,pifont,amsfonts,amsthm,enumitem,mathrsfs}
\usepackage{mathtools}
\usepackage{hyperref}
\usepackage{titlesec}
 \usepackage{lipsum}
 \usepackage{lscape}
 \usepackage{adjustbox}
\usepackage{graphicx}
 \usepackage{todonotes}
\usepackage{multirow,boldline}
\usepackage[top=30mm,right=30mm,bottom=30mm,left=30mm]{geometry}
\usepackage{makecell}
\usepackage{array}
\newcolumntype{K}[1]{>{\centering\arraybackslash}p{#1}}
\usepackage{tikz}
\usetikzlibrary{shapes,fit,decorations.pathreplacing,calligraphy,spath3}
\usepackage{scalerel}
\usepackage{relsize}
\usepackage{tikz-cd}
\usepackage{rotating}
\usepackage{dynkin-diagrams}
\usepackage{longtable} 
\usepackage{scalerel}
\usepackage[tableposition=above,width=12in]{caption}
\usepackage{arydshln} 
\usetikzlibrary{positioning}
\usepackage{comment}
\usepackage{adjustbox}
\usepackage{pdflscape}
\usepackage{float}
\usepackage{enumitem}

\let\svthefootnote\thefootnote
\newcommand\freefootnote[1]{%
  \let\thefootnote\relax%
  \footnotetext{#1}%
  \let\thefootnote\svthefootnote%
}
\newcommand{\F}{\mathbb{F}}
\newcommand{\Z}{\mathbb{Z}}
\newcommand{\N}{\mathbb{N}}
\newcommand{\R}{\mathbb{R}}
\newcommand{\C}{\mathbb{C}}

\newcommand{\Di}{\mathcal{D}i}
\newcommand{\xmark}{\text{\ding{56}}}

\newcommand{\cH}{\mathcal{H}}
\newcommand{\cG}{\mathcal{G}}
\newcommand{\cO}{\mathcal{O}}
\newcommand{\cT}{\mathcal{T}}

\newcommand{\spn}[2][]{\left\langle #2 \right\rangle_{#1}}

\makeatletter
\newcommand{\imod}[1]{\allowbreak\mkern4mu({\operator@font mod}\,\,#1)}
\makeatother

\numberwithin{equation}{subsection}

\theoremstyle{plain}

\newtheorem{Theorem}{Theorem}[section]

\newtheorem{Lemma}[Theorem]{Lemma}
\newtheorem{Corollary}[Theorem]{Corollary}

\theoremstyle{definition}

\newtheorem{Definition}[Theorem]{Definition}

\newtheorem*{maintheorem!}{Theorem \ref{maintheorem!}}
\newtheorem*{classificationthm!}{Theorem \ref{classificationthm!}}
\newtheorem*{myconjecture!}{Conjecture \ref{myconjecture!}}
\newtheorem*{mydefinition!}{Definition \ref{mydefinition!}}
\newtheorem*{existence!}{Theorem \ref{existence!}}

\newcommand{\+}{\times}

\DeclareMathOperator{\Sp}{Sp}
\DeclareMathOperator{\cd}{cd}
\DeclareMathOperator{\Gal}{Gal}
\DeclareMathOperator{\PGL}{PGL}

\DeclareMathOperator{\Iso}{Iso}

\DeclareMathOperator{\Stab}{Stab}

\DeclareMathOperator{\Alt}{Alt}

\DeclareMathOperator{\SL}{SL}

\DeclareMathOperator{\GL}{GL}

\DeclareMathOperator{\GU}{GU}

\DeclareMathOperator{\rk}{rk}
\DeclareMathOperator{\sep}{s}

\begin{document}
\title{Maximal subgroups of maximal rank in the classical algebraic groups}
\author{Vanthana Ganeshalingam, Damian Sercombe, and Laura Voggesberger}

\date{\today}

\maketitle

\begin{abstract} 

\noindent Let $k$ be an arbitrary field. We classify the maximal reductive subgroups of maximal rank in any classical simple algebraic $k$-group up to a notion of equivalence called index-conjugacy, which was first introduced in \cite{S}. This result complements \cite[Thm.\ 2]{S}, which does the same for the exceptional groups. We determine which of these subgroups may be realised over a finite field, the real numbers, or over a $\mathfrak{p}$-adic field. We also look at the asymptotics of the number of such subgroups as the rank grows large. \freefootnote{2020 \textit{Mathematics Subject Classification}. Primary 20G15; Secondary 20G07.}
\end{abstract}

\section{Introduction}\label{introduction}

Algebraic groups have been studied since (at least) the 19th century. Methods of algebraic geometry were introduced to the theory in the 1950s by Borel, Chevalley and others, and were modernised in the 1960s. The structure theory of reductive groups forms a key part of the theory, which has many applications --  including to theoretical physics, and to number theory; for instance via the Langlands program. 

\vspace{2mm}\noindent Dating back to Sophus Lie in the late 19th century \cite{Lie}, algebraists have been interested in the program of classifying the subgroups -- and in particular, the maximal subgroups -- of reductive groups. 
Over the subsequent decades various authors have contributed to the classification of maximal subgroups of reductive groups over a field $k$; most notably \cite{C,Dy,Dy1,W} for the case $k=\C$, \cite{LS,LS1,LS3,Se,Se1,Te} for $k$ algebraically closed, \cite{As,LS,LS2,LS3,LSS} for $k$ finite, and \cite{BD,K,Ka,Ta} for $k=\R$. For $k$ algebraically closed, notable works in the broader program to understand all subgroups of reductive groups include \cite{LT, Ma, St, St1}.

\vspace{2mm}\noindent We follow the program introduced in \cite{S}, and work over an arbitrary field $k$. Our overarching goal is to classify the maximal subgroups of any absolutely simple $k$-group $G$ up to a notion of equivalence called index-conjugacy (which will be clarified later).

\vspace{2mm}\noindent In this paper we look at a large and interesting class of maximal subgroups of $G$; namely, those of maximal rank in $G$ (i.e. those which contain a maximal torus of $G$). Since $G$ is smooth, its maximal tori remain so after base change by an algebraic closure $\overline{k}$ of $k$. 
Thus the subgroups of maximal rank in $G$ are $\overline{k}/k$-forms of those subgroups of $G_{\overline{k}}$ which result from Borel-de Siebenthal's algorithm \cite{BD}. An example of a maximal subgroup of $G$ which is \textbf{not} of maximal rank is the copy of $\PGL_2$ irreducibly embedded in $G=\SL_3$ via the adjoint representation.

\vspace{2mm}\noindent Our main result is Theorem \ref{classthm}, in which we classify all index-conjugacy classes of isotropic maximal reductive subgroups of maximal rank in any classical absolutely simple $k$-group $G$ in terms of combinatorial data associated to their (Tits) indices. This result complements \cite[Thm.\ 2]{S}, which does the same for the exceptional absolutely simple $k$-groups. We also establish several consequences of this result. For instance, in Theorem \ref{classthmcor2} we show that the number of such subgroups grows subexponentially in the rank of $G$.

\vspace{2mm}\noindent In this paper we make two main innovations on top of the framework and methods introduced in \cite{S}. Firstly, we have written some computer code in Julia and Oscar \cite{OSCAR} to assist with and verify the calculations involved in proving Theorem \ref{classthm}. 
This code has been made publicly available on Github \hspace{-1.2mm}\footnote[1]{The code is publicly available at \url{https://github.com/voggesbe/MaximalSubgroups.jl}}. Secondly, we introduce some new asymptotic techniques in order to prove Theorem \ref{classthmcor2}. 

\vspace{2mm}\noindent Our setup henceforth is as follows. Let $k$ be an arbitrary field. Let $\overline{k}$ be an algebraic closure of $k$. Let $k^{\sep}$ be the separable closure of $k$ in $\overline{k}$, and let $\Gamma:=\Gal(k^{\sep}/k)$ be the absolute Galois group of $k$. Let $G$ be a reductive $k$-group (that is, $G$ is a smooth connected affine group scheme over $k$ whose geometric unipotent radical $\mathscr{R}_u(G_{\overline{k}})$ is trivial).

\vspace{2mm}\noindent Note our convention is that reductive groups are connected. By a \textit{maximal reductive} (resp. \textit{maximal connected}, \textit{maximal}) subgroup of $G$ we mean a reductive (resp. smooth connected, smooth) proper subgroup scheme of $G$ which is maximal amongst all such subgroup schemes.

\vspace{2mm}\noindent Associated to our reductive $k$-group $G$ is a combinatorial invariant $\mathcal{I}(G)$ called the \textit{(Tits) index} of $G$; this consists of the Dynkin diagram of $G$ along with some additional data (refer to Definition \ref{indexdef}). For example, the simple algebraic $\R$-group $\mathrm{SU}(5,1)$ has index $\smash{\begin{tikzpicture}[baseline=-0.6ex]\dynkin[scale=0.9, ply=2,fold radius = 2mm]{A}{o***o}\end{tikzpicture}}$. This invariant $\mathcal{I}(G)$ was introduced in the 1960s by Satake \cite{Sa} 
and Tits \cite{T1, T2}. The importance of this invariant is highlighted by the following key result \cite[Thm.\ 2.7.1]{T2}: a semisimple $k$-group $G$ is uniquely determined up to isogeny by its index $\mathcal{I}(G)$ along with the isogeny class of its anisotropic kernel (i.e. the centraliser of a maximal split torus of $G$).

\vspace{2mm}\noindent Over an arbitrary field $k$ we have a complete classification of all possible indices of semisimple $k$-groups (up to combinatorial isomorphism). This result is due to the following authors \cite{Sa,T0,T1,T2,T3,Sel,We}, and is summarised in \cite[Prop.\ 3.1.2, Table 2]{T2}. 

\vspace{2mm}\noindent Now let $H$ be a reductive subgroup of maximal rank in $G$. Associated to the pair $H \subset G$ is a combinatorial invariant called the \textit{embedding of indices} of $H \subset G$; this consists of the indices of $H$ and $G$, along with a map between them which satisfies certain compatibility conditions (refer to Definition \ref{embedindices}). This notion was introduced in \cite[Def.\ 1, p4]{S}, along with a notion of combinatorial isomorphism and conjugacy.

\vspace{2mm}\noindent The significance of this invariant is as follows. Two subgroups $H$ and $H'$ of $G$ are \textit{index-conjugate} if there exists $g \in G(k^{\sep})$ such that $H^g=H'$, $g^{-1}g^{\Gamma} \subset H(k^{\sep})$ and, for any maximal split torus $S$ of $H$, $S^g$ is a maximal split torus of $H'$. One can check that index-conjugacy is an equivalence relation on the set of subgroups of $G$ that is finer than $G(k^{\sep})$-conjugacy but coarser than $G(k)$-conjugacy. It was shown in \cite[Thm.\ 1(iii)]{S} that the index-conjugacy class of a maximal rank reductive subgroup $H$ of $G$ is uniquely determined by the conjugacy class of its embedding of indices. Consequently, the $G(k)$-conjugacy class of $H$ is uniquely determined by the conjugacy class of its embedding of indices along with the $G(k)$-conjugacy class of the anisotropic kernel of $H$.

\vspace{2mm}\noindent Our main result in this paper is the following.

\begin{Theorem}\label{classthm} Let $k$ be a field, let $G$ be an absolutely simple $k$-group of classical type, and let $H$ be an isotropic maximal reductive subgroup of maximal rank in $G$. Up to conjugacy, the embedding of indices of $H \subset G$ is one of those listed in Tables \ref{A_n}, \ref{B_n}, \ref{C_n}, \ref{D_4} and \ref{D_n} for the cases where $G$ is of type $A_n$, $B_n$, $C_n$, $D_4$, and $D_n$ for $n>4$ respectively.
\end{Theorem}

\begin{Corollary}\label{classthmcor} Let $(\mathcal{G},\mathcal{H},\theta)$ be
one of the embeddings of indices that are listed in Tables \ref{A_n}, \ref{B_n}, \ref{C_n}, \ref{D_4} and \ref{D_n}. If there is a $\checkmark$ in column $8$ (resp. $9$, $10$) of the corresponding row then $(\mathcal{G},\mathcal{H},\theta)$ may be realised as the embedding of indices of a pair of reductive groups over some field $k$ with cohomological dimension $1$ (resp. $k=\R$, $k$ is $\mathfrak{p}$-adic). If there is a $\xmark$ then this is impossible.
\end{Corollary}

\noindent For the next corollary, we use the following language. A \textit{special linear group} (resp. \textit{special unitary group}) is an inner (resp. outer) $k$-form of $\SL_n$ for some $n \geq 2$. A \textit{symplectic group} is a $k$-form of $\Sp_{2n}$ for some $n \geq 1$, and a \textit{special orthogonal group} is a $k$-form of $\mathrm{SO}_n$ for some $n \geq 3$. An absolutely simple $k$-group of classical type is isogenous to (at least) one of these groups. A \textit{simple component} refers to a non-commutative simple normal subgroup.

\vspace{2mm}\noindent Now let $G$ be an isotropic absolutely simple $k$-group of classical type, which is not of trialitarian type (i.e. it is not of type $^3\hspace{-0.5mm}D_4$ or $^6\hspace{-0.5mm}D_4$). 
Associated to $G$ is a finite-dimensional division algebra $D$ over $k$ (refer to \cite[Ch.\ 24.k]{Mi} or \cite[Table II]{T2}). Let $d$ be the degree of $D$ over its center, that is, $d=\sqrt{[D:Z(D)]}$. We say that $G$ \textit{has degree $d$}. [N.B.  There may be several of such division algebras naturally associated to $G$, which arise from diagram automorphisms or low rank exceptional isomorphisms. If the rank of $G$ is at least 5 then its degree $d$ is uniquely determined, however there are a few low rank cases where $G$ has two different degrees.] 

\begin{Corollary}\label{classthmcor1} Let $k$ be a field. Let $G$ be an absolutely simple $k$-group of classical type, which is not of trialitarian type. Let $H$ be a maximal reductive subgroup of maximal rank in $G$. Suppose $H$ contains at least two simple components, all of which are isotropic. Then there exists a positive integer $d$ such that every simple component of $H$, along with $G$, is isogenous to either:
\begin{enumerate}[label=(\alph*)]
\vspace{-0.5mm}\item a special unitary group which has degree $d$,
\vspace{-1.5mm}\item a symplectic group which has degree $d$, or
\vspace{-1.5mm}\item a special orthogonal group which has degree $d$.
\end{enumerate} Moreover, if $G$ is of type $B_n$ then $d=1$.

\end{Corollary}

\noindent We conclude by establishing the following asymptotic result.

\begin{Theorem}\label{classthmcor2} Let $k$ be a field. Let $G$ be an absolutely simple $k$-group, and let $n$ be the rank of $G$. The number of index-conjugacy classes of isotropic maximal connected subgroups of maximal rank in $G$ is of the order $O(n^{c\log n})$, where $c$ is an absolute constant. If $k$ is separably closed then this bound may be sharpened to $O(n)$.
\end{Theorem}




\section{Preliminaries}

\noindent In this section we setup, define, and recall some basic properties of the index of a reductive group, and subsequently the embedding of indices of a suitable pair of reductive groups. We then present purely combinatorial analogues of these notions.

\vspace{2mm}\noindent The constructions in this section are quite technical, a necessary evil to obtain all the combinatorial data we want. We follow \cite[\S 1.4, \S 1.5, Ch.\ 2, p28]{S}.

\subsection{The index}\label{indexsubsection}

\noindent Let $k$ be an arbitrary field. Let $k^{\sep}$ be a separable closure of $k$, and let $\Gamma:=\Gal(k^{\sep}/k)$ be the absolute Galois group. 

\vspace{2mm}\noindent Let $G$ be a (connected) reductive $k$-group. We define the index of $G$ as follows.

\vspace{2mm}\noindent Fix a maximal split torus $S$ of $G$, and a maximal torus $T$ of $G$ which contains $S$. Consider the character group $X(T)$, and let $\iota:\Gamma \to \GL(X(T)_{\R})$ be the natural action of $\Gamma$ on $X(T)$. By \cite[Cor.\ 1.16]{S} we may equip $X(T)$ with a $\Gamma$-order $<$. [Recall from \cite[\S 1.3]{S} that a \textit{$\Gamma$-order} on $X(T)$ is a translation-invariant total order $<$ on $X(T)$ such that the subset $\{\chi \in X(T)\hspace{0.5mm}|\hspace{0.5mm} \chi>0, \chi|_S \neq 0\}$ of $X(T)$ is stable under $\iota(\Gamma)$.] 

\vspace{2mm}\noindent Let $\Lambda$ be the (unique) system of simple roots for $G$ with respect to $T$ that is compatible with $<$. Let $\Lambda_0$ be the subset of $\Lambda$ that vanishes on $S$. Let $W$ be the Weyl group of $G$ with respect to $T$. For every $\sigma \in \Gamma$ there exists a unique $w_{\sigma} \in W$ such that $w_{\sigma}\iota(\sigma)(\Lambda)=\Lambda$ (see for instance \cite[Lem.\ 1.17]{S}). The \textit{Tits action} $\hat{\iota}:\Gamma \to \GL(X(T)_{\R})$ is defined by $\hat{\iota}(\sigma):=w_{\sigma} \iota(\sigma)$; by construction it stabilises $\Lambda$.

\begin{Definition}\label{indexdef} The tuple $\big(X(T)_{\R}, \Lambda, \Lambda_0, \hat{\iota}(\Gamma)\big)$ is called the \textit{index} $\mathcal{I}(G)$ of $G$ (with respect to $T$ and $<$).
\end{Definition}

\noindent The index of $G$ is independent of the choice of $T$ and $<$, up to the obvious notion of combinatorial isomorphism (see \cite[Thm.\ 2.5]{S}).

\vspace{2mm}\noindent We illustrate the index $\mathcal{I}(G)$ of $G$ as follows. If $G$ is semisimple then we illustrate $\mathcal{I}(G)$ using a \textit{Tits-Satake diagram}, which is constructed by taking the Dynkin diagram of $G$, blackening each vertex in $\Lambda_0$ and linking with a solid gray bar all of the $\hat{\iota}(\Gamma)$-orbits of $\Lambda$. For arbitrary reductive $G$ we illustrate $\mathcal{I}(G)$ using the notation $\mathcal{I}(\mathscr{D}(G)) \times \mathcal{T}^d_c$; where $d$ is the dimension of $Z(G)$, $c$ is the dimension of the maximal split torus of $Z(G)$, $\mathscr{D}(G)$ is the (semisimple) derived subgroup of $G$, and $\mathcal{I}(\mathscr{D}(G))$ is represented by its Tits-Satake diagram. For example, we illustrate the index of the general unitary group $\GU_3$ (over a suitable field) by $\begin{tikzpicture}[baseline=-0.3ex]\dynkin[ply=2,rotate=90,fold radius = 2mm]{A}{oo}\end{tikzpicture} \times \mathcal{T}^1_0$. In the cases where it is not clear from this illustration, we also explicitly describe the subgroup $\hat{\iota}(\Gamma)$ of $\GL((X(T)_{\R})$.

\subsection{The embedding of indices}\label{embedindicessubsection}

\noindent Continue to use the setup from $\S \ref{indexsubsection}$. Let $H$ be a reductive subgroup of maximal rank in $G$. We define the embedding of indices of the pair $H \subset G$ as follows.

\vspace{2mm}\noindent Fix a maximal split torus $S_H$ of $H$, and a maximal torus $T_H$ of $H$ which contains $S_H$. Recall from \cite[Thm.\ 20.9(ii)]{Bo} that all maximal split tori of $G$ are conjugate by an element of $G(k)$. It follows that there exists $g \in G(k^{\sep})$ such that $(S_H)^{g} \subseteq S$ and $(T_H)^g=T$. Our total order $<$ on $X(T)$ induces a total order $<_{g^{-1}}$ on $X(T_H)$ that is given by $\chi_1 <_{g^{-1}} \chi_2$ if $\chi_1 \circ g^{-1} < \chi_2 \circ g^{-1}$ for $\chi_1, \chi_2 \in X(T_H)$. By \cite[Cor.\ 1.27]{S} we can ensure that $<$ has been chosen such that the induced total order $<_{g^{-1}}$ on $X(T_H)$ is also a $\Gamma$-order (with respect to $S_H$). We define a linear map $\theta:X(T_H)_{\R} \to X(T)_{\R}$ by $\chi \mapsto \chi \circ g^{-1}$. 

\vspace{2mm}\noindent Let $\Delta$ be the (unique) system of simple roots for $H$ with respect to $T_H$ that is compatible with $<_{g^{-1}}$. Let $\Delta_0$ be the subset of $\Delta$ that vanishes on $S_H$. Let $\hat{\iota}_H:\Gamma \to \GL(X(T_H)_{\R})$ be the associated Tits action. The index $\mathcal{I}(H)$ of $H$ is $\big(X(T_H)_{\R}, \Delta, \Delta_0, \hat{\iota}_H(\Gamma)\big)$. 

\begin{Definition}\label{embedindices} The triple $\big(\mathcal{I}(G),\mathcal{I}(H),\theta\big)$ is called the \textit{embedding of indices} of $H \subset G$ (with respect to $T$, $T_H$, $<$ and $g$).
\end{Definition}

\noindent It turns out that the embedding of indices of $H \subset G$ is independent of the choice of $T$, $T_H$, $<$ and $g$ up to a notion of combinatorial isomorphism, this is proved in \cite[Thm.\ 1]{S}.

\vspace{2mm}\noindent We describe (the combinatorial isomorphism class of) an embedding of indices $\big(\mathcal{I}(G),\mathcal{I}(H),\theta\big)$ via the aforementioned illustrations for $\mathcal{I}(G)$ and $\mathcal{I}(H)$, along with providing the necessary subsystem data required to uniquely determine $\theta$.

\subsection{Almost primitive subsystems}\label{almostprimsubsystems}

\noindent In this section we recall and discuss the notion of an almost primitive subsystem of a root system. We follow \cite[\S 2.1]{S}.

\vspace{2mm}\noindent Let $\Phi$ be a (reduced, crystallographic, abstract) root system. Let $W$ be its Weyl group, and $\Iso(\Phi)$ its isometry group. Let $p$ be either zero or a prime number.

\vspace{2mm}\noindent Associated to $\Phi$ is a set of structure constants $c_{\alpha\beta}^{mn} \in \Z$ for $\alpha, \beta \in \Phi$ and integers $m,n>0$ (for instance see $\S 11.1$ of \cite{MT}). A closed subsystem $\Psi$ of $\Phi$ is \textit{$p$-closed} if, for $\alpha, \beta \in \Psi$ and integers $m, n > 0$, we have $m\alpha+n\beta \in \Psi$ whenever $c_{\alpha\beta}^{mn}$ is non-zero mod $p$.

\vspace{2mm}\noindent Let $\Psi$ be a $p$-closed proper subsystem of $\Phi$. Associated to the pair $\Psi \subset \Phi$ are stabiliser groups $\Stab_W(\Delta)$ and $\Stab_{\Iso(\Phi)}(\Delta)$; where $\Delta$ refers to a base of $\Psi$ (these stabiliser groups are independent of the choice of $\Delta$). 

\vspace{2mm}\noindent Let $W_{\Psi}$ be the Weyl group of $\Psi$, consider it as a subgroup of $W$. If $\Psi$ is maximal amongst $p$-closed subsystems of $\Phi$ that are invariant under the normaliser $N_W(W_{\Psi})$ (resp. $N_{\Iso(\Phi)}(W_{\Psi})$) then $\Psi$ is \textit{$p$-primitive} (resp. \textit{almost $p$-primitive}) in $\Phi$. If $\Psi$ is $p$-primitive (resp. almost $p$-primitive) in $\Phi$ for some $p \geq 0$ then we say $\Psi$ is \textit{primitive} (resp. \textit{almost primitive}) in $\Phi$.

\vspace{2mm}\noindent The following classification of non-empty almost primitive subsystems is well-known, see for instance \cite[Table\ 1]{K} or \cite[Lem.\ 1, Table\ 1]{S} (however \cite[Table\ 1]{S} erroneously omits the subsystem $A_{n-1} \subset C_n$ for $p \neq 2$). 

\begin{Lemma}[Table 1, Lem.\ 2.1 of \cite{S}]\label{isoautgps} Let $\Phi$ be an irreducible root system of classical type, with Weyl group $W$. Let $\Psi$ be a non-empty almost primitive subsystem of $\Phi$. The possibilities for the pair $\Psi \subset \Phi$ are classified up to $W$-conjugacy in Table \ref{classlist}, along with their associated groups $\Stab_W(\Delta)$. Moreover, $\Stab_{\Iso(\Phi)}(\Delta)=\Stab_W(\Delta)$ unless one of the following occurs.

\vspace{2mm}\noindent $(i)$ If $\Phi=D_4$ and $\Psi=(A_1)^4$ (resp. $A_3$, $A_2$) then $\Stab_{\Iso(\Phi)}(\Delta)=S_4$ (resp. $(\Z_2)^2$, $\Z_2 \times S_3$).

\vspace{1mm}\noindent $(ii)$ If $\Phi=D_n$ ($n>4$ is even) and $\Psi=D_mD_{n-m}$ (resp. $(D_l)^{n_l}$) then $\Stab_{\Iso(\Phi)}(\Delta)=(\Z_2)^2$ (resp. $(\Z_2)^{n_l} \rtimes S_{n_l}$).

\vspace{2mm}\noindent $(iii)$ If $\Phi=A_n$ ($n>1$) or $D_n$ ($n$ is odd) then $\Stab_{\Iso(\Phi)}(\Delta)=\Z_2 \times \Stab_W(\Delta)$.
\end{Lemma}

\begin{table}[H]\centering\caption{Almost primitive subsystems $\Psi$ of classical type $\Phi$}\label{classlist}
\begin{tabular}{c | c c } 
$\Phi$ & $\Psi$ & $\Stab_W(\Delta)$  \\ [0.2ex] \hline 
$A_n$ & $A_mA_{n-m-1}$, $(A_l)^{n_l}$ & $1$, $S_{n_l}$ \topstrut \\ 
$B_n$ ($p \neq 2$) & $D_mB_{n-m}$, $B_{n-1}$ & $\Z_2$, $\Z_2$ \topstrut \\ 
$B_n$ ($p=2$) & $B_mB_{n-m}$, $D_n$, $(B_l)^{n_l}$ & $1$, $\Z_2$, $S_{n_l}$ \topstrut \\ 
$C_n$ & $C_mC_{n-m}$, $(C_l)^{n_l}$, $A_{n-1}$ ($p\neq 2$), $\widetilde{D_n}$ ($p=2$) & $1$, $S_{n_l}$, $\Z_2$, $\Z_2$ \topstrut \\ 
$D_n$ ($n>4$) & $A_{n-1}$ ($\hspace{0.5mm}(n,2)$ classes), $D_mD_{n-m}$, $(D_l)^{n_l}~$ & $~\Z_{(n,2)}$, $\Z_2$, $(\Z_2)^{n_l-1} \hspace{-0.5mm}\rtimes\hspace{-0.5mm} S_{n_l}$ \topstrut \\ 
$D_4$ & $(A_1)^4$, $A_3$ ($3$ classes), $A_2$ & $(\Z_2)^2$, $\Z_2$, $\Z_2$ \topstrut \\ 
\end{tabular}
\end{table}

\vspace{1mm}\noindent The parameters in Table \ref{classlist} are positive integers which are subject to the following constraints: $n_l \geq 2$, $ln_l=n$ (except for the case $(A_l)^{n_l} \subset A_n$ where $(l+1)n_l=n+1$), and $1 \leq m < n/2$ (except for $A_mA_{n-m-1} \subset A_n$ where $0 \leq m < (n-1)/2$, and $D_mB_{n-m} \subset B_n$ where $2 \leq m <n$ if $p \neq 2$ as well as $m=n$ for all $p$).

\subsection{Abstract indices}

\noindent In this section we recall the definition of an abstract index introduced in \cite[Def.\ 2.4]{S}. This was never formally defined in the work of Satake \cite{Sa} and Tits \cite{T1, T2}, yet they implicitly used something very similar.

\vspace{2mm}\noindent Let $E$ be a finite-dimensional real inner product space. Let $\Delta$ be a system of simple roots in $E$ (we require that $\Delta$ is reduced, but do not require that $\Delta$ spans $E$). Let $\Delta_0$ be a subset of $\Delta$. Let $\Pi$ be a subgroup of the isometry group of $E$ that stabilises both $\Delta$ and $\Delta_0$.

\begin{Definition}[Def.\ 2.4 of \cite{S}]\label{abstractindexdefn} The quadruple $\mathcal{I}:=(E,\Delta,\Delta_0,\Pi )$ is an \textit{abstract index} if the following five properties are satisfied: 

\vspace{2mm}\noindent $(i)$ \textit{(relative root system)} The projection of $\Delta$ onto $E_s$ is a system of simple roots for a not-necessarily-reduced abstract root system; where $E_s$ denotes the subspace of $E^{\Pi}$ which is perpendicular to $\Delta_0$.

\vspace{1mm}\noindent $(ii)$ \textit{(self-opposition)} The involution $-w_0$ of $\Delta$ commutes with $\Pi$ and stabilises $\Delta_0$; where $w_0$ denotes the longest element of the Weyl group of $\Delta$.

\vspace{1mm}\noindent $(iii)$ \textit{(classical)} If $\Delta$ is irreducible of classical type, $\Delta \setminus \Delta_0$ has a unique $\Pi$-orbit and $\Delta_0$ contains an irreducible component of type $A_{l-1}$ for some $l \geq 2$ then, for some $n \geq 2$, either $\Delta$ is of type $A_{n-1}$ and $l$ divides $n$, or $\Delta$ is of type $C_n$ or $D_n$ and $l$ is a power of $2$.

\vspace{1mm}\noindent $(iv)$ \textit{(exceptional)} If $\Delta$ is irreducible of exceptional type, $\Pi$ is trivial and $\Delta_0$ is irreducible then $\Delta_0$ is not of type $A_l$ or $C_l$ for any $l \in \N$.

\vspace{1mm}\noindent $(v)$ \textit{(inductive closure)} If one removes a $\Pi$-orbit $\mathcal{O}$ from $\Delta \setminus \Delta_0$ then the result $(E,\Delta \setminus \mathcal{O},\Delta_0,\Pi )$ is again an abstract index.
\end{Definition}

\vspace{1mm}\noindent Let $\mathcal{I}=(E,\Delta,\Delta_0,\Pi)$ be an abstract index. If $\Delta=\Delta_0$ and $E_s$ is trivial then $\mathcal{I}$ is \textit{anisotropic}, otherwise $\mathcal{I}$ is \textit{isotropic}. 
If $\Delta$ spans $E$ and $\Delta$ is irreducible (resp. irreducible of classical type, irreducible of exceptional type) then $\mathcal{I}$ is \textit{irreducible} (resp. \textit{of classical type}, \textit{of exceptional type}). If $\Pi$ is trivial then $\mathcal{I}$ is of \textit{inner type}, otherwise $\mathcal{I}$ is of \textit{outer type}. If $\Delta_0 = \varnothing$ then $\mathcal{I}$ is \textit{quasisplit}. If $\Delta_0=\varnothing$ and $\Pi$ is trivial then $\mathcal{I}$ is \textit{split}.

\vspace{2mm}\noindent We use the following language. The \textit{type} of $\mathcal{I}$ is the isomorphism class of $\Delta$ as a root system (sometimes depending on context we also include the order of $\Pi$). For example, we might say that $\SL_3$ is of type $A_2$, or that it is of type ${}^1\!A_2$.

\vspace{2mm}\noindent Now let $\smash{'}\mathcal{I}=(\smash{'}\hspace{-0.3mm}E,\smash{'}\hspace{-0.4mm}\Delta,\smash{'}\hspace{-0.4mm}\Delta_0,\smash{'}\Pi)$ be another abstract index. An \textit{isomorphism} from $\mathcal{I}$ to $\smash{'}\mathcal{I}$ is a bijective isometry $\psi:E \to \smash{'}\hspace{-0.3mm}E$ that satisfies $\psi(\Delta)=\smash{'}\hspace{-0.4mm}\Delta$, $\psi(\Delta_0)=\smash{'}\hspace{-0.4mm}\Delta_0$ and $\psi\Pi \psi^{-1}=\smash{'}\Pi$. We denote $\psi(\mathcal{I}):=\smash{'}\mathcal{I}$.

\vspace{2mm}\noindent Given a field $k$, the abstract index $\mathcal{I}$ is called \textit{$k$-admissible} if there exists a reductive $k$-group $G$ such that the index of $G$ is isomorphic to $\mathcal{I}$. For example, the abstract index $\dynkin{A}{*o*}$ is $\R$-admissible but not $\F_q$-admissible for any finite field $\F_q$.

\subsection{Embeddings of abstract indices}\label{embedabstractindices-section}

\noindent In this section we present a combinatorial construction called an embedding of abstract indices; it was first introduced in \cite[Def.\ 1]{S}. It should be interpreted as a ``relative" version of an abstract index; i.e. a means of embedding one abstract index into another.

\vspace{2mm}\noindent On first glance the following definition appears very complicated, however these combinatorial objects become much clearer when one represents them using a pair of Tits-Satake diagrams along with an ``orientation" -- as in the Appendix.

\begin{Definition}[Def.\ 1 of \cite{S}]\label{embedabstractindices} Let $p$ either be a prime number or $0$. A \textit{$p$-embedding of abstract indices} consists of two abstract indices $\mathcal{G}=(F,\Lambda,\Lambda_0, \Xi)$ and $\mathcal{H}=(E,\Delta,\Delta_0,\Pi)$, and a bijective isometry $\theta:E \to F$ which satisfies the following conditions:

\vspace{1.5mm}\noindent [\underline{Notation}: $\Phi:=\langle \Lambda \rangle$ (i.e. the root system in $F$ generated by $\Lambda$), $\Phi^+$ is the associated positive roots, $\Psi:=\langle \theta(\Delta) \rangle$, $\Psi_0:=\langle \theta(\Delta_0) \rangle$, $E_a$ is the span of $\Delta_0$ and $(E^{\Pi})^{\perp}$, 
$\Lambda_a:= \Lambda \cap \theta(E_a)$, $\Phi_a:=\Phi \cap \theta(E_a)$, $W$ (resp. $W_{\Lambda_a}$) is the Weyl group of $\Lambda$ (resp. $\Lambda_a$), ${}_{in}\Lambda_a$ is the union of all irreducible components of $\Lambda_a$ that are contained in $\Psi_0$.]

\vspace{3.5mm}\noindent $(A.1)$ $\Psi$ is a $p$-closed subsystem of $\Phi$, $\theta(\Delta) \subset \Phi^+$, $\Lambda_a$ is a base of $\Phi_a$ and ${}_{in}\Lambda_a \subseteq \theta(\Delta_0)$.

\vspace{1.5mm}\noindent $(A.2)$ For every $\sigma \in \Pi$ there exists a unique $w_{\sigma} \in W_{\Lambda_a}$ such that $w_{\sigma} \sigma^{\theta} \in \Xi$. Moreover, the map $\Pi \to \Xi$ given by $\sigma \mapsto w_{\sigma} \sigma^{\theta}$ is surjective.

\vspace{1.5mm}\noindent $(A.3)$ $\Lambda_a$ is $\Xi$-stable and contains $\Lambda_0$.

\vspace{1.5mm}\noindent $(A.4)$ ${}_{in}\Lambda_a$ is contained in $\Lambda_0$.
\end{Definition}

\vspace{1mm}\noindent The map $\theta$ is called the \textit{embedding} of $\mathcal{H}$ in $\mathcal{G}$. We sometimes abuse notation and identify $\mathcal{H}$ with its image under $\theta$ (so we interpret $\theta$ as just inclusion). An \textit{embedding of abstract indices} is a $p$-embedding of abstract indices for some $p \geq 0$. 

\vspace{2mm}\noindent Let $(\mathcal{G},\mathcal{H},\theta)$ be a $p$-embedding of abstract indices. If $\mathcal{H}$ is isotropic then $\mathcal{G}$ is automatically also isotropic by $(A.3)$; if this condition holds we say $(\mathcal{G},\mathcal{H},\theta)$ is \textit{isotropic}. Otherwise, $(\mathcal{G},\mathcal{H},\theta)$ is called \textit{anisotropic}. 
If $\mathcal{H}$ and $\mathcal{G}$ are both split (resp. quasisplit) then $(\mathcal{G},\mathcal{H},\theta)$ is called \textit{split} (resp. \textit{quasisplit}). If ${}_{in}\Lambda_a=\Lambda_a$ then $(\mathcal{G},\mathcal{H},\theta)$ is \textit{independent}.

\vspace{2mm}\noindent Let $V$ denote the largest subspace of $F$ that is contained in the span of $\Phi$ and is perpendicular to $\Psi$. If the fixed point subspace of $V$ under the action of $\Pi^{\theta}$ is trivial and $\Psi$ is maximal among $p$-closed $\Pi^{\theta}$-stable subsystems of $\Phi$ then $(\mathcal{G},\mathcal{H},\theta)$ is \textit{maximal}. [N.B. This condition only depends on the pair $(\mathcal{H},\Phi)$, so we will often use the synonym $\mathcal{H}$ \textit{is maximal in} $\Phi$.]

\vspace{2mm}\noindent Let $(\smash{'}\mathcal{G},\smash{'}\mathcal{H},\smash{'}\hspace{-0.1mm}\theta)$ be another $p$-embedding of abstract indices; denote the ambient vector space of $\smash{'}\mathcal{G}$ by $\smash{'}\hspace{-0.5mm}F$. An \textit{isomorphism} from $(\mathcal{G},\mathcal{H},\theta)$ to $(\smash{'}\mathcal{G},\smash{'}\mathcal{H},\smash{'}\hspace{-0.1mm}\theta)$ is a bijective isometry $\phi:F \to \smash{'}\hspace{-0.5mm}F$ such that $\phi(\mathcal{G})=\smash{'}\mathcal{G}$ and $\phi(\theta(\mathcal{H}))=\smash{'}\hspace{-0.1mm}\theta(\smash{'}\mathcal{H})$. If such an isomorphism exists then $(\mathcal{G},\mathcal{H},\theta)$ and $(\smash{'}\mathcal{G},\smash{'}\mathcal{H},\smash{'}\hspace{-0.1mm}\theta)$ are \textit{isomorphic} (a.k.a. $(\mathcal{H},\theta)$ and $(\smash{'}\mathcal{H},\smash{'}\hspace{-0.1mm}\theta)$ are \textit{isomorphic in $\mathcal{G}$}).

\vspace{2mm}\noindent Now identify $F=\smash{'}\hspace{-0.5mm}F$. Recall that $W$ denotes the Weyl group of $\Phi$. An element $w \in W$ is a \textit{conjugation} from $(\mathcal{H},\theta)$ to $(\smash{'}\mathcal{H},\smash{'}\hspace{-0.1mm}\theta)$ in $\Phi$ if $w(\theta(\mathcal{H}))=\smash{'}\hspace{-0.1mm}\theta(\smash{'}\mathcal{H})$. If in addition $w(\mathcal{G})=\smash{'}\mathcal{G}$ then $w$ is a \textit{conjugation} from $(\mathcal{H},\theta)$ to $(\smash{'}\mathcal{H},\smash{'}\hspace{-0.1mm}\theta)$ in $\mathcal{G}$. If such a conjugation in $\Phi$ (resp. $\mathcal{G}$) exists then $(\mathcal{H},\theta)$ and $(\smash{'}\mathcal{H},\smash{'}\hspace{-0.1mm}\theta)$ are \textit{conjugate} in $\Phi$ (resp. $\mathcal{G}$).

\vspace{2mm}\noindent Given a field $k$ with characteristic $p$, the $p$-embedding of abstract indices $(\mathcal{G},\mathcal{H},\theta)$ is called \textit{$k$-admissible} if there exists a pair of reductive $k$-groups $H \subset G$ such that the embedding of indices of $H \subset G$ is isomorphic to $(\mathcal{G},\mathcal{H},\theta)$.

\section{Proofs}\label{proofs}

\noindent In this section we prove Theorem \ref{classthm}, its Corollaries \ref{classthmcor} and \ref{classthmcor1}, and also Theorem \ref{classthmcor2}.

\vspace{2mm}\noindent Our strategy of proof for Theorem \ref{classthm} follows that of \cite[Thm.\ 2]{S}. It involves constructing Tables \ref{A_n}, \ref{B_n}, \ref{C_n}, \ref{D_4} and \ref{D_n} for the cases $A_n$, $B_n$, $C_n$, $D_4$, and $D_n$ for $n>4$ respectively. We present these tables along with some associated computations in the Appendix.

\vspace{2mm}\noindent\underline{Proof of Theorem \ref{classthm}.}

\vspace{2mm}\noindent We rephrase the statement of Theorem \ref{classthm} in combinatorial terms. Let $\mathbb{P}$ denote the set of prime numbers.

\begin{Theorem}\label{classthmcomb} \textit{For each $p \in \mathbb{P} \cup \{0\}$ and each (isomorphism class of) abstract index $\mathcal{G}$ of classical type, all $\mathcal{G}$-conjugacy classes of isotropic maximal $p$-embeddings of abstract indices are classified in Tables \ref{A_n}, \ref{B_n}, \ref{C_n}, \ref{D_4} and \ref{D_n} for the cases where $\mathcal{G}$ is of type $A_n$, $B_n$, $C_n$, $D_4$, and $D_n$ for $n>4$ respectively.}
\end{Theorem}

\noindent We use the setup and notation of $\S \ref{embedabstractindices-section}$. That is, consider any $p$-embedding of abstract indices $(\mathcal{G},\mathcal{H},\theta)$ for some $p \in \mathbb{P} \cup \{0\}$. We abuse notation and identify $\mathcal{H}$ with its image under $\theta$. 
We denote $\mathcal{G}=(E,\Lambda,\Lambda_0,\Xi)$, $\mathcal{H}=(E,\Delta,\Delta_0,\Pi)$, $\Phi:=\langle \Lambda \rangle$, $\Phi_0:=\langle \Lambda_0 \rangle$, $W$ is the Weyl group of $\Phi$, $\Iso(\Phi)$ is the isometry group of $\Phi$, $\Psi:=\langle \Delta \rangle$, $\Psi_0:=\langle \Delta_0 \rangle$, $E_s$ is the subspace of $E^{\Pi}$ which is perpendicular to $\Delta_0$, $\Phi_s := \Phi \cap E_s$, $E_a:=(E_s)^{\perp}$, $\Lambda_a:= \Lambda \cap E_a$, $\Phi_a := \Phi \cap E_a$, $E_{\Lambda}$ (resp. $E_{\Delta}$, $E_{\Delta_0}$) is the span of $\Lambda$ (resp. $\Delta$, $\Delta_0$) in $E$, $\overline{E_a}:=\smash{(E_{\Delta_0})^{\perp}} \cap E_a$, $V:=E_{\Lambda} \cap (E_{\Delta})^{\perp}$, and ${}_{in}\Phi_a$ is the union of all irreducible components of $\Phi_a$ that are contained in $\Psi_0$. By assumption $\Phi$ spans $E$.

\vspace{2mm}\noindent A note to the reader: our tables may be read as follows. Given a choice of parameters, each row uniquely corresponds to an isomorphism class of an embedding of abstract indices. We will show that a single $\mathcal{G}$-isomorphism class splits into multiple distinct $\mathcal{G}$-conjugacy classes only in very rare cases: which only occur when $\Psi \subset \Phi$ is either $A_{n-1} \subset D_n$ or $D_m^{n_m} \subset D_n$, for $n\geq 4$ even. For example, when $n=4$ this phenomenon occurs only for the embeddings $$\begin{tikzpicture}[baseline=-0.3ex]
\dynkin[ply=2,rotate=90,fold radius = 2mm]{A}{o*o}\end{tikzpicture}\times\mathcal{T}^1_0 \subset \dynkin{D}{*o**},~ \begin{tikzpicture}[baseline=-0.3ex]
\dynkin[ply=2,rotate=90,fold radius = 2mm]{A}{*o*}\end{tikzpicture}\times\mathcal{T}^1_0 \subset \dynkin{D}{oooo},~ \begin{tikzpicture}[baseline=-0.3ex]
\dynkin[ply=2,rotate=90,fold radius = 2mm]{A}{ooo}\end{tikzpicture}\times\mathcal{T}^1_0 \subset \dynkin{D}{oooo},~ \begin{tikzpicture}[baseline=0.5ex] 
\dynkin[name=first]{A}{o}
\node (current) at ($(first root 1)+(-0,3.4mm)$) {};
\dynkin[at=(current),name=second]{A}{o}
\node (currenta) at ($(first root 1)+(-3.4mm,0)$) {};
\dynkin[at=(currenta),name=third]{A}{o}
\node (currenta) at ($(third root 1)+(-0,3.4mm)$) {};
\dynkin[at=(currenta),name=fourth]{A}{o}
\begin{scope}[on background layer]
\draw[/Dynkin diagram/fold style]
($(first root 1)$) -- ($(third root 1)$);%
\draw[/Dynkin diagram/fold style]
($(second root 1)$) -- ($(fourth root 1)$);%
\end{scope}\end{tikzpicture} \subset \dynkin{D}{oooo},~ \textnormal{ and } \begin{tikzpicture}[baseline=0.5ex] 
\dynkin[name=first]{A}{o}
\node (current) at ($(first root 1)+(-0,3.4mm)$) {};
\dynkin[at=(current),name=second]{A}{*}
\node (currenta) at ($(first root 1)+(-3.4mm,0)$) {};
\dynkin[at=(currenta),name=third]{A}{o}
\node (currenta) at ($(third root 1)+(-0,3.4mm)$) {};
\dynkin[at=(currenta),name=fourth]{A}{*}
\begin{scope}[on background layer]
\draw[/Dynkin diagram/fold style]
($(first root 1)$) -- ($(third root 1)$);%
\draw[/Dynkin diagram/fold style]
($(second root 1)$) -- ($(fourth root 1)$);%
\end{scope}
\end{tikzpicture} \subset \dynkin{D}{oooo},$$ each of which represents a single isomorphism class which splits into 3 conjugacy classes, along with the embedding $$\begin{tikzpicture}[baseline=-0.3ex]
\dynkin[ply=2,rotate=90,fold radius = 2mm]{A}{o*o}\end{tikzpicture}\times\mathcal{T}^1_0 \subset \dynkin{D}{oo**},$$ which represents a single isomorphism class which splits into 2 conjugacy classes. 
For $n>4$ a single isomorphism class splits into at most 2 conjugacy classes. 

\begin{proof} Let $p \in \mathbb{P} \cup \{0\}$ and let $\Phi$ be an irreducible root system of classical type. We do not require that $p$ and $\Phi$ are fixed, rather we allow them to run over all possibilities.

\vspace{2mm}\noindent Our method of proof is as follows. For every abstract index $\mathcal{G}$ of type $\Phi$, we describe and implement an algorithm which we use to classify all $\mathcal{G}$-conjugacy classes of isotropic maximal $p$-embeddings of abstract indices $(\mathcal{H},\theta)$. We present our results in columns 1--7 of Tables \ref{A_n}, \ref{B_n}, \ref{C_n}, \ref{D_4} and \ref{D_n} respectively for the cases where $\Phi$ equals $A_n$, $B_n$, $C_n$, $D_4$, and $D_n$ for $n>4$.

\vspace{2mm}\noindent We first look at rows 1 and 2 (of Tables \ref{A_n}, \ref{B_n}, \ref{C_n}, \ref{D_4} and \ref{D_n}).

\vspace{2mm}\noindent By the argument in \cite[Prop.\ 2.15, proof]{S}, $\Delta$ is a (non-empty) almost primitive subsystem (of simple roots) of $\Phi$. As such, in columns $1$ and $2$ respectively, we use Table \ref{classlist} and Lemma \ref{isoautgps} to list all $W$-conjugacy classes of almost primitive subsystems $\Delta$ of $\Phi$ along with their associated groups $\Stab_{\Iso(\Phi)}(\Delta)$. The only case where there are multiple distinct $W$-conjugacy classes of $\Delta$ in $\Phi$ which are all $\Iso(\Phi)$-conjugate, is $A_{n-1} \subset D_n$ for $n \geq 4$ even. If $n=4$ (resp. $n>4$ is even) there are 3 (resp. 2) of such $W$-conjugacy classes. We add this information to the relevant rows.

\vspace{2mm}\noindent We now move on to column 3. Recall from \cite[Lem.\ 2.10]{S} that the property of $(\mathcal{G},\mathcal{H},\theta)$ being maximal (a.k.a. $\mathcal{H}$ being maximal in $\Phi$) depends only on the $W$-conjugacy class of $\Delta$ in $\Phi$ and the conjugacy class of $\Pi$ in $\Stab_{\Iso(\Phi)}(\Delta)$. As such, for each $\Delta$ in column $1$, our next task is to list -- in column $3$ -- all conjugacy classes of subgroups $\Pi$ of $\Stab_{\Iso(\Phi)}(\Delta)$ which satisfy this condition. These fall into three categories.

\vspace{2mm}\noindent Firstly, suppose that $\Psi=\langle \Delta \rangle$ is a maximal $p$-closed subsystem of maximal rank in $\Phi$. By inspection of Table \ref{classlist} this holds when $\Psi \subset \Phi$ is one of the following: $D_mB_{n-m} \subset B_n$ for $2 \leq m <n$ and $p \neq 2$, $D_n \subset B_n$ for all $p$, $B_mB_{n-m} \subset B_n$ for $1 \leq m \leq n/2$ and $p=2$, $C_mC_{n-m} \subset C_n$ for $1 \leq m \leq n/2$ and all $p$, $\smash{\widetilde{D_n}} \subset C_n$ for $p=2$, $D_mD_{n-m} \subset D_n$ for $1 \leq m \leq n/2$ and all $p$ (which includes $A_1^4 \subset D_4$). For these cases we observe that $V$ is trivial, since $\Psi$ has maximal rank in $\Phi$. So we simply list all conjugacy classes of subgroups of $\Stab_{\Iso(\Phi)}(\Delta)$. These groups are very small, so this is no trouble.

\vspace{2mm}\noindent Secondly, suppose that $\Psi$ is a maximal $p$-closed subsystem of $\Phi$ which is not of maximal rank in $\Phi$. Again by inspection of Table \ref{classlist}, this holds when $\Psi \subset \Phi$ is one of the following: $A_mA_{n-m-1} \subset A_n$ for $0 \leq m \leq (n-1)/2$ and for all $p$, $B_{n-1} \subset B_n$ for $p \neq 2$, $A_{n-1} \subset C_n$ for $p \neq 2$, $A_{n-1} \subset D_n$ for all $p$ (of which there are 3 conjugacy classes if $n=4$, 2 conjugacy classes if $n>4$ is even, and 1 conjugacy class if $n$ is odd). For each of these cases we list all conjugacy classes of subgroups of $\Stab_{\Iso(\Phi)}(\Delta)$, and then remove those ones which fix a non-trivial subspace of $V$. This is again easy, as the groups $\Stab_{\Iso(\Phi)}(\Delta)$ are small.

\vspace{2mm}\noindent Thirdly, suppose that $\Psi$ is not a maximal $p$-closed subsystem of $\Phi$. These are all of the remaining cases in Table \ref{classlist}. We consider them individually. \begin{itemize} \item Say $\Psi \subset \Phi$ is $A_2 \subset D_4$, for all $p$. Then $\Stab_{\Iso(\Phi)}(\Delta)$ is the dihedral group $\mathscr{D}i_{12}$. The only proper subsystems of $D_4$ which strictly contain a copy of $A_2$ are of type $A_3$. An element $\psi \in \Stab_{\Iso(\Phi)}(\Delta)$ stabilises one of the $A_3$ subsystems containing $\Delta$ if and only if the order of $\psi$ is not divisible by $3$. Since $\dim V =2$, any $\psi \in \Stab_{\Iso(\Phi)}(\Delta)$ whose order is divisible by 3 does not fix any non-trivial subspace of $V$. So the possibilities for $\Pi$ are those (conjugacy classes of) subgroups of $\mathscr{D}i_{12}$ whose order is divisible by 3; there are five of them.

\item Say $\Psi \subset \Phi$ is either $B_m^{n_m} \subset B_n$ for $mn_m=n$, $n_m > 2$ and $p=2$, or it is $C_m^{n_m} \subset C_n$ for $mn_m=n$, $n_m > 2$ and for all $p$. In both cases $\Stab_{\Iso(\Phi)}(\Delta)$ is the symmetric group $S_{n_m}$. Note that $V$ is trivial as the rank of $\Psi$ equals that of $\Phi$. So it suffices to find those conjugacy classes of subgroups $\Pi$ of $\Stab_{\Iso(\Phi)}(\Delta)$ for which there does not exist a $p$-closed $\Pi$-stable proper subsystem of $\Phi$ which strictly contains $\Psi$. We observe that the $W$-conjugacy classes of $p$-closed subsystems of $\Phi$ which contain $\Psi$ are in $1-1$ correspondence with partitions of $n_m$, so the aforementioned condition holds if and only if $\Pi$ is a primitive subgroup of $S_{n_m}$. [N.B. Of course we make no attempt to classify primitive permutation groups for arbitrary $n_m$; for an overview of the theory refer to \cite{Ca}.]

\item Say $\Psi \subset \Phi$ is $A_m^{n_m} \subset A_n$ for $n+1=(m+1)n_m$, $n_m > 2$ and for all $p$. Then $\Stab_{\Iso(\Phi)}(\Delta) \cong \Z_2 \times S_{n_m}$. A similar argument as for the previous case tells us that the conjugacy classes of subgroups $\Pi$ of $\Stab_{\Iso(\Phi)}(\Delta)$ for which there does not exist a $p$-closed $\Pi$-stable proper subsystem of $\Phi$ which strictly contains $\Psi$, are precisely those which project onto a primitive subgroup of its $S_{n_m}$ factor. Let $\Pi$ be such a subgroup of $\Stab_{\Iso(\Phi)}(\Delta)$. To ensure that $\mathcal{H}$ is maximal in $\Phi$, we need to check that $V^{\Pi}$ is trivial. Assume otherwise. Then the pointwise stabiliser $\smash{Z_W(V^{\Pi})}$ is a proper parabolic subgroup of $W$. 
Choose a corank 1 parabolic subgroup $W_1$ of $W$ which contains $\smash{Z_W(V^{\Pi})}$. Let $\Phi_1$ be the parabolic subsystem of $\Phi$ associated to $W_1$. Observe that $\Phi_1$ is a $p$-closed proper subsystem of $\Phi$, it contains $\Psi$, and it is $\Pi$-stable, 
hence it equals $\Psi$. But $\Phi_1$ has corank 1 in $\Phi$, which contradicts the fact that $n_m >2$. So indeed $\mathcal{H}$ is maximal in $\Phi$ for such choices of $\Pi$.

\item Finally, say $\Psi \subset \Phi$ is $D_m^{n_m} \subset D_n$ for $mn_m=n$, $n_m > 2$ and for all $p$. Then $\Stab_{\Iso(\Phi)}(\Delta) \cong \Z_2^{n_m} \rtimes S_{n_m}$. Note that $V$ is trivial as the rank of $\Psi$ equals that of $\Phi$. So, once again, the conjugacy classes of subgroups $\Pi$ of $\Stab_{\Iso(\Phi)}(\Delta)$ for which $\mathcal{H}$ is maximal in $\Phi$ are precisely those which project onto a primitive subgroup of its $S_{n_m}$ factor. \end{itemize} This exhausts all of the entries in Table \ref{classlist}. To refine our classification from $\mathcal{G}$-isomorphism to $\mathcal{G}$-conjugacy, we need to know when a given conjugacy class of subgroups $\Pi$ of $\Stab_{\Iso(\Phi)}(\Delta)$ splits into multiple distinct classes with respect to conjugation by an element of $\Stab_W(\Delta)$. By inspection of the groups $\Stab_{\Iso(\Phi)}(\Delta)$, this occurs only when $\Psi \subset \Phi$ is $D_m^{n_m} \subset D_n$ for $n \geq 4$ even, $m,n_m \geq 2$, and $mn_m=n$, in which case $\Stab_{\Iso(\Phi)}(\Delta) \cong \Z_2^{n_m} \rtimes S_{n_m}$ for $n>4$ and $\cong S_4$ for $n=2m=4$. Except for the $n=4$ case, we do not give an explicit description of the set of $\Stab_W(\Delta)$-conjugacy classes of such subgroups of $\Stab_{\Iso(\Phi)}(\Delta)$.

\vspace{2mm}\noindent This completes column 3 of our tables. We next look at columns 4 and 5.

\vspace{2mm}\noindent We use the classification of isomorphism classes of irreducible abstract indices found in \cite[Table II]{T2} to list -- in column 4 -- all conjugacy classes in $\Phi$ of isotropic (i.e. we exclude those for which $\Delta =\Delta_0$) abstract indices $\mathcal{H}$ that are associated to the pair $(\Delta,\Pi)$.

\vspace{2mm}\noindent There is a subtlety here in that certain abstract indices $\mathcal{H}$ may be oriented in multiple ways in $\Phi$ (i.e. for a given $(\Delta,\Pi)$ there may exist multiple choices for $\Delta_0$ which are conjugate by an element of $\Iso(\Phi)$, but not by an element of $W$). 
By careful consideration of the subsystem structure of the classical irreducible root systems, we find that this phenomenon occurs only in the following situation: $\Phi=D_n$ for $n$ even (including $n=4$), $\Delta = D_mD_{n-m}$ for $m>0$ even, and $$\mathcal{H} \cong \dynkin[scale=1,labels={,,d,,,(r-1)d,,,,rd},rotate=0]{D}{*...*o*......*o*...**o} \times \dynkin[scale=1,labels={,,d',,,(r'-1)d',,,,r'd'},rotate=0]{D}{*...*o*......*o*...**o};$$ where $d,d' \geq 2$ are both powers of 2, $m=rd$, and $n-m=r'd'$. This occurs because $\Phi=D_n$ for $n>4$ even (resp. $n=4$) contains 2 (resp. 3) $W$-conjugacy classes of subsystems of type $(A_{d-1})^r(A_{d'-1})^{r'}$, all of which are $\Iso(\Phi)$-conjugate. For such a $\mathcal{H}$, it turns out that the only possibility for $\mathcal{G}$ is the abstract index $$\dynkin[scale=1,labels={,,d,,,(r-1)d,,,,n},rotate=0]{D}{*...*o*......*o*...**o},$$ and moreover that $d=d'$. For $n>4$, since this abstract index $\mathcal{G}$ is not invariant under the non-trivial diagram automorphism of $D_n$, it contains only one of the two orientations of $\mathcal{H}$. 
 The situation is similar for $n=4$, but for 3 orientations instead of 2. In summary, this phenomenon does not result in distinct $\mathcal{G}$-conjugacy classes of $\mathcal{H}$ which are $\mathcal{G}$-isomorphic -- so we ignore it in our tables.

\vspace{2mm}\noindent Next, for each abstract index $\mathcal{H}$ in column $4$, we find the $W$-conjugacy class of $\Psi_0$ in $\Phi$ and put it in column 5. This information can be read straight off of $\mathcal{H}$.

\vspace{2mm}\noindent We now turn our attention to column 6. Recall from \cite[Lem.\ 2.10]{S} that the $W$-conjugacy class of $\Phi_a$ in $\Phi$ is independent of the choice of representative of the conjugacy class of $\mathcal{H}$ in $\Phi$. As such, for each abstract index $\mathcal{H}$ in column $4$, we compute the $W$-conjugacy class of $\Phi_a$ in $\Phi$ and put it in column $6$. We perform these computations in Appendices $A$, $B$, $C$ and $D$ for the cases where $\Phi=A_n$, $B_n$, $C_n$ and $D_n$ respectively. These computations are quite laborious, so we have written some computer code in Julia and Oscar \cite{OSCAR} to assist with and verify these calculations.

\vspace{2mm}\noindent Finally, we move on to column 7. We inspect all (isomorphism classes of) abstract indices $\mathcal{G}$ of type $\Phi$ (these are listed in \cite[Table II]{T2}). For each abstract index $\mathcal{H}$ in column $4$, we exclude all abstract indices of type $\Phi$ which violate at least one of the conditions $(A.2)$, $(A.3)$ and $(A.4)$ of Definition \ref{embedabstractindices}. We list all remaining abstract indices $\mathcal{G}$ of type $\Phi$ in column 7 (if there are none we leave it blank).

\vspace{2mm}\noindent Using $(A.2)$, if $\Pi \leq \Stab_W(\Delta)$ then we exclude those abstract indices $\mathcal{G}$ of outer type. Otherwise, we exclude those abstract indices $\mathcal{G}$ of inner type (this condition is vacuous when $\Phi$ is $B_n$ or $C_n$). If $\Phi=D_4$ then we must also distinguish between those $\Pi$ which are associated to an order 2 outer automorphism of $D_4$, and those which are associated to a triality automorphism.

\vspace{2mm}\noindent Using $(A.3)$, we exclude $\mathcal{G}$ unless $\Phi$ contains a $\Xi$-stable parabolic subsystem of type $\Phi_a$ which contains $\Phi_0$ (by ``of type $\Phi_a$" we mean that it lies in the $W$-conjugacy class of subsystems of $\Phi$ that we computed in column 6). [N.B. For the cases where there are multiple parabolic subsystems $\Phi_a$ of $\Phi$ which are $\Iso(\Phi)$-conjugate but not $W$-conjugate, we must take care to distinguish between them. In our tables this occurs precisely when $\mathcal{G}$ is of type ${}^1\hspace{-0.5mm}D_n$ for even $n\geq 4$, and $\smash{\Phi_a \cong A_{j-1}^{n/j}}$ for some even $j \geq 2$.]

\vspace{2mm}\noindent Next, we compute the $W$-conjugacy class of ${}_{in}\Phi_a$ in $\Phi$ (by \cite[Lem.\ 2.10]{S} these are independent of the choice of representative of the conjugacy class of $\mathcal{H}$ in $\Phi$). Using $(A.4)$, we exclude those abstract indices $\mathcal{G}$ where $\Phi_0$ does not contain a parabolic subsystem of type ${}_{in}\Phi_a$ (in most cases ${}_{in}\Phi_a$ is trivial, and so this restriction is vacuous). [N.B. To be precise we must do something slightly stronger. We must also compute the $W$-conjugacy class of $\Phi \cap \overline{E_a}$ in $\Phi$, and exclude those abstract indices $\mathcal{G}$ where $\Phi_a$ does not contain a parabolic subsystem of type ${}_{in}\Phi_a \big(\Phi \cap \overline{E_a}\big)$ such that the ${}_{in}\Phi_a$ component is contained in $\Phi_0$. It turns out that this slightly stronger condition is relevant in our tables only when $\Psi \subset \Phi$ is $A_{n-1} \subset D_n$ for some even $n \geq 4$, and $\mathcal{H}$ is either $\smash{\dynkin[scale=1, ply=2,rotate=0,fold radius = 2mm]{A}{oo...oo*oo...oo} \+ \cT_0^1}$ or $\smash{\dynkin[scale=1, ply=2,rotate=0,fold radius = 2mm]{A}{*o*...*o*o*...*o*} \+ \cT_0^1}$.] 

\vspace{2.5mm}\noindent This completes the justification of the entries in the first seven columns of Tables \ref{A_n}, \ref{B_n}, \ref{C_n}, \ref{D_4} and \ref{D_n}. For each pair of indices $(\mathcal{H},\mathcal{G})$ listed in columns $4$ and $7$, the embedding $\theta$ of $\mathcal{H}$ in $\mathcal{G}$ has been uniquely determined (up to $\mathcal{G}$-conjugacy) by the data present in columns $1-4$. The application of this procedure, along with the computations in the Appendix, completes the proof of Theorem \ref{classthmcomb}.
\end{proof}

\noindent Now recall the statement of \cite[Thm.\ 1(iii)]{S}: given a reductive $k$-group $G$, the index-conjugacy class of a maximal rank reductive subgroup $H$ of $G$ is uniquely determined by the conjugacy class of its embedding of indices. So Theorem \ref{classthm} follows immediately from combining \cite[Thm.\ 1(iii)]{S} with Theorem \ref{classthmcomb}.

\vspace{2mm}\noindent We now move on to the proof of Corollary \ref{classthmcor}.

\vspace{2mm}\noindent For $k$ a field, let $\Gamma$ denote the absolute Galois group of $k$ and $\mathrm{cd}(k)$ the cohomological dimension of $\Gamma$. Let $\R$ denote the real numbers. A \textit{$\mathfrak{p}$-adic field} refers to a non-Archimedean local field of characteristic zero.

\vspace{2mm}\noindent\underline{Proof of Corollary \ref{classthmcor}.}

\vspace{3mm}\noindent We restate Corollary \ref{classthmcor} in a slightly more precise form than in the introduction.

\begin{Corollary}\label{classthmcombprecise} Let $(\mathcal{G},\mathcal{H},\theta)$ be one of the embeddings of abstract indices that are listed in Tables \ref{A_n}, \ref{B_n}, \ref{C_n}, \ref{D_4} and \ref{D_n}. If there is a $\checkmark$ in column $8$ (resp. $9$, $10$) of the corresponding row then $(\mathcal{G},\mathcal{H},\theta)$ is $k$-admissible for some field $k$ with $\cd(k)=1$ (resp. $k=\R$, $k$ is $\mathfrak{p}$-adic). If there is a $\xmark$ then $(\mathcal{G},\mathcal{H},\theta)$ is not $k$-admissible for such a field.

\vspace{2mm}\noindent \textit{[In the cases where we list more than one value for $\Pi$ in a single row of one of the tables, we put a $\checkmark$ in the corresponding column if at least one of these values for $\Pi$ corresponds to a $k$-admissible $(\mathcal{G},\mathcal{H},\theta)$.]}
\begin{proof} We describe and implement an algorithm which we use to determine whether or not $(\mathcal{G},\mathcal{H},\theta)$ is $k$-admissible for some field $k$ with $\cd(k)=1$ (resp. $k=\R$, $k$ is $\mathfrak{p}$-adic). We continue to use the notation $\mathcal{G}=(E,\Lambda,\Lambda_0,\Xi)$ and $\mathcal{H}=(E,\Delta,\Delta_0,\Pi)$.

\vspace{2mm}\noindent \underline{\textit{Column $8$}}: $\cd(k) =1$.

\vspace{2mm}\noindent If $(\mathcal{G},\mathcal{H},\theta)$ is $k$-admissible for some field $k$ with $\cd(k)=1$, then by $\S 3.3.1$ of \cite{T2} it is quasisplit (i.e. both $\Lambda_0$ and $\Delta_0$ are empty). We claim that the converse holds also.

\vspace{2mm}\noindent Suppose $(\mathcal{G},\mathcal{H},\theta)$ is quasisplit. Choose any $p \in \mathbb{P} \cup \{0\}$ for which $(\mathcal{G},\mathcal{H},\theta)$ is a $p$-embedding of abstract indices. Let $K$ be an algebraically closed field of characteristic $p$, and consider the field of rational functions $K(t)$ in one variable over $K$. It is known that the absolute Galois group of $K(t)$ is a free profinite group of rank equal to the cardinality of $K$. For $p=0$ this result is due to Douady \cite[Thm.\ 3.4.8]{Sz}, for $p >0$ it was proved independently by Harbater \cite[Cor. 4.2(b)]{Ha} and Pop \cite[Cor. of Thm.\ A]{Po}. Consequently, every finite group may be realised as a quotient of the absolute Galois group of $K(t)$. 
Since $\mathcal{G}$ is quasisplit, it follows that $\mathcal{G}$ is $K(t)$-admissible. 
So we can apply \cite[Cor.\ 2]{S} which tells us that $(\mathcal{G},\mathcal{H},\theta)$ is $K(t)$-admissible. Any free profinite group has cohomological dimension 1. So the claim holds.

\vspace{2mm}\noindent As such we put a $\checkmark$ in column $8$ if $\Lambda_0$ and $\Delta_0$ are both empty, otherwise we put an $\xmark$.

\vspace{2mm}\noindent We cross-check the entries in column $8$ using \cite[Ch.\ 27, 28]{MT}.

\vspace{2mm}\noindent \underline{\textit{Column $9$}}: $k=\R$.

\vspace{2mm}\noindent Suppose $k=\R$. We first check if $(\mathcal{G},\mathcal{H},\theta)$ is \textit{not} a $0$-embedding of abstract indices. By inspection of Table \ref{classlist} this occurs precisely when the pair $(\Delta,\Lambda)$ is one of: $\smash{\widetilde{D_n}} \subset C_n$, $B_mB_{n-m} \subset B_n$ for $1 \leq m < n/2$, or $(B_l)^{n_l} \subset B_n$ for $ln_l=n$ and $n_l \geq 2$. If $(\mathcal{G},\mathcal{H},\theta)$ corresponds to one of these cases we put a $\xmark$ in column 9, since $\R$ has characteristic 0. Otherwise, we proceed as follows. 

\vspace{2mm}\noindent We next check if $\mathcal{G}$ is $k$-admissible, using \cite[Table II]{T2}. If not, we put a $\xmark$ in column $9$. If so, we proceed as follows. 

\vspace{2mm}\noindent If $(\mathcal{G},\mathcal{H},\theta)$ is independent then we apply \cite[Cor.\ 1]{S} and put a $\checkmark$ in column $9$. If $(\mathcal{G},\mathcal{H},\theta)$ is quasisplit but not independent then we apply \cite[Cor.\ 2]{S} and put a $\checkmark$ (resp. $\xmark$) in column $9$ if $\Pi=1$ or $\Z_2$ (resp. $|\Pi|>2$). If $(\mathcal{G},\mathcal{H},\theta)$ is neither quasisplit nor independent then we resort to using \cite[Table 3]{K}. Namely, we put a $\checkmark$ in column $9$ if the corresponding embedding exists in \cite[Table 3]{K}, otherwise we put a $\xmark$.

\vspace{2mm}\noindent [N.B. One could independently reproduce \cite[Table 3]{K} via the following method. Take a chain of reductive $\C$-groups $T \subset H \subset G$ where $T$ is a torus and $G$ is simple, such that the root system of $G$ (resp. $H$) with respect to $T$ is of type $\Lambda$ (resp. $\Delta$). Classify all conjugacy classes of holomorphic involutions of $G$ which stabilise both $H$ and $T$. This turns out to be equivalent to classifying $\R$-descents of the chain $T \subset H \subset G$.]

\vspace{2mm}\noindent \underline{\textit{Column $10$}}: $k$ is $\mathfrak{p}$-adic.

\vspace{2mm}\noindent As with column 9, we first check if $(\mathcal{G},\mathcal{H},\theta)$ is a $0$-embedding of abstract indices. If it is not, we put a $\xmark$ in column 10 since $\mathfrak{p}$-adic fields have characteristic 0. Otherwise, we proceed as follows. 

\vspace{2mm}\noindent We check if $\mathcal{G}$ is $k$-admissible for some $\mathfrak{p}$-adic field $k$, using \cite[Table II]{T2}. If not, we put a $\xmark$ in column $10$. If so, we proceed as follows.

\vspace{2mm}\noindent If $(\mathcal{G},\mathcal{H},\theta)$ is independent then we apply \cite[Cor.\ 1]{S} and put a $\checkmark$ in column $10$.

\vspace{2mm}\noindent Next, consider the cases where $(\mathcal{G},\mathcal{H},\theta)$ is quasisplit but not independent. By \cite[Cor.\ 2]{S} we put a $\checkmark$ in column $10$ if we can show there exists a Galois $\mathfrak{p}$-adic field extension with Galois group isomorphic to $\Pi$, if we can show such a field extension does not exist we put a $\xmark$. To complete this task, we need some understanding of which finite groups may be realised as quotients of the absolute Galois group of some $\mathfrak{p}$-adic field. 

\vspace{2mm}\noindent It is well-known that the Galois group of any finite $\mathfrak{p}$-adic field extension is solvable, so we can exclude all cases where $\Pi$ is not solvable. Moreover, the absolute Galois group of any $\mathfrak{p}$-adic field has a quotient isomorphic to $\smash{\hat{\Z}}$; the profinite completion of the integers (see for instance \cite[II.7]{GMS}). 
So we put a $\checkmark$ in column $10$ for all cases where $\Pi$ is cyclic. In particular, for any prime $q$, we note that the cyclic group $\Z_q$ is a primitive permutation group of degree $q$. 

\vspace{2mm}\noindent The finite groups $S_3$, $(\Z_2)^2$ and $\mathrm{Alt}_4$ each may be realised as the Galois group of a $\mathfrak{p}$-adic field extension. To see this, we give examples of such extensions using \cite{LMFDB}. The extension of $\mathbb{Q}_3$ with defining polynomial $x^6+3$ has Galois group $S_3$. 
The extension of $\mathbb{Q}_2$ with defining polynomial $x^4+2x^2+4x+2$ has Galois group $(\Z_2)^2$. 
The extension of $\mathbb{Q}_2$ with defining polynomial $$x^{12} + 6x^{11} + 22x^{10} + 56x^9 + 126x^8 + 240x^7 + 332x^6 - 18x^5 - 459x^4 - 394x^3 - 344x^2 + 138x + 423$$ has Galois group $\mathrm{Alt}_4$. 
As such we put a $\checkmark$ in column $10$ for the cases where $\Pi \cong S_3$, $(\Z_2)^2$ and $\mathrm{Alt}_4$.

\vspace{2mm}\noindent Next, we consider the cases where $\Delta_0$ is empty, $\Pi \cong \Z_2$ and $\Lambda_a=\Lambda_0 \cong A_1$. Consider the embedding of abstract indices $(\mathcal{L},\mathcal{H}_m,\theta)$; where $\mathcal{L}:=(E,\Lambda_a,\Lambda_0,\Xi)$ and $\mathcal{H}_m:=(E,\Delta_0,\Delta_0,\Pi)$. Recall that $\mathcal{G}$ is $k$-admissible for some  $\mathfrak{p}$-adic field $k$. So let $G$ be an absolutely simple $k$-group with index isomorphic to $\mathcal{G}$. Let $L$ be a Levi subgroup of a minimal parabolic subgroup of $G$, and let $T$ be any maximal torus of $L$. Observe that the embedding of indices of the pair $T \subset L$ is isomorphic to $(\mathcal{L},\mathcal{H}_m,\theta)$ (since any anisotropic absolutely simple $k$-group of type $A_1$ contains a $1$-dimensional anisotropic torus, for which there is only one possible Galois action). Then applying \cite[Thm.\ 3]{S} tells us that $(\mathcal{G},\mathcal{H},\theta)$ is $k$-admissible. So we put a $\checkmark$ in column $10$ for all such cases.

\vspace{2mm}\noindent For all remaining cases, we are unsure of whether or not $(\mathcal{G},\mathcal{H},\theta)$ is $k$-admissible for some $\mathfrak{p}$-adic field $k$. So we either put a $?$ symbol, or leave the corresponding entry empty, in column $10$.

\vspace{2mm}\noindent The application of this algorithm completes the proof of the corollary.
\end{proof}
\end{Corollary}

\vspace{2mm}\noindent\underline{Proof of Corollary \ref{classthmcor1}.}

\begin{proof} This is by inspection of the restrictions on the parameter $d$ which appear in column 7 of Tables \ref{A_n}, \ref{B_n}, \ref{C_n}, \ref{D_4} and \ref{D_n}. We interpret these restrictions group-theoretically using the discussions alongside \cite[Table 2]{T2}, as well as \cite[Ch.\ 24]{Mi}.
\end{proof}

\vspace{2mm}\noindent We now move on to the proof of Theorem \ref{classthmcor2}. All logarithms are taken to be of base 2.

\vspace{2mm}\noindent\underline{Proof of Theorem \ref{classthmcor2}.}

\begin{proof} Let $G$ be an absolutely simple $k$-group.

\vspace{2mm}\noindent We first note that, if our field $k$ is perfect, then a maximal connected subgroup $H$ of $G$ is either parabolic or reductive by \cite[Cor.\ 3.7]{BT} (Borel-Tits). Over an arbitrary field $k$ this fact no longer holds, however it remains true if in addition we require that $H$ has maximal rank in $G$. It follows from \cite[Thm.\ 20.9(iii)]{Bo} that the number of $G(k)$-conjugacy classes of maximal parabolic subgroups of $G$ is bounded above by the rank of $G$. So we need to find an upper bound for the number of index-conjugacy classes of isotropic maximal reductive subgroups of maximal rank in $G$.

\vspace{2mm}\noindent Let $\Phi$ be an irreducible root system. Henceforth we use the same combinatorial notation as in the proof of Theorem \ref{classthmcomb}. Namely, $W$ is the Weyl group of $\Phi$, $\mathcal{G}$ is an irreducible abstract index of type $\Phi$, $p$ is either a prime or 0, $(\mathcal{G},\mathcal{H},\theta)$ is a $p$-embedding of abstract indices, and $\mathcal{H}=(E,\Delta,\Delta_0,\Pi)$. Let $n$ denote the rank of $\Phi$.

\vspace{2mm}\noindent We define $N(\Phi)$ to be the supremum of the number of $\mathcal{G}$-conjugacy classes of isotropic maximal $p$-embeddings of abstract indices, where the supremum is taken as $\mathcal{G}$ runs over all (isomorphism classes of) irreducible abstract indices of type $\Phi$ and $p$ runs over all prime numbers along with 0. Applying \cite[Thm.\ 1]{S} tells us that the number of index-conjugacy classes of isotropic maximal reductive subgroups of maximal rank in $G$ is bounded above by $N(\Phi(G))$; where $\Phi(G)$ denotes the root system of $G$. So, for each irreducible root system of classical type $\Phi$, it suffices to find an upper bound for $N(\Phi)$ as the rank $n$ of $\Phi$ grows large.

\vspace{2mm}\noindent Now suppose $\Phi$ is of classical type, and take $p$ and $\mathcal{G}$ to be fixed. By definition $N(\Phi)$ is bounded above by the number of possibilities for the triple $(\Delta,\Delta_0,\Pi)$ up to $W$-conjugacy. 
We use Theorem \ref{classthm} and the corresponding tables in Appendices \ref{An}--\ref{Dn} to find an appropriate upper bound for the number of such triples $(\Delta,\Delta_0,\Pi)$. Note that $\Delta$ must be non-empty, as otherwise $\mathcal{H}$ would be anisotropic. 

\vspace{2mm}\noindent A close inspection of the tables in Appendices \ref{An}--\ref{Dn} reveals that, for fixed $\Pi$, the $W$-conjugacy class of the pair $(\Delta,\Delta_0)$ depends on at most 3 free parameters, each of which can take at most $n$ possible values. So, for fixed $\Pi$, the number of possible $W$-conjugacy classes for the pair $(\Delta,\Delta_0)$ is of the order $O(n^3)$. In particular, it grows polynomially in $n$.

\vspace{2mm}\noindent It turns out that the number of possibilities for $\Pi$ up to $W$-conjugacy does not grow polynomially in $n$. Again by inspection of the tables in Appendices \ref{An}--\ref{Dn}, the worst case scenario is realised when $\Phi=C_n$, $\Delta =(C_1)^n$ and $\Delta_0=\varnothing$ (for any $p$). In this case, the possibilities for $\Pi$ are in 1--1 correspondence with conjugacy classes of primitive subgroups of the symmetric group $S_n$. By \cite[Thm.\ I]{PS}, the number of conjugacy classes of primitive subgroups of $S_n$ is bounded above by $n^{c\log n}$, where $c$ is an absolute constant. 
In summary, we have shown that $N(\Phi)=O(n^{c\log n})$. This proves the first assertion.

\vspace{2mm}\noindent For the second assertion, suppose $k$ is separably closed. All reductive groups over $k$ are split, hence $\Delta_0$ is empty and $\Pi$ is trivial. So $N(\Phi)$ is bounded above by the number of $W$-conjugacy classes of almost primitive subsystems $\Delta$ of $\Phi$, which -- a brief inspection of Table \ref{classlist} tells us -- grows as $O(n)$. 
This completes the proof of the theorem.
\end{proof}

\renewcommand{\thesubsection}{\Alph{subsection}}
\titleformat{\subsection}{\normalfont\large\bfseries}{Appendix \thesubsection}{1em}{}
\section*{Appendix}\label{appendix}
\noindent In the appendix we complete the statement and proof of Theorem \ref{classthm}. We perform the remaining calculations needed to construct Tables \ref{A_n}, \ref{B_n}, \ref{C_n}, \ref{D_4} and \ref{D_n}, then we present the tables. Explicitly, for each abstract index $\mathcal{H}$ listed in the fourth column of our tables, we compute the $W$-conjugacy class of $\Phi_a$ in $\Phi$ (the entry in the sixth column of our tables). 

\vspace{2mm}\noindent We consider the case where $\Phi$ is $A_n$, $B_n$, $C_n$, and $D_n$ in Appendices $A$, $B$, $C$, and $D$ respectively. 

\vspace{2mm}\noindent We use the setup and notation of $\S \ref{embedabstractindices-section}$. Let $p$ either be a prime number or $0$. Let $(\mathcal{G}, \mathcal{H}, \theta)$ be a $p$-embedding of abstract indices, say $\mathcal{G}=(F,\Lambda,\Lambda_0, \Xi)$ and $\mathcal{H}=(E,\Delta,\Delta_0,\Pi)$. We abuse notation and identify $\mathcal{H}$ with its image under $\theta$ (i.e. we identify $E=F$). 

\vspace{2mm}\noindent We denote $\Phi:=\langle \Lambda \rangle$, $\Phi_0:=\langle \Lambda_0 \rangle$, $W$ is the Weyl group of $\Phi$, $\Iso(\Phi)$ is the isometry group of $\Phi$, $\Psi:=\langle \Delta \rangle$, $\Psi_0:=\langle \Delta_0 \rangle$, $E_s$ is the subspace of $E^{\Pi}$ which is perpendicular to $\Delta_0$, $\Phi_s := \Phi \cap E_s$, $E_a:=(E_s)^{\perp}$, $\Phi_a := \Phi \cap E_a$, and $E_{\Delta}$ (resp. $E_{\Delta_0}$) is the span of $\Delta$ (resp. $\Delta_0$) in $E$. By assumption $\Phi$ is an irreducible root system of classical type, and it spans $E$. 

\vspace{2mm}\noindent We introduce some additional notation. Let $\alpha_0$ be the highest root of $\Psi$. 
If $\Phi$ is not simply-laced then we use a $\widetilde{\textcolor{white}{m}}$ to denote a subsystem of $\Phi$ that consists only of short roots. We follow the convention for root systems that $A_0=A_{-1}=B_0=C_0=D_0=D_1=\varnothing$, $A_1=B_1=C_1$, $D_2=A_1^2$, $B_2=C_2$, and $D_3=A_3$.

\addcontentsline{toc}{section}{Appendix}
\subsection{\texorpdfstring{\hspace*{-0.45cm}. $\Phi=A_n =\dynkin[labels={\alpha_1, \alpha_2,\alpha_{n-1},\alpha_n},scale=2.2,
edge length=.5cm, root radius = 0.03cm] A{**...**}$}{An}}\label{An}

Let $\Phi= A_n$ for some integer $n \geq 1$. We use the following standard parametrisation in $\R^{n+1}$ of $\Phi$. Let $\alpha_i = e_i-e_{i+1}$ for all $i\in\{1,\ldots,n\}$, where $e_i$ denotes the $i$'th standard basis vector.\\

Following \cite[Table II]{T2}, the isotropic indices of type $\Phi$ are given as follows:\\
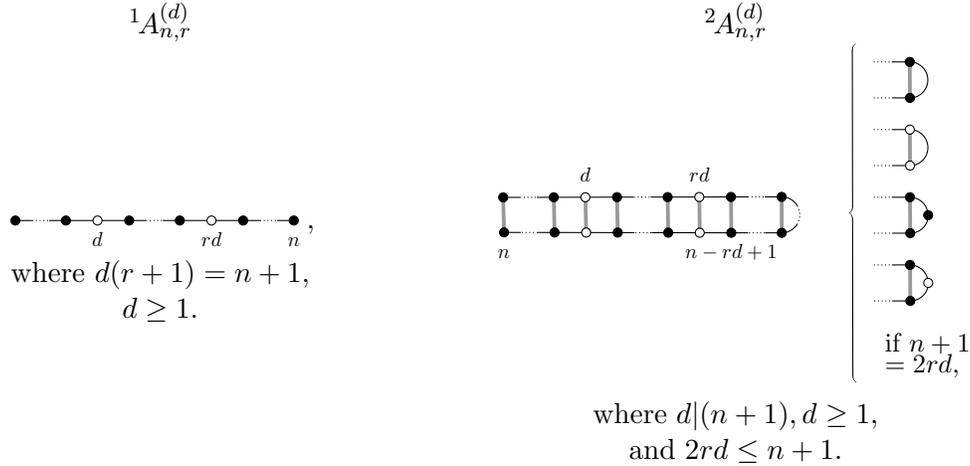
\begin{center}
\begin{tabular}{cp{1.5cm}c}
$^{1\!}A^{(d)}_{n,r}$ & &$^{2\!}A^{(d)}_{n,r}$ \\
 \multirow{2}{*}{\makecell{ $\dynkin[scale=1.2,labels={,,d,,,rd,,n}]{A}{*...*o*......*o*...*},$\\ where $d(r+1) = n+1$,\\ $d\geq 1$.}} & &\makecell{%
\begin{tikzpicture}[scale=1.2] 
\node (up) at (0,2) {};
\dynkin[at=(up),name=upper,labels={,,d,,,rd,,,,{n-rd +1},,,,,,n}, ply=2,fold radius=2mm]{A}{*.*o*.*o*.*.*.*o*.*o*.*};
\node (up1) at ($(up)+(4.5,1.5)$) {};
\node (up2) at ($(up1)+(0,-0.75)$) {};
\node (up3) at ($(up2)+(0,-0.75)$) {};
\node (up4) at ($(up3)+(0,-0.75)$) {};
\dynkin[at=(up1),name=upper_right, labels={,},ply=2,fold radius=2mm]{A}{.**.};
\draw[densely dotted] ($(upper_right root 1)+(-0.4,0)$) -- ($(upper_right root 1)+(-0.2,0)$);
\draw ($(upper_right root 1)+(-0.2,0)$) -- ($(upper_right root 1)+(-0.05,0)$);
\draw[densely dotted] ($(upper_right root 2)+(-0.4,0)$) -- ($(upper_right root 2)+(-0.2,0)$);
\draw ($(upper_right root 2)+(-0.2,0)$) -- ($(upper_right root 2)+(-0.05,0)$);
\dynkin[at=(up2),name=middle1_right, labels={,}, ply=2,fold radius=2mm]{A}{.oo.};
\draw[densely dotted] ($(middle1_right root 1)+(-0.4,0)$) -- ($(middle1_right root 1)+(-0.2,0)$);
\draw ($(middle1_right root 1)+(-0.2,0)$) -- ($(middle1_right root 1)+(-0.05,0)$);
\draw[densely dotted] ($(middle1_right root 2)+(-0.4,0)$) -- ($(middle1_right root 2)+(-0.2,0)$);
\draw ($(middle1_right root 2)+(-0.2,0)$) -- ($(middle1_right root 2)+(-0.05,0)$);
\dynkin[at=(up3),name=middle_right, labels={,,}, ply=2,fold radius=2mm]{A}{.***.};
\draw[densely dotted] ($(middle_right root 1)+(-0.4,0)$) -- ($(middle_right root 1)+(-0.2,0)$);
\draw ($(middle_right root 1)+(-0.2,0)$) -- ($(middle_right root 1)+(-0.05,0)$);
\draw[densely dotted] ($(middle_right root 3)+(-0.4,0)$) -- ($(middle_right root 3)+(-0.2,0)$);
\draw ($(middle_right root 3)+(-0.2,0)$) -- ($(middle_right root 3)+(-0.05,0)$);
\dynkin[at=(up4),name=lower_right, labels={,,},ply=2,fold radius=2mm]{A}{.*o*.};
\draw[densely dotted] ($(lower_right root 1)+(-0.4,0)$) -- ($(lower_right root 1)+(-0.2,0)$);
\draw ($(lower_right root 1)+(-0.2,0)$) -- ($(lower_right root 1)+(-0.05,0)$);
\draw[densely dotted] ($(lower_right root 3)+(-0.4,0)$) -- ($(lower_right root 3)+(-0.2,0)$);
\draw ($(lower_right root 3)+(-0.2,0)$) -- ($(lower_right root 3)+(-0.05,0)$);
\draw [decorate, decoration = {calligraphic brace}] ($(lower_right root 3)+(-6mm,-9mm)$) --  ($(upper_right root 1)+(-6mm,2mm)$);
\node[align=left] at ($(up4) + (2mm,-8mm)$) {\small{if $n+1$}\\ \small{$ = 2rd$,}};
\end{tikzpicture}}\\
 & & \makecell{where $d|(n+1), d\geq 1,$\\ and $2rd\leq n+1$.}
\end{tabular}
\end{center}

We follow the convention that anisotropic indices of type ${}^1\hspace{-0.5mm}A_n$ (resp. ${}^2\hspace{-0.5mm}A_n$) satisfy $r=0$ and $d=n+1$ (resp. $r=0$ and $d$ may be any positive integer). 

Recall from Table \ref{classlist} that the $W$-conjugacy classes of almost primitive subsystems $\Delta$ of $\Phi$ are as follows:
\begin{enumerate}
\item $A_m^{n_m}$ with $(m+1)n_m = n+1$ and $n_m\geq 2$,
\item $A_m\+A_{n-m-1}$ for $0\leq m< (n-1)/2$.
\end{enumerate}

In the fourth column of Table \ref{A_n}, we list all abstract indices $\cH$ of type $\Delta$ that are maximal in $\Phi$.\\

We will first consider the case where $\Delta$ is of type $A_m^{n_m}$ with $(m+1)n_m = n+1$ and $n_m\geq 2$. 
By Table \ref{classlist}, we have $\Stab_{\Iso(\Phi)}(\Delta)=\Z_2\+S_{n_m}$. Let $\rho:\Z_2\+ S_{n_m} \to S_{n_m}$ denote the natural projection.
Since we are interested in maximal rank subgroups, we only need to consider indices of form $\bigtimes_{j=1}^{n_m} {}^{1\!}A_{m,r}^{(d)}$ where $\Pi \subseteq W$ and $\rho(\Pi)$ is a primitive subgroup of $S_{n_m}$, or $\bigtimes_{j=1}^{n_m} {}^{2\!}A_{m,r}^{(d)}$ where $\Pi \nsubseteq W$ and $\rho(\Pi)$ is a primitive subgroup of $S_{n_m}$.\\

Our goal is to compute the $W$-conjugacy class of $\Phi_a$ in $\Phi$ for each abstract index $\cH$ of type $\Delta$. We proceed in several steps:
\begin{enumerate}
    \item\label{E_0} Find the subspace $E_{\Delta_0}$ of $E$, i.e. the span of $\Delta_0$. 
    \item\label{E_0'} Compute the perp space $E_{\Delta_0}^{\perp}$ of $E_{\Delta_0}$ in $E$.
    \item\label{E_s} Compute the fixed point subspace $E_s := (E_{\Delta_0}^{\perp})^{\Pi}$. Note that $\dim E_s$ equals the number of $\Pi$-orbits of white vertices in $\cH$. 
    \item\label{E_a} Compute $E_a:=E_s^{\perp}$. 
    \item\label{Phi_a} Compute $\Phi_a:=\Phi\cap E_a$. We have $\rk(\Phi_a) \leq \dim E_a$ with equality if the index of $G$ is of inner type. In general, by counting dimensions we see that 
    $$\rk (\Phi_a) = \dim E_a-\sum_{\cO}(|\cO|-1),$$ where $\cO$ runs over all $\Xi$-orbits of $\cG$ which are contained in $\Phi_a$. 
\end{enumerate}
We will usually follow this approach. Since all classical groups have a substantial part of their root system that is of type $A_n$, we use this approach as a template for the other Dynkin types.\\

Let $E'$ be the $n+1$-dimensional $\R$-vector space  spanned by $\{e_1,...,e_{n+1}\}$. To streamline our computations, we embed $\Phi=A_n$ in the ambient space $E'$ and let $E=\spn{e_i-e_{i+1}\mid 1\leq i\leq n}\subseteq E'$.\\

First, consider the case where $\Pi \subseteq W$ and $\cH=\bigtimes_{j=1}^{n_m} {}^{1\!}A_{m,r}^{(d)}$. That is, we get the following index:\\
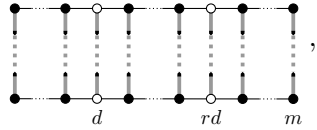
\begin{center}
    \begin{tikzpicture}[scale=1.2,baseline=-5mm] 
\dynkin[name=upper]{A}{*.*o*.*o*.*}
\node (current) at ($(upper root 1)+(0,-10mm)$) {};
\dynkin[at=(current),name=middle,labels={,,d,,,rd,,m}]{A}{*.*o*.*o*.*}
\begin{scope}[on background layer]
\foreach \x in {1,...,8} 
\draw[/Dynkin diagram/fold style]
($(upper root \x)$) -- ($(upper root \x)+(0,-3mm)$);
\foreach \x in {1,...,8} 
\draw[/Dynkin diagram/fold style]
($(middle root \x)+(0,3mm)$) -- ($(middle root \x)$);
\foreach \x in {1,...,8} 
\draw[dotted,black!40,line width=.5mm]
($(upper root \x)+(0,-3mm)$) -- ($(middle
root \x)+(0,3mm)$);
\end{scope}
\end{tikzpicture},\\ where $d(r+1) = m+1$.
\end{center}
\begin{enumerate}
    \item We see that $\Delta_0=\left\{\alpha_{i+j(m+1)}\left| \begin{array}{l}
    1\leq i \leq m, d\nmid i \\
    0\leq j \leq n_m-1
  \end{array}\right.\right\} $, so 
  \begin{align*}
    E'_{\Delta_0}&=\spn[\R]{e_{i+j(m+1)}-e_{i+1+j(m+1)}\left| \begin{array}{l}
    1\leq i \leq m, d\nmid i \\
    0\leq j \leq n_m-1
  \end{array}\right.}\\
  &=\spn[\R]{e_{i}-e_{i+1} \mid
    1\leq i \leq n, d\nmid i }.
  \end{align*}
\item Let 
\begin{align}\label{v_id}
v_{i,j}^{(d)}:= \sum_{l=1}^d e_{(i-1)d+l+j(m+1)}.
\end{align}
We compute that
$$(E'_{\Delta_0})^{\perp}=\spn[\R]{v_{i,j}^{(d)} \left| \begin{array}{l}
    1\leq i \leq r+1,\\
    0\leq j \leq n_m-1
  \end{array}\right.}.$$
  \item The action of $\Pi$ fixes the vectors 
  \begin{align}\label{beta_i}
      \beta_i := \sum_{j=0}^{n_m-1} \alpha_{i+j(m+1)} = \sum_{j=0}^{n_m-1} e_{i+j(m+1)}- e_{i+1+j(m+1)}
      \end{align}
      for all $i\in\{1,\ldots,m\},$
  hence
  $$E^{\Pi}= \spn[\R]{\sum_{j=0}^{n_m-1} e_{i+j(m+1)}- e_{i+1+j(m+1)} \Bigg|  1\leq i\leq m }.$$ Note that 
  \begin{align*}
\sum_{i=1}^{d} i\beta_{(l_1-1)d+i} &= \sum_{j=0}^{n_m-1} (v_{l_1,j}^{(d)}- de_{l_1d+1+j(m+1)})\\
\sum_{i=1}^{(l_2-l_1-1)d-1} d\beta_{l_1d+i} & = \sum_{j=0}^{n_m-1} (de_{l_1d+1+j(m+1)} - de_{(l_2-1)d+j(m+1)}) \\
\sum_{i=0}^{d-1} (d-i)\beta_{(l_2-1)d+i}&= \sum_{j=0}^{n_m-1} (-v_{l_2,j}^{(d)}- de_{(l_2-1)d+j(m+1)}),
\end{align*}
  so
  $$\sum_{i=1}^{d} i\beta_{(l_1-1)d+i} +\sum_{i=1}^{(l_2-l_1-1)d-1} d\beta_{l_1d+i}+\sum_{i=0}^{d-1} (d-i)\beta_{(l_2-1)d+i}= \sum_{j=0}^{n_m-1} (v_{l_1,j}^{(d)}-v_{l_2,j}^{(d)})\in E_s,$$
  where $1\leq l_1<l_2\leq r+1$.
  In particular, since $\dim E_s = r,$ and $E_s = E^{\Pi}\cap (E'_{\Delta_0})^{\perp}$, we have
  $$E_s = \spn[\R]{\sum_{j=0}^{n_m-1} (v_{i,j}^{(d)}-v_{i+1,j}^{(d)})\Bigg\vert 1\leq i \leq r }.$$
  \item Let
\begin{align*}
    B_1 & := \left\{e_{i+j(m+1)}-e_{i+1+j(m+1)}\middle| \begin{array}{l}
    1\leq i \leq m, d \nmid i \\
    0 \leq j \leq n_m-1
\end{array}\right\}\\
B_2 & := \left\{e_{id+1+ (j-1)(m+1)}-e_{id+1+j(m+1)}\middle| \begin{array}{l}
    0\leq i \leq r \\
    1 \leq j \leq n_m-1
\end{array}\right\}\\
B & := B_1\cup B_2.
\end{align*}
Then $B$ is a basis of $E_a$, since its elements are linearly independent and $$\dim E_a = n-\dim E_s= n - r = n_m(m-r)+(r+1)(n_m-1) = |B|.$$
\item Finally, we compute $\Phi_a := E_a\cap \Phi$.\\
As $B_1, B_2 \subseteq \Phi^+,$ we can take the simple roots of $\Phi_a$ to be given by $B_1\cup B_2$. In particular, $\Phi_a = A_{n_md-1}^{r+1}$. 
We can see this as follows:\\
Each set $B_{(1,j)}:=\{e_{i+jd}-e_{i+1+jd}\mid 1\leq i \leq d-1\}\subseteq B_1$ for $0\leq j\leq (n+1){d}-1$ gives rise to a root system of type $A_{d-1}$. We have $(n+1)/d = n_m(r+1)$ of such sets. Now, because of the roots in $B_2$ the sets $B_{(1,j)}$ together with the sets $B_{(1,j+l(r+1))}$, for $1 \leq l\leq n_m-1$, form an irreducible root system for each $0\leq j\leq r$. Therefore, we get $r+1$ irreducible root systems of type $A_q$, where $q := n_m(d-1)+n_m-1 = dn_m-1$.
\end{enumerate}
Next, let $\Pi \nsubseteq W$ such that $\rho(\Pi)$ is a primitive subgroup of $S_{n_m}$ and $\cH=\bigtimes_{j=1}^{n_m} {}^{2\!}A_{m,r}^{(d)}$. That is, we get the following index:
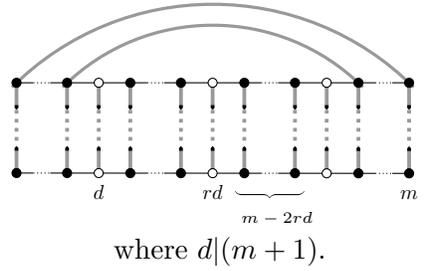
\begin{center}
    \begin{tikzpicture}[scale=1.2,baseline=-5mm] 
\dynkin[name=upper]{A}{*.*o*.*o*.*o*.*}
\node (current) at ($(upper root 1)+(0,-10mm)$) {};
\dynkin[at=(current),name=middle,labels={,,d,,,rd, , , , ,m}]{A}{*.*o*.*o*.*o*.*}
\begin{scope}[on background layer]
\foreach \x/\y in {1/11, 2/10}
\draw[very thick,black!40] ($(upper root \x)+(-0.5mm,-0.5mm)$) to [out=45, in=135] ($(upper root \y)+(0.5mm,-0.5mm)$);

\foreach \x in {1,...,11} 
\draw[/Dynkin diagram/fold style]
($(upper root \x)$) -- ($(upper root \x)+(0,-3mm)$);
\foreach \x in {1,...,11} 
\draw[/Dynkin diagram/fold style]
($(middle root \x)+(0,3mm)$) -- ($(middle root \x)$);
\foreach \x in {1,...,11} 
\draw[dotted,black!40,line width=.5mm]
($(upper root \x)+(0,-3mm)$) -- ($(middle root \x)+(0,3mm)$);
\end{scope}
\draw [decorate,
	decoration = {calligraphic brace}] ($(middle root 8)+(1mm,-2mm)$) --  ($(middle root 7)+(-1mm,-2mm)$);
 \node at ($(middle root 8)+(-2mm,-5mm)$) {\tiny $m-2rd$};
\end{tikzpicture},\\
where $d|(m+1)$.
\end{center}
We will again follow steps \ref{E_0} - \ref{Phi_a} above.
\begin{enumerate}
    \item In this case, 
    \begin{align*}
    E'_{\Delta_0} = &\spn[\R]{e_{i+j(m+1)}-e_{i+1+j(m+1)}\left| \begin{array}{l}
    1\leq i \leq rd, d\nmid i \\
    0\leq j \leq n_m-1
  \end{array}\right.}\\
  \oplus &\spn[\R]{e_{i+j(m+1)}-e_{i+1+j(m+1)}\middle| \begin{array}{l}
    rd+1\leq i \leq m-rd\\
    0\leq j \leq n_m-1
  \end{array}}\\
  \oplus &\spn[\R]{e_{i+j(m+1)}-e_{i+1+j(m+1)}\left| \begin{array}{l}
    m-rd+1\leq i \leq m, d\nmid i \\
    0\leq j \leq n_m-1
  \end{array}\right.}
\end{align*}
\item We get \begin{align*}
    (E'_{\Delta_0})^{\perp} = &\spn[\R]{v_{i,j}^{(d)} \left| \begin{array}{l}
    1\leq i \leq r,\\
    0\leq j \leq n_m-1
  \end{array}\right.} \oplus \spn[\R]{v_{i,j}^{(d)} \left| \begin{array}{l}
    \frac{m+1}{d}-r+1\leq i \leq \frac{m+1}{d},\\
    0\leq j \leq n_m-1
  \end{array}\right.}\\
  \oplus & \spn[\R]{\sum_{i=rd+1}^{m-rd+1}e_{i+j(m+1)}\bigg|  0\leq j \leq n_m-1 },
  \end{align*}
  where we use the same notation as in \eqref{v_id}.
  \item Recall from Step \ref{E_s} that $\dim E_s = r$. We have 
  $$ E^{\Pi} = \spn[\R]{\beta_i+\beta_{m-i+1}\middle| 1\leq i\leq \lfloor\frac{m}{2}\rfloor} \oplus \spn[\R]{\beta_{\frac{m+1}{2}}\middle| \text{if }\frac{m+1}{2}\in\N},$$ using the notation in \eqref{beta_i}. In particular,
\begin{align*}   
   &\sum_{j=id}^{\lfloor \frac{m+1}{2}\rfloor} d(\beta_{j}+\beta_{m-j+1}) +\sum_{j=1}^{d-1} (d-j)(\beta_{id-j}+\beta_{m-id+j+1})\\
   =&\sum_{j=0}^{n_m-1} (-v_{i,j}^{(d)}+v_{\frac{m+1}{d}-i+1,j}^{(d)})\in E_s = (E'_{\Delta_0})^{\perp}\cap E^{\Pi} ,   
\end{align*}
  where $1\leq i \leq r$. Therefore,
  $$ \spn[\R]{\sum_{j=0}^{n_m-1} (v_{i,j}^{(d)}-v_{\frac{m+1}{d}-i+1,j}^{(d)})\bigg|1\leq i\leq r }= E_s.$$
  \item Let \begin{align*}
    B_{1,1} & := \left\{e_{i+j(m+1)}-e_{i+1+j(m+1)}\middle| \begin{array}{l}
    1\leq i \leq rd, d \nmid i \\
    0 \leq j \leq n_m-1
\end{array}\right\}\\
B_{1,2} & := \left\{e_{m-i+1+j(m+1)}-e_{m-i+2+j(m+1)}\middle| \begin{array}{l}
    1\leq i \leq rd, d \nmid i \\
    0 \leq j \leq n_m-1
\end{array}\right\}\\
B_{2,1} & := \left\{e_{id+1+ (j-1)(m+1)}-e_{id+1+j(m+1)}\middle| \begin{array}{l}
    0\leq i \leq r-1 \\
    1 \leq j \leq n_m-1
\end{array}\right\}\\
B_{2,2} & := \left\{e_{m-id+(j-1)(m+1)}-e_{m-id+j(m+1)}\middle| \begin{array}{l}
    0\leq i \leq r-1 \\
    1 \leq j \leq n_m-1
\end{array}\right\}\\
B_3 & := \left\{e_{i+j(m+1)} \middle|\begin{array}{l}
    rd +1 \leq i \leq m-rd+1 \\
    0 \leq j \leq n_m-1
\end{array}\right\}\\
B_4 & := \left\{\sum_{j=0}^{n_m-1} (v_{i,j}^{(d)}+v_{\frac{m+1}{d}-i+1,j}^{(d)})\middle|1\leq i\leq r \right\}\\
B & := \bigcup_{i,j=1}^2 B_{i,j}\cup\bigcup_{i=3}^4 B_{i}.
\end{align*}
We have that $B\subseteq E_s^{\perp}=E'_a$, where we take the perp space in $E'$. As the elements in $B$ are all linearly independent, and $|B|= n+1-r =\dim E'_a$, $B$ is a basis of $E'_a$.
\item Note that $B_{i,j}\subseteq \Phi$ for $i,j \in \{1,2\}$ but $B_3,B_4\nsubseteq \Phi$. Also, note that $\Phi$ only contains elements of the form $e_i-e_j$ for some $i,j \in \{1,...,n+1\}$, therefore 
$$\Phi_a := \Phi\cap E'_a = \scaleleftright[1.75ex]{<}
{\left.\begin{array}{l}
B_{1,1}, B_{1,2},B_{2,1}, B_{2,2},\\
e_{i+j(m+1)}-e_{i+1+j(m+1)},\\
e_{i+j(m+1)}-e_{i+(j+1)(m+1)}
\end{array}\right| \begin{array}{l}
    rd+1\leq i \leq m-rd \\
    0 \leq j \leq n_m-1
\end{array}}{>}.$$
The sets $B_{i,1}$ together with $B_{i,2}$, $i=1,2$ give rise to subsystems of type $A_{dn_m-1}^r$, as in the previous case. With the same argument, one can see that $B_3$ and $B_4$ give rise to a subsystem of type $A_{n_m(m-2rd+1)-1}$. Together, we get $\Phi_a =A_{dn_m-1}^r\+A_{n_m(m-2rd+1)-1}$.
\end{enumerate}

Finally, consider subsystems of type $A_m\+A_{n-m-1}$ for $0\leq m< (n-1)/2$. In this case, $\Stab_{\Iso(\Phi)}(\Delta)=\Z_2$ by Table \ref{classlist}. We may ignore the case $\Pi = 1$, as any associated index $\mathcal{H}$ is not maximal in $\Phi$ (i.e. it wouldn't give rise to a maximal connected subgroup of a simple $k$-group $G$ with root system $\Phi$). So we need only consider the case $\Pi = \Z_2$. 
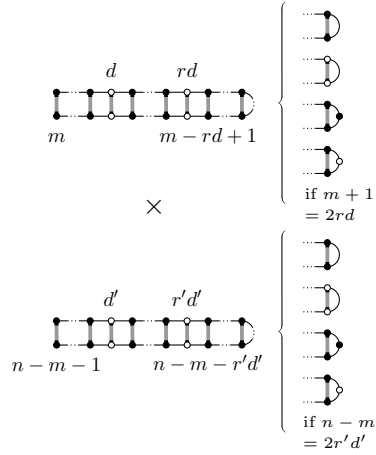
\begin{center}
\begin{tikzpicture}[scale=0.8,baseline=-5mm] 
\node (up) at (0,2) {};
\dynkin[at=(up),name=upper,labels={,,d,,,rd,,,,{m-rd +1},,,,,,m}, ply=2,fold radius=2mm]{A}{*.*o*.*o*.*.*.*o*.*o*.*};
\node (up1) at ($(up)+(4.5,1.3)$) {};
\node (up2) at ($(up1)+(0,-0.75)$) {};
\node (up3) at ($(up2)+(0,-0.75)$) {};
\node (up4) at ($(up3)+(0,-0.75)$) {};
\dynkin[at=(up1),name=upper_right, labels={,},ply=2,fold radius=2mm]{A}{.**.};
\draw[densely dotted] ($(upper_right root 1)+(-0.4,0)$) -- ($(upper_right root 1)+(-0.2,0)$);
\draw ($(upper_right root 1)+(-0.2,0)$) -- ($(upper_right root 1)+(-0.05,0)$);
\draw[densely dotted] ($(upper_right root 2)+(-0.4,0)$) -- ($(upper_right root 2)+(-0.2,0)$);
\draw ($(upper_right root 2)+(-0.2,0)$) -- ($(upper_right root 2)+(-0.05,0)$);
\dynkin[at=(up2),name=middle1_right, labels={,}, ply=2,fold radius=2mm]{A}{.oo.};
\draw[densely dotted] ($(middle1_right root 1)+(-0.4,0)$) -- ($(middle1_right root 1)+(-0.2,0)$);
\draw ($(middle1_right root 1)+(-0.2,0)$) -- ($(middle1_right root 1)+(-0.05,0)$);
\draw[densely dotted] ($(middle1_right root 2)+(-0.4,0)$) -- ($(middle1_right root 2)+(-0.2,0)$);
\draw ($(middle1_right root 2)+(-0.2,0)$) -- ($(middle1_right root 2)+(-0.05,0)$);
\dynkin[at=(up3),name=middle_right, labels={,,}, ply=2,fold radius=2mm]{A}{.***.};
\draw[densely dotted] ($(middle_right root 1)+(-0.4,0)$) -- ($(middle_right root 1)+(-0.2,0)$);
\draw ($(middle_right root 1)+(-0.2,0)$) -- ($(middle_right root 1)+(-0.05,0)$);
\draw[densely dotted] ($(middle_right root 3)+(-0.4,0)$) -- ($(middle_right root 3)+(-0.2,0)$);
\draw ($(middle_right root 3)+(-0.2,0)$) -- ($(middle_right root 3)+(-0.05,0)$);
\dynkin[at=(up4),name=lower_right, labels={,,},ply=2,fold radius=2mm]{A}{.*o*.};
\draw[densely dotted] ($(lower_right root 1)+(-0.4,0)$) -- ($(lower_right root 1)+(-0.2,0)$);
\draw ($(lower_right root 1)+(-0.2,0)$) -- ($(lower_right root 1)+(-0.05,0)$);
\draw[densely dotted] ($(lower_right root 3)+(-0.4,0)$) -- ($(lower_right root 3)+(-0.2,0)$);
\draw ($(lower_right root 3)+(-0.2,0)$) -- ($(lower_right root 3)+(-0.05,0)$);
\draw [decorate, decoration = {calligraphic brace}] ($(lower_right root 3)+(-7mm,-5mm)$) --  ($(upper_right root 1)+(-7mm,2mm)$);
\node[align=left] at ($(up4) + (2mm,-7mm)$) {\tiny{if $m+1$}\\[-2mm] \tiny{$ = 2rd$}};
\node at ($(upper root 12)+(-2mm,-15mm)$) {$\+$};
\node (current) at ($(upper root 1)+(0,-40mm)$) {};
\dynkin[at=(current),name=middle,labels={,,d',,,r'd',,,,n-m-r'd',,,,,,n-m-1}, ply=2,fold radius=2mm]{A}{*.*o*.*o*.*.*.*o*.*o*.*};
\node (low1) at ($(current)+(4.5,1.3)$) {};
\node (low2) at ($(low1)+(0,-0.75)$) {};
\node (low3) at ($(low2)+(0,-0.75)$) {};
\node (low4) at ($(low3)+(0,-0.75)$) {};
\dynkin[at=(low1),name=upper_right2, labels={,},ply=2,fold radius=2mm]{A}{.**.};
\draw[densely dotted] ($(upper_right2 root 1)+(-0.4,0)$) -- ($(upper_right2 root 1)+(-0.2,0)$);
\draw ($(upper_right2 root 1)+(-0.2,0)$) -- ($(upper_right2 root 1)+(-0.05,0)$);
\draw[densely dotted] ($(upper_right2 root 2)+(-0.4,0)$) -- ($(upper_right2 root 2)+(-0.2,0)$);
\draw ($(upper_right2 root 2)+(-0.2,0)$) -- ($(upper_right2 root 2)+(-0.05,0)$);
\dynkin[at=(low2),name=upper2_right2, labels={,},ply=2,fold radius=2mm]{A}{.oo.};
\draw[densely dotted] ($(upper2_right2 root 1)+(-0.4,0)$) -- ($(upper2_right2 root 1)+(-0.2,0)$);
\draw ($(upper2_right2 root 1)+(-0.2,0)$) -- ($(upper2_right2 root 1)+(-0.05,0)$);
\draw[densely dotted] ($(upper2_right2 root 2)+(-0.4,0)$) -- ($(upper2_right2 root 2)+(-0.2,0)$);
\draw ($(upper2_right2 root 2)+(-0.2,0)$) -- ($(upper2_right2 root 2)+(-0.05,0)$);
\dynkin[at=(low3),name=middle_right2, labels={,,}, ply=2,fold radius=2mm]{A}{.***.};
\draw[densely dotted] ($(middle_right2 root 1)+(-0.4,0)$) -- ($(middle_right2 root 1)+(-0.2,0)$);
\draw ($(middle_right2 root 1)+(-0.2,0)$) -- ($(middle_right2 root 1)+(-0.05,0)$);
\draw[densely dotted] ($(middle_right2 root 3)+(-0.4,0)$) -- ($(middle_right2 root 3)+(-0.2,0)$);
\draw ($(middle_right2 root 3)+(-0.2,0)$) -- ($(middle_right2 root 3)+(-0.05,0)$);
\dynkin[at=(low4),name=lower_right2, labels={,,},ply=2,fold radius=2mm]{A}{.*o*.};
\draw[densely dotted] ($(lower_right2 root 1)+(-0.4,0)$) -- ($(lower_right2 root 1)+(-0.2,0)$);
\draw ($(lower_right2 root 1)+(-0.2,0)$) -- ($(lower_right2 root 1)+(-0.05,0)$);
\draw[densely dotted] ($(lower_right2 root 3)+(-0.4,0)$) -- ($(lower_right2 root 3)+(-0.2,0)$);
\draw ($(lower_right2 root 3)+(-0.2,0)$) -- ($(lower_right2 root 3)+(-0.05,0)$);
\draw [decorate, decoration = {calligraphic brace}] ($(lower_right2 root 3)+(-7mm,-5mm)$) --  ($(upper_right2 root 1)+(-7mm,2mm)$);
\node[align=left] at ($(low4) + (2mm, -7mm)$) {\tiny{if $n-m$}\\[-2mm] \tiny{$ = 2r'd'$}};
\end{tikzpicture}
\end{center}
\begin{enumerate}
    \item Here, 
\begin{align*}
E'_{\Delta_0} &= \spn[\R]{e_i-e_{i+1} \mid 1\leq i\leq rd, d\nmid i}\\
 & \oplus \spn[\R]{e_{m-i+1}-e_{m-i+2} \mid 1\leq i\leq rd, d\nmid i}\\
 & \oplus \spn[\R]{e_i-e_{i+1} \mid rd+1\leq i\leq m-rd}\\
 & \oplus\spn[\R]{e_{i+m+1}-e_{i+m+2} \mid 1\leq i\leq r'd', d'\nmid i}\\
 & \oplus \spn[\R]{e_{n-i+1}-e_{n-i+2} \mid 1\leq i\leq r'd', d'\nmid i}\\
 & \oplus \spn[\R]{e_i-e_{i+1} \mid r'd'+1+m+1\leq i\leq n-r'd'+1}.
\end{align*}
Therefore, $\dim E'_{\Delta_0} = 2r(d-1)+2r'(d'-1)+m-2rd+n-m-1-2r'd' =n-1-2(r+r')$.
\item It follows, that $\dim (E'_{\Delta_0})^{\perp} = 2(r+r'+1)$. \\
 Let $v_i^{(d)} := v_{i,0}^{(d)} = \sum_{j=1}^d e_{(i-1)d+j}$ and $w_i^{(d,m)}:=\sum_{j=1}^d e_{(i-1)d+j+m+1}$. We have  \begin{align*}
   (E'_{\Delta_0})^{\perp} = &\spn[\R]{v_{i}^{(d)} \middle| 
    1\leq i \leq r } \oplus \spn[\R]{v_{i}^{(d)} \middle| 
    \frac{m+1}{d}-r+1\leq i \leq \frac{m+1}{d}}\\
  \oplus & \spn[\R]{\sum_{i=rd+1}^{m-rd+1}e_{i}}
  \oplus \spn[\R]{w_{i}^{(d',m)} \middle| 
    1\leq i \leq r'd'} \\
    \oplus &\spn[\R]{w_{i}^{(d',m)} \middle| 
    \frac{n-m}{d'}-r'+1\leq i \leq \frac{n-m}{d'}}
  \oplus \spn[\R]{\sum_{i=r'd'+m+2}^{n-r'd'+1}e_{i}}, 
  \end{align*}
  as then $\dim (E'_{\Delta_0})^{\perp} = 2r+2r'+2$.
  \item Similar to the cases before, we get 
  \begin{align*}
      E^{\Pi} 
      = &\spn[\R]{\alpha_i+\alpha_{m-i+1}\middle| 1\leq i \leq \lfloor\frac{m}{2}\rfloor} \oplus \spn[\R]{\alpha_{i+m+1}+\alpha_{n-i+1}\middle| 1\leq i \leq \lfloor\frac{n-m-1}{2}\rfloor}\\
      \oplus &\spn[\R]{\alpha_{\frac{m+1}{2}}, \alpha_{\frac{n-m}{2}}\middle| \text{ if } \frac{m+1}{2}\in \N, \frac{n-m}{2}\in\N}.
  \end{align*} 
  Then
  \begin{align*}   
   &\sum_{j=id}^{\lfloor \frac{m+1}{2}\rfloor} d(\alpha_{j}+\alpha_{m-j+1}) +\sum_{j=1}^{d-1} (d-j)(\alpha_{id-j}+\alpha_{m-id+j+1})\\
   =&(-v_{i}^{(d)}+v_{\frac{n+1}{d}-i+1}^{(d)})\in E_s,   
\end{align*}
and 
\begin{align*}   
   &\sum_{j=id'+(m+1)}^{\lfloor \frac{n-m}{2}\rfloor} d'(\alpha_{j}+\alpha_{n-j+1}) +\sum_{j=1}^{d'-1} (d'-j)(\alpha_{id'-j+(m+1)}+\alpha_{n-id'+j+1})\\
   =&(-w_{i}^{(d',m)}+w_{\frac{n-m}{d'}-i+1}^{(d',m)})\in E_s,   
\end{align*}
  where $1\leq i \leq r$. Therefore,
  $$ \scaleleftright[1.75ex]{<}
{\left.\begin{array}{l} v_{i}^{(d)}-v_{\frac{m+1}{d}-i+1}^{(d)}\\[0.3cm]
  w_{j}^{(d)}-w_{\frac{n-m}{d'}-j+1}^{(d',m)}
  \end{array} \right|\begin{array}{l} 1\leq i\leq r\\ 1\leq j \leq r'\end{array}}{>}= E_s.$$
  As expected, $\dim E_s =r+r'$.
  \item Again, we proceed as in the previous case and obtain
  \begin{align*}
      B_{1,1} & := \left\{e_{i}-e_{i+1}\middle| 1\leq i \leq rd, d \nmid i \right\}\\
B_{1,2} & := \left\{e_{m-i+1}-e_{m-i+2}\middle| 1\leq i \leq rd, d \nmid i\right\}\\
 B_{1,3} & := \left\{e_{id+1}+e_{m-id}\middle| 0\leq i \leq r-1\right\}\\
B_{2,1} & := \left\{e_{i+m+1}-e_{i+1+m+1}\middle| 1\leq i \leq r'd', d' \nmid i \right\}\\
B_{2,2} & := \left\{e_{n-i+1}-e_{n-i+2}\middle| 1\leq i \leq r'd', d' \nmid i\right\}\\
 B_{2,3} & := \left\{e_{id'+m+2}+e_{n-id'}\middle| 0\leq i \leq r'-1\right\}\\
B_3 & := \left\{e_{i} \middle|rd +1 \leq i \leq m-rd+1 \right\}\\
B_4 & := \left\{e_{i+m+1} \middle|r'd' +1 \leq i \leq n-m-r'd'+1 \right\}\\
B & := B_{1,1} \cup B_{1,2} \cup B_{1,3} \cup B_{2,1} \cup B_{2,2} \cup B_{2,3} \cup B_3 \cup B_4.
  \end{align*}
  The elements of $B$ are linearly independent and $$|B| = 2r(d-1)+r+2r'(d'-1)+r'+m-2rd+n-m-2r'd' = n-(r+r')=\dim E'_a,$$ so $B$ is a basis for $E'_a$.
  \item As previously, we note that 
  \begin{align*}
      \Phi_a := \Phi\cap E'_a = \scaleleftright[1.75ex]{<}
{\left.\begin{array}{l}
B_{1,1}, B_{1,2},B_{2,1}, B_{2,2},\\
e_{i}-e_{i+1},\\
e_{j+m+1}-e_{j+1+(m+1)}\\
e_{rd+1}-e_{r'd'+m+2}
\end{array}\right| \begin{array}{l}
    rd+1\leq i \leq m-rd \\
    r'd'+1 \leq j \leq n-m-2r'd'+1
\end{array}}{>}.
  \end{align*}
  Again, $B_{1,1}$ and $B_{1,2}$ give rise to a root system of type $A_{d-1}^{2r}$ and analogously, $B_{2,1}$ and $B_{2,2}$ give rise to a root system of type $A_{d'-1}^{2r'}$. The set $$\left\{ \begin{array}{l}e_{i}-e_{i+1},\\
e_{j+m+1}-e_{j+1+(m+1)},\\
e_{rd+1}-e_{r'd'+m+2}
\end{array}\middle| \begin{array}{l}
    rd+1\leq i \leq m-rd+1 \\
    r'd'+1 \leq j \leq n-m-2r'd'+1
\end{array}\right\}$$ gives rise to a root system of type $A_{n-2(rd+r'd')}$. 
\end{enumerate}

We write the results in the following table.
 \newpage\newgeometry{top=0.6cm,bottom=0.6cm,left=1cm, right=1cm,foot=0cm}
\begin{landscape}
\vspace*{0.8cm}
\begin{small}
\begin{longtable}{| c | c | c | c | c | c | c | c | c | c | c |}\caption{Isotropic maximal embeddings of abstract indices of type $A_n$}\label{A_n} \\ \hline

\multicolumn{1}{|c|}{} & \multicolumn{1}{c|}{} & \multicolumn{1}{c|}{} & \multicolumn{1}{c|}{} & \multicolumn{1}{c|}{} & \multicolumn{1}{c|}{} & \multicolumn{1}{c|}{} & \multicolumn{3}{c|}{\small{Special fields}} \\

  $\Delta$ & \small{$\hspace{-0.5mm}\!\Stab_{\Iso(\Phi)}(\Delta)\hspace{-0.5mm}$\!\hspace{-0.5mm}} & $\Pi$ & $\mathcal{H}$ & $\Psi_0$ & $\Phi_a$ & $\mathcal{G}$ & \small{$\hspace{-0.8mm}\!\cd 1 \!\hspace{-1mm}$} & $\!\R\!$ & $\!\!Q_{\mathfrak{p}}\!\!$  \\ \hline \hline \endhead
\hline \multicolumn{11}{|c|}{{Continued on next page}} \\ \hline \endfoot 
 \endlastfoot
 &&&&&&&&&\\[-0.37cm]
\makecell{\!$A_m A_{n-m-1}$\!\\ \!$0\!\leq\! m\!\leq\!\frac{n-1}{2}$\!} & \makecell{$\Z_2$ \\ \!(if $m\!<\!\frac{n-1}{2}$)\!} & \makecell{$\Z_2$} & %
\makecell{\begin{tikzpicture}[scale=0.8,baseline=-5mm] 
\node (up) at (0,2) {};
\dynkin[at=(up),name=upper,labels={,,d,,,rd,,,,{m\!-\!rd\!+\!1},,,,,,m}, ply=2,fold radius=2mm]{A}{*.*o*.*o*.*.*.*o*.*o*.*};
\node (up1) at ($(up)+(4.5,1.3)$) {};
\node (up2) at ($(up1)+(0,-0.75)$) {};
\node (up3) at ($(up2)+(0,-0.75)$) {};
\node (up4) at ($(up3)+(0,-0.75)$) {};
\dynkin[at=(up1),name=upper_right, labels={,},ply=2,fold radius=2mm]{A}{.**.};
\draw[densely dotted] ($(upper_right root 1)+(-0.4,0)$) -- ($(upper_right root 1)+(-0.2,0)$);
\draw ($(upper_right root 1)+(-0.2,0)$) -- ($(upper_right root 1)+(-0.05,0)$);
\draw[densely dotted] ($(upper_right root 2)+(-0.4,0)$) -- ($(upper_right root 2)+(-0.2,0)$);
\draw ($(upper_right root 2)+(-0.2,0)$) -- ($(upper_right root 2)+(-0.05,0)$);
\dynkin[at=(up2),name=middle1_right, labels={,}, ply=2,fold radius=2mm]{A}{.oo.};
\draw[densely dotted] ($(middle1_right root 1)+(-0.4,0)$) -- ($(middle1_right root 1)+(-0.2,0)$);
\draw ($(middle1_right root 1)+(-0.2,0)$) -- ($(middle1_right root 1)+(-0.05,0)$);
\draw[densely dotted] ($(middle1_right root 2)+(-0.4,0)$) -- ($(middle1_right root 2)+(-0.2,0)$);
\draw ($(middle1_right root 2)+(-0.2,0)$) -- ($(middle1_right root 2)+(-0.05,0)$);
\dynkin[at=(up3),name=middle_right, labels={,,}, ply=2,fold radius=2mm]{A}{.***.};
\draw[densely dotted] ($(middle_right root 1)+(-0.4,0)$) -- ($(middle_right root 1)+(-0.2,0)$);
\draw ($(middle_right root 1)+(-0.2,0)$) -- ($(middle_right root 1)+(-0.05,0)$);
\draw[densely dotted] ($(middle_right root 3)+(-0.4,0)$) -- ($(middle_right root 3)+(-0.2,0)$);
\draw ($(middle_right root 3)+(-0.2,0)$) -- ($(middle_right root 3)+(-0.05,0)$);
\dynkin[at=(up4),name=lower_right, labels={,,},ply=2,fold radius=2mm]{A}{.*o*.};
\draw[densely dotted] ($(lower_right root 1)+(-0.4,0)$) -- ($(lower_right root 1)+(-0.2,0)$);
\draw ($(lower_right root 1)+(-0.2,0)$) -- ($(lower_right root 1)+(-0.05,0)$);
\draw[densely dotted] ($(lower_right root 3)+(-0.4,0)$) -- ($(lower_right root 3)+(-0.2,0)$);
\draw ($(lower_right root 3)+(-0.2,0)$) -- ($(lower_right root 3)+(-0.05,0)$);
\draw [decorate, decoration = {calligraphic brace}] ($(lower_right root 3)+(-7mm,-5mm)$) --  ($(upper_right root 1)+(-7mm,2mm)$);
\node[align=left] at ($(up4) + (0,-5mm)$) {\tiny{if $j\!=\!-1$}};
\node[align=center] at ($(upper root 12)+(-5mm,-18mm)$) {\small{$d|(m+1), d\geq 1,$}\\ \small{$j=m-2rd \geq -1$}};
\node at ($(upper root 12)+(-2mm,-25mm)$) {$\+$};
\node (current) at ($(upper root 1)+(0,-45mm)$) {};
\dynkin[at=(current),name=middle,labels={,,d',,,r'd',,,,n\!-\!m\!-\!r'd',,,,,,n\!-\!m\!-\!1}, ply=2,fold radius=2mm]{A}{*.*o*.*o*.*.*.*o*.*o*.*};
\node (low1) at ($(current)+(4.5,1.3)$) {};
\node (low2) at ($(low1)+(0,-0.75)$) {};
\node (low3) at ($(low2)+(0,-0.75)$) {};
\node (low4) at ($(low3)+(0,-0.75)$) {};
\dynkin[at=(low1),name=upper_right2, labels={,},ply=2,fold radius=2mm]{A}{.**.};
\draw[densely dotted] ($(upper_right2 root 1)+(-0.4,0)$) -- ($(upper_right2 root 1)+(-0.2,0)$);
\draw ($(upper_right2 root 1)+(-0.2,0)$) -- ($(upper_right2 root 1)+(-0.05,0)$);
\draw[densely dotted] ($(upper_right2 root 2)+(-0.4,0)$) -- ($(upper_right2 root 2)+(-0.2,0)$);
\draw ($(upper_right2 root 2)+(-0.2,0)$) -- ($(upper_right2 root 2)+(-0.05,0)$);
\dynkin[at=(low2),name=upper2_right2, labels={,},ply=2,fold radius=2mm]{A}{.oo.};
\draw[densely dotted] ($(upper2_right2 root 1)+(-0.4,0)$) -- ($(upper2_right2 root 1)+(-0.2,0)$);
\draw ($(upper2_right2 root 1)+(-0.2,0)$) -- ($(upper2_right2 root 1)+(-0.05,0)$);
\draw[densely dotted] ($(upper2_right2 root 2)+(-0.4,0)$) -- ($(upper2_right2 root 2)+(-0.2,0)$);
\draw ($(upper2_right2 root 2)+(-0.2,0)$) -- ($(upper2_right2 root 2)+(-0.05,0)$);
\dynkin[at=(low3),name=middle_right2, labels={,,}, ply=2,fold radius=2mm]{A}{.***.};
\draw[densely dotted] ($(middle_right2 root 1)+(-0.4,0)$) -- ($(middle_right2 root 1)+(-0.2,0)$);
\draw ($(middle_right2 root 1)+(-0.2,0)$) -- ($(middle_right2 root 1)+(-0.05,0)$);
\draw[densely dotted] ($(middle_right2 root 3)+(-0.4,0)$) -- ($(middle_right2 root 3)+(-0.2,0)$);
\draw ($(middle_right2 root 3)+(-0.2,0)$) -- ($(middle_right2 root 3)+(-0.05,0)$);
\dynkin[at=(low4),name=lower_right2, labels={,,},ply=2,fold radius=2mm]{A}{.*o*.};
\draw[densely dotted] ($(lower_right2 root 1)+(-0.4,0)$) -- ($(lower_right2 root 1)+(-0.2,0)$);
\draw ($(lower_right2 root 1)+(-0.2,0)$) -- ($(lower_right2 root 1)+(-0.05,0)$);
\draw[densely dotted] ($(lower_right2 root 3)+(-0.4,0)$) -- ($(lower_right2 root 3)+(-0.2,0)$);
\draw ($(lower_right2 root 3)+(-0.2,0)$) -- ($(lower_right2 root 3)+(-0.05,0)$);
\draw [decorate, decoration = {calligraphic brace}] ($(lower_right2 root 3)+(-7mm,-5mm)$) --  ($(upper_right2 root 1)+(-7mm,2mm)$);
\node[align=left] at ($(low4) + (0,-5mm)$) {\tiny{if $j'=-1$}};
\node[align=center] at ($(middle root 12)+(4mm,-24mm)$) {\small{$d'|(n\!-\!m), d'\geq 1,$}\\\small{$j'=n-m-1-2r'd' \geq -1$,}\\\small{$r+r'\geq 1$,}\\
if $m=(n-1)/2$ then $r \leq r'$ \vspace{1mm}};
\end{tikzpicture}} &%
\makecell{$A_{d-1}^{2r}$\\ $\+$ \\ $A_{d'-1}^{2r'}$\\ $\+$ \\ $A_{j}$ \\ $\+$ \\ $A_{j'}$ \\ $\!\!(A_{-1}\!:=\!\varnothing)\!\!$} &
\makecell{$A_{d-1}^{2r}$\\ $\+$ \\ $ A_{d'-1}^{2r'}$\\ $\+$ \\ $\hspace{-0.5mm}\!A_{j+j'+1}\!\hspace{-0.5mm}$} & \makecell{$\begin{tikzpicture}[scale=0.8,baseline=-5mm] 
\node (up) at (0,2) {};
\dynkin[at=(up),name=upper,labels={,,d,,,id,,,,{n-id +1},,,,,,n}, ply=2,fold radius=2mm]{A}{*.*o*.*o*.*.*.*o*.*o*.*};
\node (up1) at ($(up)+(4.5,1.3)$) {};
\node (up2) at ($(up1)+(0,-0.75)$) {};
\node (up3) at ($(up2)+(0,-0.75)$) {};
\node (up4) at ($(up3)+(0,-0.75)$) {};
\dynkin[at=(up1),name=upper_right, labels={,},ply=2,fold radius=2mm]{A}{.**.};
\draw[densely dotted] ($(upper_right root 1)+(-0.4,0)$) -- ($(upper_right root 1)+(-0.2,0)$);
\draw ($(upper_right root 1)+(-0.2,0)$) -- ($(upper_right root 1)$);
\draw[densely dotted] ($(upper_right root 2)+(-0.4,0)$) -- ($(upper_right root 2)+(-0.2,0)$);
\draw ($(upper_right root 2)+(-0.2,0)$) -- ($(upper_right root 2)$);
\dynkin[at=(up2),name=middle1_right, labels={,}, ply=2,fold radius=2mm]{A}{.oo.};
\draw[densely dotted] ($(middle1_right root 1)+(-0.4,0)$) -- ($(middle1_right root 1)+(-0.2,0)$);
\draw ($(middle1_right root 1)+(-0.2,0)$) -- ($(middle1_right root 1)+(-0.05,0)$);
\draw[densely dotted] ($(middle1_right root 2)+(-0.4,0)$) -- ($(middle1_right root 2)+(-0.2,0)$);
\draw ($(middle1_right root 2)+(-0.2,0)$) -- ($(middle1_right root 2)+(-0.05,0)$);
\dynkin[at=(up3),name=middle_right, labels={,,}, ply=2,fold radius=2mm]{A}{.***.};
\draw[densely dotted] ($(middle_right root 1)+(-0.4,0)$) -- ($(middle_right root 1)+(-0.2,0)$);
\draw ($(middle_right root 1)+(-0.2,0)$) -- ($(middle_right root 1)$);
\draw[densely dotted] ($(middle_right root 3)+(-0.4,0)$) -- ($(middle_right root 3)+(-0.2,0)$);
\draw ($(middle_right root 3)+(-0.2,0)$) -- ($(middle_right root 3)$);
\dynkin[at=(up4),name=lower_right, labels={,,},ply=2,fold radius=2mm]{A}{.*o*.};
\draw[densely dotted] ($(lower_right root 1)+(-0.4,0)$) -- ($(lower_right root 1)+(-0.2,0)$);
\draw ($(lower_right root 1)+(-0.2,0)$) -- ($(lower_right root 1)$);
\draw[densely dotted] ($(lower_right root 3)+(-0.4,0)$) -- ($(lower_right root 3)+(-0.2,0)$);
\draw ($(lower_right root 3)+(-0.2,0)$) -- ($(lower_right root 3)$);
\draw [decorate, decoration = {calligraphic brace}] ($(lower_right root 3)+(-7mm,-5mm)$) --  ($(upper_right root 1)+(-7mm,2mm)$);
\node[align=left] at ($(up4) + (0,-5mm)$) {\tiny{if $n\!+\!1$}\\ \tiny{$ = 2id$}};
\end{tikzpicture}$\\[-0.5cm] only if $d=d'$, \\ $(r+r')d \leq id \leq \lceil\frac{n}{2}\rceil$, \\ if $j=-1$ or $j'=-1$ \\ then $i=r+r'$} & \makecell{\checkmark \\ iff $d\!=\!1$, \\ $r\!=\!\lceil \frac{m}{2} \rceil$, \\ \!$r'\!=\!\lfloor\hspace{-0.3mm}\frac{n-m}{2}\hspace{-0.3mm}\rfloor$\! \\ and \\ $i\!=\!\lceil \frac{n}{2}\rceil$} & \makecell{$\checkmark$ \\ iff $d\!=\!1$ \\ \!\!and either\!\! \\ \!\!$i\!=\!r\!+\!r'$ or\!\! \\ $\!\!\min\{t\!+\!1,n\!-\!t\}$;\!\! \\ where \\ $\!\!t\!:=\!m\!+\!r'\!-\!r\!\!$} & \makecell{$\checkmark$ \\ \!\!if $d\!=\!1$,\!\! \\ \!\!$n\!-\!2(r\!+\!r')$\!\! \\ \!\!$=\!\pm 1$ or $0$,\!\! \\ $?$ \\ if $d \!=\!1$, \\ \!\!$n\!-\!2(r\!+\!r')$\!\!\\ $=\!2$ or $3$, \\$\xmark$ \\otherwise} 
\\ \hline
 &&&&&&&&&\\[-0.3cm]
\multirow{10}{*}{{\makecell{$A_m^{n_m}$ \\ $n\!+\!1\!=$ \\
 $(m\!+\!1)n_m,$ \\ $n_m\geq 2$}}} & \multirow{10}{*}{\makecell{$\Z_2\+ S_{n_m}$; \\ where \\ \!\!$\rho:\Z_2\hspace{-0.5mm}\+\hspace{-0.5mm} S_{n_m} $ \\ \!\!$\to S_{n_m}\hspace{-0.5mm}$ is the\!\! \\ \!\!projection\!\!}} & \makecell{\!\!$\Pi\!\subseteq\!W$,\!\!\\ \hspace{-0.5mm}\!\!$\rho(\Pi)$ is a \!\!\hspace{-0.5mm} \\ \hspace{-0.5mm}\!\! primit.\!\! \\ subgrp.\\ of $S_{n_m}$} & \makecell{$\begin{tikzpicture}[scale=0.8,baseline=-5mm] 
\dynkin[name=upper]{A}{*.*o*.*o*.*}
\node (current) at ($(upper root 1)+(0,-10mm)$) {};
\dynkin[at=(current),name=middle,labels={,,d,,,rd,,m}]{A}{*.*o*.*o*.*}
\begin{scope}[on background layer]
\foreach \x in {1,...,8} 
\draw[/Dynkin diagram/fold style]
($(upper root \x)$) -- ($(upper root \x)+(0,-3mm)$);
\foreach \x in {1,...,8} 
\draw[/Dynkin diagram/fold style]
($(middle root \x)+(0,3mm)$) -- ($(middle root \x)$);
\foreach \x in {1,...,8} 
\draw[dotted,black!40,line width=.5mm]
($(upper root \x)+(0,-3mm)$) -- ($(middle
root \x)+(0,3mm)$);
\end{scope}
\end{tikzpicture}$ \\ $d,r\geq 1$, $d(r+1)=m+1$}
 & $A_{d-1}^{n_m(r+1)}$ & $A_{dn_m-1}^{r+1}$ & \makecell{$\dynkin[scale=0.8,labels={,,c,,,ic,,n}]{A}{*...*o*......*o*...*}$\\ where $c|dn_m$, 
 \\  $c(i+1) = n+1$} & \makecell{\checkmark \\ iff \\\!\!$c\!=\!d\!=\!1$\!\!} & \makecell{\checkmark \\ iff $d\!=\!1$, \\ $n_m\!=\!2$} &  \makecell{\checkmark \\ \!\!if $c\!=\!d\!=\!1$,\!\! \\ \!\!$n_m$ is prime\!\!} 
 \\
  &&&&&&&&&\\[-0.3cm] \cline{3-11}
 &&&&&&&&&\\[-0.3cm]
  &  & \makecell{\!\!$\Pi\!\nsubseteq\!W$,\!\!\\ \!$\rho(\Pi)$ is\! \\ \hspace{-0.5mm}\!\! a primit.\!\! \\ \!subgrp.\!\\ of $S_{n_m}$} & \makecell{\begin{tikzpicture}[scale=0.8,baseline=-5mm] 
\dynkin[name=upper]{A}{*.*o*.*o*.*o*.*}
\node (current) at ($(upper root 1)+(0,-10mm)$) {};
\dynkin[at=(current),name=middle,labels={,,d,,,rd, , , , ,m}]{A}{*.*o*.*o*.*o*.*}
\begin{scope}[on background layer]
\foreach \x/\y in {1/11, 2/10}
\draw[very thick,black!40] ($(upper root \x)+(-0.5mm,-0.5mm)$) to [out=45, in=135] ($(upper root \y)+(0.5mm,-0.5mm)$);
\foreach \x in {1,...,11} 
\draw[/Dynkin diagram/fold style]
($(upper root \x)$) -- ($(upper root \x)+(0,-3mm)$);
\foreach \x in {1,...,11} 
\draw[/Dynkin diagram/fold style]
($(middle root \x)+(0,3mm)$) -- ($(middle root \x)$);
\foreach \x in {1,...,11} 
\draw[dotted,black!40,line width=.5mm]
($(upper root \x)+(0,-3mm)$) -- ($(middle root \x)+(0,3mm)$);
\end{scope}
\draw [decorate,
	decoration = {calligraphic brace}] ($(middle root 8)+(1mm,-2mm)$) --  ($(middle root 7)+(-1mm,-2mm)$);
 \node at ($(middle root 8)+(-2mm,-5mm)$) {\tiny $j$};
\end{tikzpicture}\\ $d,r\geq 1$, $d|(m+1)$,\\ $j=m-2rd\geq -1$} & %
\makecell{$A_{d-1}^{2rn_m}$ \\ $\+$ \\ $A_j^{n_m}$} & %
\makecell{$A_{dn_m-1}^{2r}$ \\ $\+$ \\ \!\hspace{-0.5mm}$A_{n-2n_mrd}$\!} &  %
\makecell{\\[-1em]$\begin{tikzpicture}[scale=0.8,baseline=-5mm] 
\node (up) at (0,2) {};
\dynkin[at=(up),name=upper,labels={,,c,,,ic,,,,{n-ic +1},,,,,,n}, ply=2,fold radius=2mm]{A}{*.*o*.*o*.*.*.*o*.*o*.*};
\node (up1) at ($(up)+(4.5,1.3)$) {};
\node (up2) at ($(up1)+(0,-0.75)$) {};
\node (up3) at ($(up2)+(0,-0.75)$) {};
\node (up4) at ($(up3)+(0,-0.75)$) {};
\dynkin[at=(up1),name=upper_right, labels={,},ply=2,fold radius=2mm]{A}{.**.};
\draw[densely dotted] ($(upper_right root 1)+(-0.4,0)$) -- ($(upper_right root 1)+(-0.2,0)$);
\draw ($(upper_right root 1)+(-0.2,0)$) -- ($(upper_right root 1)$);
\draw[densely dotted] ($(upper_right root 2)+(-0.4,0)$) -- ($(upper_right root 2)+(-0.2,0)$);
\draw ($(upper_right root 2)+(-0.2,0)$) -- ($(upper_right root 2)$);
\dynkin[at=(up2),name=middle1_right, labels={,}, ply=2,fold radius=2mm]{A}{.oo.};
\draw[densely dotted] ($(middle1_right root 1)+(-0.4,0)$) -- ($(middle1_right root 1)+(-0.2,0)$);
\draw ($(middle1_right root 1)+(-0.2,0)$) -- ($(middle1_right root 1)+(-0.05,0)$);
\draw[densely dotted] ($(middle1_right root 2)+(-0.4,0)$) -- ($(middle1_right root 2)+(-0.2,0)$);
\draw ($(middle1_right root 2)+(-0.2,0)$) -- ($(middle1_right root 2)+(-0.05,0)$);
\dynkin[at=(up3),name=middle_right, labels={,,}, ply=2,fold radius=2mm]{A}{.***.};
\draw[densely dotted] ($(middle_right root 1)+(-0.4,0)$) -- ($(middle_right root 1)+(-0.2,0)$);
\draw ($(middle_right root 1)+(-0.2,0)$) -- ($(middle_right root 1)$);
\draw[densely dotted] ($(middle_right root 3)+(-0.4,0)$) -- ($(middle_right root 3)+(-0.2,0)$);
\draw ($(middle_right root 3)+(-0.2,0)$) -- ($(middle_right root 3)$);
\dynkin[at=(up4),name=lower_right, labels={,,},ply=2,fold radius=2mm]{A}{.*o*.};
\draw[densely dotted] ($(lower_right root 1)+(-0.4,0)$) -- ($(lower_right root 1)+(-0.2,0)$);
\draw ($(lower_right root 1)+(-0.2,0)$) -- ($(lower_right root 1)$);
\draw[densely dotted] ($(lower_right root 3)+(-0.4,0)$) -- ($(lower_right root 3)+(-0.2,0)$);
\draw ($(lower_right root 3)+(-0.2,0)$) -- ($(lower_right root 3)$);
\draw [decorate, decoration = {calligraphic brace}] ($(lower_right root 3)+(-7mm,-5mm)$) --  ($(upper_right root 1)+(-7mm,2mm)$);
\node[align=left] at ($(up4) + (0,-5mm)$) {\tiny{if $n\!+\!1$}\\ \tiny{$ = 2ic$}};
\end{tikzpicture}$\\[-0.5cm] where $c|dn_m$, 
\\ $rdn_m \leq ic \leq \lceil \frac{n}{2}\rceil$} & %
\makecell{\checkmark \\ iff \\\!\!$c\!=\!d\!=\!1$,\!\! \\ \!\!$i\!=\!\lceil \frac{n}{2}\rceil$\!\!} & $\xmark$ & \makecell{$\checkmark$ \\ \hspace{-0.5mm}\!\!if\hspace{-0.2mm} $c\!=\!d\!=\!1$,\!\!\hspace{-0.5mm} \\ \!\!$j\!\in\!\{0,-\hspace{-0.5mm}1\}$,\hspace{-0.5mm}\!\! \\ \!\!$n\!\leq\!2i$,\!\! 
\\ \!\!$n_m$ is prime,\hspace{0.5mm}\!\!\\ $\xmark$ \\ \!\!if $c,d,j$ or \!\! \\ $n\!-\!2ic$ \\ are $>1$} 
\\[-0.3cm]
 &&&&&&&&&\\ \hline 
\end{longtable}
\end{small}
\end{landscape}
\restoregeometry
\subsection{\texorpdfstring{ \hspace*{-0.65cm}. $\Phi=B_n=\dynkin[labels={\alpha_1, \alpha_2,\alpha_{n-1},\alpha_n},scale=2.2,
edge length=.5cm, root radius = 0.03cm] B{**...**}$, $n\geq 2$}{Bn}} \label{Bn}

Let $\Phi=B_n$, $n\geq 2$. We use the following standard parametrisation in $\R^n$ of $\Phi$. Let $\alpha_i=e_i-e_{i+1}$ for $1 \leq i \leq n-1$ and $\alpha_n=e_n$, where $e_i$ denotes the $i$'th standard basis vector.\\

By Table II of \cite{T2}, the isotropic indices of type $\Phi$ consist of the following diagrams for $0 \leq j \leq n$:
$$\begin{tikzpicture}[scale=0.8,baseline=-2.2ex] \dynkin[ply=1,rotate=0,fold radius = 2mm,labels={,n-j,,}]{B}{o...o*...**} \end{tikzpicture}$$

Recall from Table \ref{classlist} that $\Delta$ is one of the following types: 
\begin{enumerate}
\item $B_m^{n_m}$ with $mn_m = n$ and $n_m \geq 2$, if $p=2$.
\item $B_mB_{n-m}$ for $1 \leq m \leq n-1$, if $p=2$.
\item $D_n$, for all $p$ 
\item $D_mB_{n-m}$ with $2 \leq m < n$, if $p \neq 2$.
\item $B_{n-1}$, if $p \neq 2$.
\end{enumerate}

In the fourth column of Table \ref{B_n}, we list all indices $\mathcal{H}$ of type $\Delta$ that are maximal in $\Phi$. Hence, we will also be needing the Tits isotropic indices of type $D_m$. It depends on whether there is an order 2 graph automorphism swapping the roots $\alpha_{m-1}, \alpha_m$. Assume first that there is no graph automorphism, then the isotropic indices $\mathcal{H}$ consists of the following \cite[Table II]{T2}:

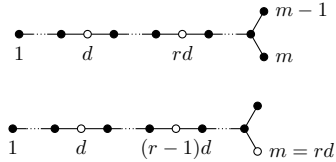
$$\dynkin[scale=1,labels={1,,d,,,rd,,,m-1,m}, rotate=0]{D}{*...*o*......*o*...***}$$

$$\dynkin[scale=1,labels={1,,d,,,(r-1)d,,,,m=rd},rotate=0]{D}{*...*o*......*o*...**o}$$
where $m=rd$ and $r,d \geq 1$. Now assume that there is a graph automorphism that swaps the roots $\alpha_{m-1}$ and $\alpha_{m}$. Then the isotropic indices $\mathcal{H}$ consists of the following \cite[Table II]{T2}:

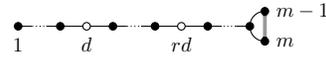
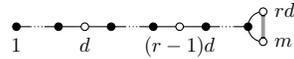
$$\dynkin[scale=1,labels={1,,d,,,rd,,,m-1,m}, ply=2,rotate=0,fold radius = 2mm]{D}{*...*o*......*o*...***}$$
where $r,d \geq 1$.

$$\dynkin[scale=1,labels={1,,d,,,(r-1)d,,,rd,m}, ply=2,rotate=0,fold radius = 2mm]{D}{*...*o*......*o*...*oo}$$
where $m=rd+1$, $r \geq 1$, $d \in \{1, 2 \}$. \\

Suppose $\Delta=D_mB_{n-m}$ where $2 \leq m \leq n$ and $p \neq 2$ or $\Delta=D_n$ for all $p$. Then we have that $S=\Z_2$ by Table \ref{classlist} and suppose first that $\Pi = 1$. If $\mathcal{H}$ is $\dynkin[scale=1,labels={m-1,,m-d,,,m-rd,,,1,}, rotate=180, vertical shift =-.5ex]{D}{*...*o*......*o*...***}$ $\times$ $\dynkin[labels={m,,n-i+1,,n}]{B}{o...o*...**}$ and $\Psi_0=D_jA_{d-1}^rB_i$, where $j=m-rd \geq 2$, $d=2^a|2m ,d \geq 1$ and $0 \leq i \leq n-m$. Hence: $$\Delta_0=\{\alpha_0, \alpha_1,...,\alpha_{j-1} \} \cup_{i=0}^{r-1} \{\alpha_{j+id+1},... \alpha_{j+(i+1)d-1} \} \cup \{\alpha_{n-i+1},...,\alpha_n \}.$$  

We can compute $\dim E_s=n-(m+i-r)$ by counting the number of white orbits in $\mathcal{H}$. The following set of $\Pi$-stable simple roots of cardinality $n-(m+i-r)$ then must span $E_s$ hence $\Phi_s=\Phi \cap E_s$:
$$\Phi_s= \{\alpha_{cd} \text{  $|$ } c=1,2,...,r \} \cup \{\alpha_{m},...,\alpha_{n-i}  \}$$

Therefore $\Phi_a=(\Phi_s)^\perp \cong A_{d-1}^rB_l$ where $l=i+j$.\\

Suppose now $\Pi = \Z_2$ and $\mathcal{H}$ is $\dynkin[scale=1,labels={m-1,,m-d,,,m-rd,,,,}, ply=2,rotate=180,fold radius = 2mm, vertical shift =-.5ex]{D}{*...*o*......*o*...***}$ $\times$ $\dynkin[labels={m,,n-i+1,,n}]{B}{o...o*...**}$ and $\Psi_0=D_jA_{d-1}^rB_i$, where $j=m-rd \geq 2$, $d=2^a|2m, d \geq 1$ and $0 \leq i \leq n-m$. Then an identical argument to the case where $\Pi=1$ shows that $\Phi_a=(\Phi_s)^\perp \cong A_{d-1}^rB_l$ where $l=i+j$.\\

If $\mathcal{H}$ is $\dynkin[scale=1,labels={m-1,,m-d,,,d+1,,,,}, ply=2,rotate=180,fold radius = 2mm, vertical shift =-.5ex]{D}{*...*o*......*o*...*oo}$ $\times$ $\dynkin[labels={,,n-i+1,,}]{B}{o...o*...**}$ and $\Psi_0=A_{d-1}^rB_i$, where $m-rd=1$, $d \in \{1,2 \}$ and $0 \leq i \leq n-m$. Then $S=\Z_2$ and $\Pi=\Z_2$. Hence: $$\Delta_0=\cup_{c=1}^{r-1} \{\alpha_{cd+1},... \alpha_{(c+1)d} \} \cup \{\alpha_{n-c+1},...,\alpha_n \}.$$  

We compute that $\dim E_s=n-(m+i-r)$. The following set of $\Pi$-stable roots of cardinality $n-(m+i-r)$ then must span $E_s$ hence $\Phi_s=\Phi \cap E_s$:

$$\Phi_s=\{\tfrac{1}{2}(\alpha_0+\alpha_1)\} \cup \{\alpha_{cd+1} \text{  $|$ } c=1,2,...,r-1 \} \cup \{\alpha_{m},...,\alpha_{n-i}  \}$$

Therefore $\Phi_a=(\Phi_s)^\perp \cong A_{d-1}^rB_{i+1}$. \\

If $\mathcal{H}$ is $\dynkin[scale=1,labels={,,(r-1)d,,,d,,,1,},rotate=180, vertical shift =-.5ex]{D}{*...*o*......*o*...**o}$ $\times$ $\dynkin[labels={,,n-i+1,,}]{B}{o...o*...**}$ and $\Psi_0=A_{d-1}^rB_i$, where $m=rd$, $d=2^a|2m, d \geq 1$ and $0 \leq i \leq n-m$. Then $S=1$ and therefore $\Pi=1$. Hence: $$\Delta_0=\cup_{c=0}^{r-1} \{\alpha_{cd+1},... \alpha_{(c+1)d-1} \} \cup \{\alpha_{n-i+1},...,\alpha_n \}.$$  

We compute that $\dim E_s=n-(m+i-r)$. The following set of $\Pi$-stable roots of cardinality $n-(m+i-r)$ then must span $E_s$ hence $\Phi_s=\Phi \cap E_s$:

$$\Phi_s= \{\alpha_{cd} \text{  $|$ } c=0,1,2,...,r-1 \} \cup \{\alpha_{m},...,\alpha_{n-i}  \}$$

Therefore $\Phi_a=(\Phi_s)^\perp \cong A_{d-1}^rB_i$. \\

Suppose $\Delta=B_{n-1}$ and $p \neq 2$. We have that $S=\Z_2$ by Table \ref{classlist}. We may ignore the case $\Pi = 1$, as any associated index $\mathcal{H}$ is not maximal in $\Phi$. So we need only consider the case $\Pi = \Z_2$. Then $\mathcal{H}$ is $\dynkin[labels={,,n-i,,}]{B}{o...o*...**} \times {\mathcal{T}_0}^1$  and $\Psi_0=B_i$, where $0 \leq i \leq 2$. Hence: $$\Delta_0=\{ \alpha_{n-i},...,\alpha_{n-1} \}.$$  

We compute that $\dim E_s=n-(i+1)$. The following set of $\Pi$-stable simple roots of cardinality $n-(i+1)$ then must span $E_s$ hence $\Phi_s=\Phi \cap E_s$:

$$\Phi_s= \{ \alpha_1,...,\alpha_{n-(i+1)} \}.$$

Therefore $\Phi_a=(\Phi_s)^\perp \cong B_{i+1}$.\\

Suppose $\Delta=B_{m}^{n_m}$ where $mn_m=n$, $n_m \geq 2$ and $p=2$. We have that $S=S_{n_m}$ by Table \ref{classlist} and $\Pi$ is a primitive subgroup of $S_{n_m}$. If $\mathcal{H}$ is $(\dynkin[labels={,,m-i+1,,}]{B}{o...o*...**})^{n_m}$ and $\Psi_0=B_{i}^{n_m}$ where $0 \leq i \leq m-1$.\\

One checks that: $$\Delta_0= \cup_{c=0}^{{n_m}-2} \{ -(\alpha_{cm+1} + \alpha_{cm+2}+...+\alpha_n),\alpha_{cm+1},...,\alpha_{cm+i-1} \} \cup \{ \alpha_{n-i+1},...,\alpha_{n} \} \text{ if } i>0,$$  
$$\Delta_0= \emptyset \text{ if } i=0.$$  

The subspace $E_{\Delta_0}$ has dimension $in_m$ and $E_{\Delta_0}^{\perp}$ has dimension $n-in_m$. We compute that: $$E_{\Delta_0}^\perp =  \oplus_{r=1}^{m-i} \oplus_{c=0}^{n_m-2} \langle  e_{(c+1)m-r+1}) \rangle \oplus_{r=1}^{m-i} \langle e_{n-(m-r)} \rangle.$$

The action of $\Pi$ on the basis vectors $e_i$ is given here: $\Pi(e_i)=e_{i+m}$ for all $1 \leq i \leq m(n_m-2)$, $\Pi(e_i)=-e_{2n+1-(i+2m)}$ for all $m(n_m-2) < i \leq n-m$ and $\Pi(e_i)=-e_{n-i+1}$ for all $n-m < i \leq n$. Hence the fixed point space $E_s$ with dimension $m-i$ is computed: 
$$E_s= \oplus_{r=1}^{m-i} \langle (\Sigma_{c=0}^{n_m-2} e_{(c+1)m-r+1}) - (e_{n-(m-r)} \rangle.$$
Then the dimension of $E_a$ is $n-(m-i)$ and $E_a=B_1 + B_2$ where:
$$B_1=\oplus_{r=1}^{m-i} \oplus_{c=0}^{n_m-2} \langle e_{(c+1)m-r+1} + (e_{n-(m-r)} \rangle,$$
$$B_2= \oplus_{c=0}^{n_m-2} \oplus_{r=1}^{i} \langle e_{cm+r} \rangle \oplus_{r=1}^{i} \langle e_{n-r+1}\rangle.$$

Again, $B_1$ gives rise to a root system of type $A_{n_m-1}^{m_i}$ and analogously $B_2$ gives rise to a root system $B_{in_m}$. Then it follows that: $$\Phi_a=E_a \cap \Phi \cong A_{n_m-1}^{m-i}B_{in_m}.$$

Suppose $\Delta=B_{m}B_{n-m}$ where $1 \leq m \leq n-1$ and $p=2$. We have that $S=1$ by Table \ref{classlist} and therefore $\Pi = 1$. If $\mathcal{H}$ is $\dynkin[labels={,,m-i+1,,}]{B}{o...o*...**}$ $\times$ $\dynkin[labels={,,n-m-j+1,,}]{B}{o...o*...**}$ and $\Psi_0=B_{i}B_j$ where $0 \leq i \leq m$, $0 \leq j \leq n-m$ and $i+j < n$. Hence: $$\Delta_0= \{ \alpha_{1}-\alpha_{0},\alpha_2,...,\alpha_{i} \} \cup \{ \alpha_{n-j+1},...,\alpha_{n} \}.$$  

We compute that $\dim E_s=n-(i+j)$. The following set of $\Pi$-stable simple roots of cardinality $n-(i+j)$ then must span $E_s$ hence $\Phi_s=\Phi \cap E_s$:

$$\Phi_s= \{ \alpha_{i+1},...,\alpha_m \} \cup \{ \alpha_{m+1},...,\alpha_{n-j} \}, \text{ if } i> 0$$

$$\Phi_s= \{ \alpha_1-\alpha_0,\alpha_{2},...,\alpha_m \} \cup \{ \alpha_{m+1},...,\alpha_{n-j} \}, \text{ if } i=0$$

Therefore $\Phi_a=(\Phi_s)^\perp \cong B_{i+j}$.\\

\newpage\newgeometry{top=0.5cm,bottom=1cm,left=0.5cm, right=0.5cm,foot=0.5cm}
\begin{landscape}
\begin{longtable}{| c | c | c | c | c | c | c | c | c | c | c |}\caption{Isotropic maximal embeddings of abstract indices of type $B_n$} \label{B_n} \\ \hline

\multicolumn{1}{|c|}{} & \multicolumn{1}{c|}{} & \multicolumn{1}{c|}{} & \multicolumn{1}{c|}{} & \multicolumn{1}{c|}{} & \multicolumn{1}{c|}{} & \multicolumn{1}{c|}{} & \multicolumn{3}{c|}{\small{Special fields}} \\

$\Delta$ & \small{$\!\!\Stab_W(\Delta)\!\!$} & $\Pi$ & $\textcolor{white}{iiaaa} \mathcal{H} \textcolor{white}{aaaii}$ & $\Psi_0$ & $\Phi_a$ & $\!\!\mathcal{G}\!\!$ & \small{$\!\!\cd 1 \!\!$} & $\!\R\!$ & $\!\!Q_{\mathfrak{p}}\!\!$ \\ \hline \hline 
\endfirsthead
\hline
\multirow{19}{*}{\makecell{$D_{m}B_{n-m}$ \\ $2 \leq m < n$ \\if $p \neq 2$, \\ and $m=n$ \\ for all $p$}} & \multirow{19}{*}{$\Z_2$} & \multirow{7}{*}{$1$} & \makecell{$\dynkin[scale=1,labels={m-1,,m-d,,,m-rd,,,1,}, rotate=180]{D}{*...*o*......*o*...***}$ \\ $\dynkin[labels={,n-m-i,,,}]{B}{o...o*...**}$\\  $j=m-rd \geq 2$, \\ $d=2^a|2m, d \geq 1,$ \\ \!\!\!$0 \leq i \leq n-m$, $i+j<n$\!\!\!} & {$D_{j}A_{d-1}^rB_i$} & \makecell{$A_{d-1}^rB_{l}$ \\ $l=i+j$} & \makecell{$\begin{tikzpicture}[scale=0.8,baseline=-2.2ex] \dynkin[ply=1,rotate=0,fold radius = 2mm,labels={,n-c,,}]{B}{o...o*...**} \end{tikzpicture}$ \\ $0 \leq c\leq l$,\\ only if $d=1$} & $\xmark$ & \makecell{$\checkmark$ \\ iff $j$ is even and \\ either $c\!=\!l$ or $c=$\\ $\max\{i\!-\!j,j\!-\!i\!-\!1\}$} & \makecell{$\xmark$ \\unless \\ $i,c\!\leq \!1$ \\ \!and $j\!=\!2$\!} 
\topstrut \\ \cdashline{4-11}

&&& \makecell{$\dynkin[scale=1,labels={,,(r-1)d,,,d,,,1,},rotate=180]{D}{*...*o*......*o*...**o}$ \\ $\dynkin[labels={,n-m-i,,,}]{B}{o...o*...**}$  \\ $m=rd$, $d=2^a|2m, d \geq 1$, \\ $0 \leq i \leq n-m$} & $A_{d-1}^rB_i$ & $A_{d-1}^rB_{i}$ & \makecell{$\begin{tikzpicture}[scale=0.8,baseline=-2.2ex] \dynkin[ply=1,rotate=0,fold radius = 2mm,labels={,n-i,,}]{B}{o...o*...**} \end{tikzpicture}$ \\ only if $d=1$} & \makecell{$\checkmark$ \\ iff $i\!=\!0$} & $\checkmark$ & \makecell{$\checkmark$ \\ iff \\ \!$i\!\in\!\{0,1\}$\!\!}   \topstrut \\ \cline{3-11}

& & \multirow{7}{*}{$\Z_2$} & \!\makecell{$\dynkin[scale=1,labels={m-1,,m-d,,,m-rd,,,,}, ply=2,rotate=180,fold radius = 2mm]{D}{*...*o*......*o*...***}$ \\ $\dynkin[labels={,n-m-i,,,}]{B}{o...o*...**}$ \\ $j=m-rd \geq 2$,\\ $d=2^a|2m, d \geq 1$, \\ \!\!\!$0 \leq i \leq n-m$, $i+j<n$\!\!\!}\! & $D_{j}A_{d-1}^rB_i$ & \makecell{$A_{d-1}^rB_{l}$ \\ $l=i+j$} & \makecell{$\begin{tikzpicture}[scale=0.8,baseline=-2.2ex] \dynkin[ply=1,rotate=0,fold radius = 2mm,labels={,n-c,,}]{B}{o...o*...**}  \end{tikzpicture}$ \\ $0 \leq c \leq l$, \\ only if $d=1$ } & $\xmark$ & \makecell{$\checkmark$ \\ iff $j$ is odd and \\ either $c\!=\!l$ or $c=$\\ $\max\{i\!-\!j,j\!-\!i\!-\!1\}$} & $\xmark$ 
\topstrut \\ \cdashline{4-11}

&&& \!\makecell{$\dynkin[scale=1,labels={m-1,,m-d,,,d+1,,,,}, ply=2,rotate=180,fold radius = 2mm]{D}{*...*o*......*o*...*oo}$ \\ $\dynkin[labels={,n-m-i,,,}]{B}{o...o*...**}$ \\ \!\!$m-rd=1$, $d \in \{1, 2\}$,\!\!\\ $0 \leq i \leq n-m$}\! & \makecell{$A_{d-1}^rB_i$} & \makecell{$A_{d-1}^rB_{i+1}$} & \makecell{$\begin{tikzpicture}[scale=0.8,baseline=-2.2ex] \dynkin[ply=1,rotate=0,fold radius = 2mm,labels={,n-c,,}]{B}{o...o*...**}  \end{tikzpicture}$ \\ $0 \leq c \leq i+1$, \\ only if $d=1$} & \makecell{$\checkmark$ \\ iff \\ \!$i\!=\!c\!=\!0$\!} & \makecell{$\checkmark$ \\ iff $c\!=\!i\!+\!1$ or \\ $c\!=\!\max\{0,i\!-\!1\}$} & \makecell{$\checkmark$ \\if $i\!=\!0$,\\ $\xmark$ \\if $i\!>\!1$ \\ or $c\!>\!1$} \topstrut \\ \cdashline{4-11}

\hline {\makecell{$B_{n-1}$ \\ $p \neq 2$}} & {$\Z_2$} & $\Z_2$ &  \makecell{ $\dynkin[labels={,,n-1-i,,}]{B}{o...o*...**} \times {\mathcal{T}_0}^1$ \\ $0 \leq i \leq n-2$} & $B_{i}$ & $B_{i+1}$ & \makecell{$\begin{tikzpicture}[scale=0.8,baseline=-2.2ex] \dynkin[ply=1,rotate=0,fold radius = 2mm,labels={,n-c,,}]{B}{o...o*...**}  \end{tikzpicture}$ \\ $0 \leq c \leq i+1$} & \makecell{$\checkmark$ \\ iff \\ \!$i\!=\!c\!=\!0$\!} & \makecell{$\checkmark$ \\ iff $c\!=\!i\!+\!1$ or \\ $c\!=\!\max\{0,i\!-\!1\}$} & \makecell{$\checkmark$ \\ if $i\!=\!0$, \\\!$\xmark$ \hspace{-0.2mm}if \hspace{-0.2mm}$i\!>\!1$\!\!\\ or $c\!>\!1$} \topstrut \\ \hline 

{\makecell{$B_{m}B_{n-m}$ \\ $1 \!\leq\! m \!\leq\! \frac{n}{2}$ \\ $p=2$}} & \makecell{$1$ \\ \!\!(if $m\!<\!\frac{n}{2}$)\!\!} 
& {$1$} & \makecell{$\dynkin[labels={,m-i,,,}]{B}{o...o*...**}$ \\ $\dynkin[labels={,n-m-j,,,}]{B}{o...o*...**}$ \\ \!\!\!$0 \hspace{-0.1mm}\leq \hspace{-0.1mm}i \hspace{-0.1mm}\leq \hspace{-0.1mm}m$, \hspace{-0.1mm}$0 \hspace{-0.1mm}\leq \hspace{-0.1mm}j \hspace{-0.1mm}\leq \hspace{-0.1mm} n\hspace{-0.1mm}-\hspace{-0.1mm}m$,\!\!\! \\ $i+j<n$, \\ if $m=n/2$ then $i \leq j$} & $B_iB_j$ & \makecell{$B_{l}$ \\ $l=i+j$} & \makecell{$\begin{tikzpicture}[scale=0.8,baseline=-2.2ex] \dynkin[ply=1,rotate=0,fold radius = 2mm,labels={,n-c,,}]{B}{o...o*...**}  \end{tikzpicture}$ \\ $0 \leq c \leq l$, \\ $c=l$ if $i=0$ \\ or if $j=0$} & \makecell{$\checkmark$ \\ iff \\ $i\!=\!j\!=\!0$} & $\xmark$ & $\xmark$  \topstrut \\ \hline

{\makecell{$B_m^{n_m}$ \\ $mn_m=n$ \\ $n_m \geq 2$ \\ $p=2$}} & {$S_{n_m}$} & \makecell{primit. \\ \!subgrp.\! \\ of $S_{n_m}$} & \makecell{$\begin{tikzpicture}[scale=0.8,baseline=-5mm] 
\dynkin[name=upper]{B}{o.o*.**}
\node (current) at ($(upper root 1)+(0,-10mm)$) {};
\dynkin[at=(current),name=middle,labels={,m-i,,,m}]{B}{o.o*.**}
\begin{scope}[on background layer]
\foreach \x in {1,...,5} 
\draw[/Dynkin diagram/fold style]
($(upper root \x)$) -- ($(upper root \x)+(0,-3mm)$);
\foreach \x in {1,...,5} 
\draw[/Dynkin diagram/fold style]
($(middle root \x)+(0,3mm)$) -- ($(middle root \x)$);
\foreach \x in {1,...,5} 
\draw[dotted,black!40,line width=.5mm]
($(upper root \x)+(0,-3mm)$) -- ($(middle
root \x)+(0,3mm)$);
\end{scope}
\end{tikzpicture}$ \\ $0 \leq i \leq m-1$} & $B_{i}^{n_m}$ & $\!A_{n_m-1}^{m-i}B_{in_m}\!$ & \makecell{$\begin{tikzpicture}[scale=0.8,baseline=-2.2ex] \dynkin[ply=1,rotate=0,fold radius = 2mm,labels={,n-c,,}]{B}{o...o*...**}  \end{tikzpicture}$ \\ $\!\!0 \leq c \leq in_m\!\!$} & \makecell{$\checkmark$ \\ iff $i\!=\!0$} & $\xmark$ & $\xmark$ 
\topstrut \\ \hline
\end{longtable}
\end{landscape}
\newpage \restoregeometry
\subsection{\texorpdfstring{ \hspace*{-0.65cm}. $\Phi=C_n=\dynkin[labels={\alpha_1, \alpha_2,\alpha_{n-1},\alpha_n},scale=2.2,
edge length=.5cm,root radius = 0.03cm] C{**...**}$, $n\geq 3$}{Cn}} \label{Cn}
Let $\Phi=C_n$, $n \geq 3$. We use the following standard parametrisation in $\mathbb{R}^n$ of $\Phi$. Let $\alpha_i=e_i-e_{i+1}$ for $1 \leq i \leq n-1$ and $\alpha_n=2e_n$, where $e_i$ denotes the $i$'th standard basis vector.\\

By Table II of \cite{T2}, the isotropic indices of type $\Phi$ consist of the following diagrams:
$$\dynkin[scale=1,labels={1,,d,,,rd,,,,n}, rotate=0]{C}{*...*o*......*o*...***}$$
where $n-rd \geq 1$, $r \geq 1$, $d \geq 2, d| 2n$. 
$$\dynkin[scale=1,labels={1,,d,,,(r-1)d,,,,n=rd},rotate=0]{C}{*...*o*......*o*...**o}$$
where $n=rd$ and $r,d \geq 1, d | 2n$.

We follow the convention that anisotropic indices of type $C_n$ satisfy $r=0$, and $d$ may be any positive integer other than 1. 

Recall from Table \ref{classlist} that $\Delta$ is one of the following types: 
\begin{enumerate}
\item $C_m^{n_m}$ with $mn_m = n$ and $n_m \geq 2$,
\item $C_mC_{n-m}$ for $1 \leq m \leq n-1$,
\item $A_{n-1}$,
\item $\widetilde{D_n}$ if $p=2$.

\end{enumerate}

Suppose $\Delta=C_{m}^{n_m}$ where $mn_m=n$, $n_m \geq 2$. We have that $S=S_{n_m}$ by Table \ref{classlist} and $\Pi$ is a primitive subgroup $S_{n_m}$. If $\mathcal{H}$ is $(\dynkin[labels={,,d,,,rd,,,,m}]{C}{*...*o*......*o*...***})^{n_m}$ where $j=m-rd, r,d \geq 1$ and if $d=1$ then $j=0$. So $\Psi_0=A_{d-1}^{rn_m}{C_{j}}^{n_m}$.\\

One checks that: $$\Delta= \cup_{b=0}^{{n_m}-2} \{ -(2\alpha_{bm+1}+2\alpha_{bm+2}+...+2\alpha_{n-1}+\alpha_n) \} \cup_{c=1}^{m-1} \{ \cup_{b=1}^{{n_m}-1} \{ \alpha_{bm-c} \} \cup_{c=0}^{m-1} \{ \alpha_{n-c}\} \}.$$  

If $j=0$ then $$\Delta_0=\Delta -(\cup_{b=1}^{n_m-1}\{ \cup_{c=1}^{r-1} \{\alpha_{bm-cd} \} \} \cup_{c=0}^{{n_m}-2} \{ -(2\alpha_{cm+1}+2\alpha_{cm+2}+...+2\alpha_{n-1}+\alpha_n) \} $$ $$\cup_{c=1}^{r} \{\alpha_{(n_m-1)m+cd} \}).$$ 

If $j \geq 1$ then $$\Delta_0=\Delta -(\cup_{b=1}^{n_m-1}\{ \cup_{c=1}^r \{\alpha_{bm-cd} \} \} \cup_{c=1}^r \{\alpha_{(n_m-1)m+cd} \}).$$ 

For all $j \geq 0$, the subspace $E_{\Delta_0}$ has dimension $n-rn_m$ and $E_{\Delta_0}^{\perp}$ has dimension $rn_m$ and we compute that: $$E_{\Delta_0}^\perp = \oplus_{c=1}^{r} \oplus_{b=1}^{n_m-1} \langle \sum_{f=1}^{d} e_{bm-cd+f} \rangle \oplus_{c=1}^{r} \langle \sum_{f=1}^{d} e_{(n_m-1)m+(c-1)d+f} \rangle.$$

The action of $\Pi$ on the basis vectors $e_i$ is given here: $\Pi(e_i)=e_{i+m}$ for all $1 \leq i \leq m(n_m-2)$, $\Pi(e_i)=-e_{2n+1-(i+2m)}$ for all $m(n_m-2) < i \leq n-m$ and $\Pi(e_i)=-e_{n-i+1}$ for all $n-m < i \leq n$. Hence the fixed point space $E_s$ with dimension $r$ is computed: $$E_s= \oplus_{c=1}^{r} \langle \sum_{f=1}^{d}  (\sum_{b=1}^{n_m-1} e_{bm-cd+f})  - e_{(n_m-1)m+(c-1)d+f} \rangle.$$
Then it follows that $E_a$ of dimension $n-r$ is given by $E_a=B_1+B_2$ where:
$$B_1=\oplus_{c=1}^{r} (\oplus_{f=2}^{d} \langle e_{m-cd+1}-e_{m-cd+f} \rangle \oplus_{b=2}^{n_m-1} \oplus_{f=1}^{d} \langle e_{m-cd+1}-e_{bm-cd+f}\rangle$$ $$ \oplus_{f=1}^{d} \langle e_{m-cd+1} + e_{(n_m-1)m+(c-1)d+f} \rangle),$$
$$ B_2 = \oplus_{f=1}^{j} \oplus_{b=1}^{n_m-1} \langle e_{(b-1)m+f} \rangle \oplus_{f=1}^{j} \langle e_{n-j+f} \rangle.$$
Again, $B_1$ gives rise to a root system of type $A_{dn_{m}-1}^{r}$ and analogously $B_2$ gives rise to a root system $C_{jn_m}$. Then it follows that $$\Phi_a=E_a \cap \Phi \cong A_{dn_{m}-1}^{r}C_{jn_m} \text{ if } j \geq 1 \text{ and } \Phi_a=E_a \cap \Phi \cong A_{dn_{m}-1}^{r} \text{ if } j=0.$$

Suppose $\Delta=C_{m}C_{n-m}$ where $1 \leq m \leq n-1$. We have that $S=1$ by Table \ref{classlist} and therefore $\Pi = 1$. If $\mathcal{H}$ is $\dynkin[scale=1,labels={,,d,,,rd,,,m}, rotate=0]{C}{*...*o*......*o*...**} \times$  $\dynkin[scale=0.8,labels={,,d',,,r'd',,,n-m}, rotate=0]{C}{*...*o*......*o*...**}$ and $\Psi_0=A_{d-1}^rC_jA_{d'-1}^{r'}C_{j'}$ where $j=m-rd,j'=n-m-r'd'$ and $r, r',d,d' \geq 1$ and if $d=1$ then $m=r$ and if $d'=1$ then $n-m=r'$. Hence: $$\Delta = \{ \alpha_0, \alpha_1,...,\alpha_{m-1} \} \cup \{\alpha_{m+1},...\alpha_n \}$$
$$\Delta_0= \Delta - (\{ \alpha_{m-d},\alpha_{m-2d},...,\alpha_{m-rd}\} \cup \{ \alpha_{m+d'},\alpha_{m+2d'},...,\alpha_{m+r'd'}\} ).$$  

We compute that $\dim E_s=(r+r')$. The following set of $\Pi$-stable roots of cardinality $(r+r')$ then must span $E_s$ hence:

$$\Phi_s= \{ \alpha_{m-d},\alpha_{m-2d},...,\alpha_{m-rd}\} \cup \{ \alpha_{m+d'},\alpha_{m+2d'},...,\alpha_{m+r'd'}\} $$

Therefore $\Phi_a=(\Phi_s)^\perp \cong A_{d-1}^rA_{d'-1}^{r'}C_{j+j'}$.\\

Suppose $\Delta=A_{n-1}$ where $p \neq 2$. We have that $S=\Z_2$ by Table \ref{classlist}. We may ignore the case $\Pi = 1$, as any associated index $\mathcal{H}$ is not maximal in $\Phi$. So we need only consider the case $\Pi = \Z_2$. Then $\mathcal{H}$ is of the following form:
$$\begin{tikzpicture}[scale=0.8,baseline=-5mm] 
\node (up) at (0,2) {};
\dynkin[at=(up),name=upper,labels={,,d,,,rd,,,,{n-rd},,,,,,n-1}, ply=2,fold radius=2mm]{A}{*.*o*.*o*.*.*.*o*.*o*.*};
\node (up1) at ($(up)+(4.5,1.5)$) {};
\node (up2) at ($(up1)+(0,-0.75)$) {};
\node (up3) at ($(up2)+(0,-1)$) {};
\node (up4) at ($(up3)+(0,-0.75)$) {};
\dynkin[at=(up1),name=upper_right, labels={,},ply=2,fold radius=2mm]{A}{.**.};
\draw[densely dotted] ($(upper_right root 1)+(-0.4,0)$) -- ($(upper_right root 1)+(-0.2,0)$);
\draw ($(upper_right root 1)+(-0.2,0)$) -- ($(upper_right root 1)+(0.01,0)$);
\draw[densely dotted] ($(upper_right root 2)+(-0.4,0)$) -- ($(upper_right root 2)+(-0.2,0)$);
\draw ($(upper_right root 2)+(-0.2,0)$) -- ($(upper_right root 2)+(0.01,0)$); 
\dynkin[at=(up2),name=middle1_right, labels={,}, ply=2,fold radius=2mm]{A}{.oo.};
\draw[densely dotted] ($(middle1_right root 1)+(-0.4,0)$) -- ($(middle1_right root 1)+(-0.2,0)$);
\draw ($(middle1_right root 1)+(-0.2,0)$) -- ($(middle1_right root 1)+(0.01,0)$);
\draw[densely dotted] ($(middle1_right root 2)+(-0.4,0)$) -- ($(middle1_right root 2)+(-0.2,0)$);
\draw ($(middle1_right root 2)+(-0.2,0)$) -- ($(middle1_right root 2)+(0.01,0)$);
\node[align=left] at ($(up2) + (0,-5mm)$) {\tiny{if $m=0$}};
\dynkin[at=(up3),name=middle_right, labels={,,}, ply=2,fold radius=2mm]{A}{.***.};
\draw[densely dotted] ($(middle_right root 1)+(-0.4,0)$) -- ($(middle_right root 1)+(-0.2,0)$);
\draw ($(middle_right root 1)+(-0.2,0)$) -- ($(middle_right root 1)+(0.01,0)$);
\draw[densely dotted] ($(middle_right root 3)+(-0.4,0)$) -- ($(middle_right root 3)+(-0.2,0)$);
\draw ($(middle_right root 3)+(-0.2,0)$) -- ($(middle_right root 3)+(0.01,0)$);
\dynkin[at=(up4),name=lower_right, labels={,,},ply=2,fold radius=2mm]{A}{.*o*.};
\draw[densely dotted] ($(lower_right root 1)+(-0.4,0)$) -- ($(lower_right root 1)+(-0.2,0)$);
\draw ($(lower_right root 1)+(-0.2,0)$) -- ($(lower_right root 1)+(0.01,0)$);
\draw[densely dotted] ($(lower_right root 3)+(-0.4,0)$) -- ($(lower_right root 3)+(-0.2,0)$);
\draw ($(lower_right root 3)+(-0.2,0)$) -- ($(lower_right root 3)+(0.01,0)$);
\draw [decorate, decoration = {calligraphic brace}] ($(lower_right root 3)+(-7mm,-5mm)$) --  ($(upper_right root 1)+(-7mm,2mm)$);
\node[align=left] at ($(up4) + (0,-5mm)$) {\tiny{if $m=-1$}}; \end{tikzpicture} \times \mathcal{T}_0^1 $$ and $\Psi_0=A_{d-1}^{2r}A_m$ where $m=n-1-2rd \geq -1$, $d \geq 1$ and $d|n$. Hence: $$\Delta = \{ \alpha_1,...,\alpha_{n-1} \}$$

$$\Delta_0= \Delta - (\{ \alpha_{d},\alpha_{2d},...,\alpha_{2rd}\}).$$  

The subspace $E_{\Delta_0}$ has dimension $n-2r$ and therefore $E_{\Delta_0}^\perp$ has dimension $2r$ and we compute that:
$$E_{\Delta_0}^\perp = \oplus_{b=0}^{2r-1} \langle \Sigma_{c=1}^{d} e_{bd + c} \rangle$$
We compute that $\dim E_s=r$. The action of $\Pi$ on the basis vectors $e_i$ is given here: $\Pi(e_1)=e_1, \Pi(e_n)=e_n$ and $\Pi(e_i)=e_1+e_n-e_{n-i+1}$ for all $2 \leq i \leq n-1$ and hence the fixed point space $E_s$ is computed:
$$E_s = \oplus_{b=0}^{r-1} \langle \Sigma_{c=1}^{d} e_{bd+c}-e_{(2r-1-b)d+c}\rangle.$$
Then it follows that $E_a$ of dimension $n-r$ is given by $E_a=B_1+B_2$ where:
$$B_1 = \oplus_{b=0}^{r-1} \oplus_{c=2}^{d} \langle e_{bd+1}-e_{bd+c}\rangle \oplus_{b=0}^{r-1} \oplus_{c=1}^{d} \langle e_{bd+1}-e_{(2r-1-b)d+c}\rangle, $$
$$B_2=\oplus_{c=1}^{m+1} \langle e_{rd+c} \rangle.$$
Again, $B_1$ gives rise to a root system of type $A_{2d-1}^{r}$ and analogously $B_2$ gives rise to a root system $C_{m+1}$. Therefore $\Phi_a=E_a \cap \Phi \cong A_{2d-1}^{r}C_{m+1}$.\\

Suppose $\Delta=D_n$ where $p=2$. We have that $S=\Z_2$ by Table \ref{classlist} and assume first that $\Pi = 1$. Let $\mathcal{H}$ be $\dynkin[scale=1,labels={,,d,,,rd,,,n-1,n}, rotate=0]{D}{*...*o*......*o*...***}$ where $j=n-rd \geq 2$, $d \geq 1$, $d | 2n$ and $\Psi_0=A_{d-1}^rD_j$. Furthermore, label the highest short root $\alpha_{s}$. Hence: $$\Delta = \{ \alpha_s, \alpha_1,...,\alpha_{n-1} \}$$
$$\Delta_0= \Delta - (\{ \alpha_{n-d},\alpha_{n-2d},...,\alpha_{n-rd}\} ).$$  

We compute that $\dim E_s =r$. The following set of $\Pi$-stable roots of cardinality $r$ then must span $E_s$ hence $\Phi_s=\Phi \cap E_s$:

$$\Phi_s=\{ \alpha_{n-d},\alpha_{n-2d},...,\alpha_{n-rd}\}$$

Therefore $\Phi_a=(\Phi_s)^\perp \cong A_{d-1}^rC_{j}$.\\

Assume now that $\mathcal{H}$ is $\dynkin[scale=1,labels={,,d,,,(r-1)d,,,,n=rd},rotate=0]{D}{*...*o*......*o*...**o}$  where $n=rd$, $d \geq 1$ and $\Psi_0=A_{d-1}^{r}$. Hence:
$$\Delta_0= \Delta - (\{ \alpha_{n-d},\alpha_{n-2d},...,\alpha_{n-(r-1)d}, \alpha_s\} ).$$ 

We compute that $\dim E_s=r$. The following set of $\Pi$-stable roots of cardinality $r$ then must span $E_s$ hence $\Phi_s=\Phi \cap E_s$:

$$\Phi_s=\{ \alpha_{n-d},\alpha_{n-2d},...,\alpha_{n-(r-1)d}, \alpha_s\}$$

Hence $\Phi_a=(\Phi_s)^\perp \cong A_{d-1}^{r}$.\\

Assume now that $\Pi=\Z_2$ and $\mathcal{H}$ is $\dynkin[scale=1,labels={,,d,,,rd,,,n-1,n}, ply=2,rotate=0,fold radius = 2mm]{D}{*...*o*......*o*...***}$ where $j=n-rd \geq 2$, $d \geq 1$, $d | 2n$ and $\Psi_0=A_{d-1}^rD_j$. Hence:
$$\Delta_0= \Delta - (\{ \alpha_{n-d},\alpha_{n-2d},...,\alpha_{n-rd}\} ).$$ 

We compute that $\dim E_s=r$. The following set of $\Pi$-stable roots of cardinality $r$ then must span $E_s$ hence $\Phi_s=\Phi \cap E_s$:

$$\Phi_s=\{ \alpha_{n-d},\alpha_{n-2d},...,\alpha_{n-rd}\}$$

Then $\Phi_a=(\Phi_s)^\perp \cong A_{d-1}^{r}C_j$.\\

Next assume that $\mathcal{H}$ is $\dynkin[scale=1,labels={,,d,,,(r-1)d,,,n-1,n}, ply=2,rotate=0,fold radius = 2mm]{D}{*...*o*......*o*...*oo}$ where $n-rd=1$, $d=1$ or $2$ and $\Psi_0=A_{d-1}^r$. Hence:
$$\Delta_0= \Delta - (\{ \alpha_{n-d},\alpha_{n-2d},...,, \alpha_{n-(r-1)d}, \alpha_{1}, \alpha_s\})$$

We compute that $\dim E_s=r$. The following set of $\Pi$-stable roots of cardinality $r$ then must span $E_s$ hence $\Phi_s=\Phi \cap E_s$: 

$$\Phi_s=\{ \alpha_{n-d},\alpha_{n-2d},...,\alpha_{n-(r-1)d}, \alpha_*\} \text{ where 
 } \alpha_*=-(\alpha_1-\alpha_{s})=(0,2,...,2,1) $$
Therefore $\Phi_a=(\Phi_s)^\perp \cong A_{d-1}^{r}A_1$.\\

\newpage\newgeometry{top=0.5cm,bottom=0.5cm,left=0.5cm, right=0.5cm,foot=0.5cm}
\begin{landscape}
\begin{small}
\begin{longtable}{| c | c | c | c | c | c | c | c | c | c | c |}\caption{Isotropic maximal embeddings of abstract indices of type $C_n$}\label{C_n} \\ \hline

\multicolumn{1}{|c|}{} & \multicolumn{1}{c|}{} & \multicolumn{1}{c|}{} & \multicolumn{1}{c|}{} & \multicolumn{1}{c|}{} & \multicolumn{1}{c|}{} & \multicolumn{1}{c|}{} & \multicolumn{3}{c|}{\small{Special fields}} \\

$\Delta$ & \small{\!\!$\Stab_W(\Delta)$\!\!} & $\Pi$ & $\textcolor{white}{iiaaa} \mathcal{H} \textcolor{white}{aaaii}$ & $\Psi_0$ & $\Phi_a$ & $\mathcal{G}$ & \small{$\hspace{-0.8mm}\!\cd 1 \!\hspace{-1mm}$} & $\!\R\!$ & $\!\!Q_{\mathfrak{p}}\!\!$ \\ \hline \hline 
\endfirsthead
\hline &&&&&&&&&\\[-0.3cm]

\makecell{$C_{m}C_{n-m}$ \\ $\!1 \!\leq\! m \!\leq\! \frac{n}{2}\!$} & \makecell{$1$ \\ \!\!(if $m\!<\!\frac{n}{2}$)\!\!} 
& {$1$} & \makecell{$\dynkin[scale=0.8,labels={,,d,,,rd,,,m}, rotate=0]{C}{*...*o*......*o*...**}$ \\ $\dynkin[scale=0.8,labels={,,d',,,r'd',,,n-m}, rotate=0]{C}{*...*o*......*o*...**}$ \\ $j=m-rd \geq 0$, \\ $j'=n-m-r'd' \geq 0$, \\ $d,d',r+r' \geq 1$, \\ $d=2^a|2m$, $d'=2^{b}|2(n-m)$, \\ \!\!\!$j\hspace{-0.1mm}=\hspace{-0.1mm}0$ if $d\hspace{-0.1mm}=\hspace{-0.1mm}1$, $j'\hspace{-0.1mm}=\hspace{-0.1mm}0$ if $d'\hspace{-0.1mm}=\hspace{-0.1mm}1$,\!\!\! \\ if $m=n/2$ then $r \leq r'$} & \makecell{$A_{d-1}^rC_j$ \\ $\times$ \\ $A_{d'-1}^{r'}C_{j'}$} & \makecell{\!$A_{d-1}^rA_{d'-1}^{r'}$\! \\ $\times$ \\ $C_{j+j'}$}& \makecell{$\dynkin[scale=0.8,labels={,,d,,,id,,,}, rotate=0]{C}{*...*o*......*o*...**}$ \\ only if $d=d'$, \\ $(r+r')d \leq id \leq n$, \\ if $d=1$ then $i=n$, \\ if $j=0$ or $j'=0$ \\ then $i=r+r'$} & \makecell{$\checkmark$ \\ iff $d\!=\!1$} & \makecell{$\checkmark$ \\ iff $d\!=\!1$, or \\$d\!=\!2$ and \\ either \\ \!$i\!=\!r\!+\!r'$ or\! \\\!\!$\min\{t,n\!-\!t\}$;\!\! \\ where \\\!\!$t\!:=\!m\!+\!r'\!-\!r$\!\!} & \makecell{$\checkmark$ \\ if $d\!\in\!\{1,2\}$, \\ and either \\ $j=0$ or \\ $j'=0$, \\ $\xmark$ \\ if $d\!>\!2$} 
\topstrut \\
\hline &&&&&&&&&\\[-0.3cm]

\multirow{5}{*}{\makecell{$C_m^{n_m}$ \\ $\!mn_m=n\!$ \\ $n_m \geq 2$}} & \multirow{5}{*}{$S_{n_m}$} & \multirow{5}{*}{\makecell{primit. \\ \!subgrp.\!\\ of $S_{n_m}$}} &  \makecell{$\begin{tikzpicture}[scale=0.8,baseline=-5mm]
\dynkin[name=upper]{C}{*...*o*......*o*...***}
\node (current) at ($(upper root 1)+(0,-10mm)$) {};
\dynkin[at=(current),name=middle,labels={,,d,,,rd,,,,m}]{C}{*...*o*......*o*...***}
\begin{scope}[on background layer]
\foreach \x in {1,...,10} 
\draw[/Dynkin diagram/fold style]
($(upper root \x)$) -- ($(upper root \x)+(0,-3mm)$);
\foreach \x in {1,...,10} 
\draw[/Dynkin diagram/fold style]
($(middle root \x)+(0,3mm)$) -- ($(middle root \x)$);
\foreach \x in {1,...,10} 
\draw[dotted,black!40,line width=.5mm]
($(upper root \x)+(0,-3mm)$) -- ($(middle
root \x)+(0,3mm)$);
\end{scope}
\end{tikzpicture}$ \\ $j=m-rd \geq 1$, $r \geq 1$,\\ $d \geq 2$, $d=2^a|2m$} & $\!A_{d-1}^{rn_m}C_{j}^{n_m}\!$ & \!$A_{dn_m\!-1}^{r}C_{jn_m}\!$  & \makecell{$\dynkin[scale=0.8,labels={,,c,,,ic,,,}, rotate=0]{C}{*...*o*......*o*...**}$ \\ where $c|dn_m$,\\ $rdn_m \leq ic \leq n$, \\ if $c=1$ then $i=n$} & $\xmark$ & $\xmark$ & \makecell{$\xmark$ \\ unless $d\!=\!2$, \\ $j\!=\!1$ and \\ \!\!either \hspace{-0.2mm}$c\!=\!1$\hspace{-0.2mm} or\!\!\! \\ \!\!$(c,i)\!=\!(2,\hspace{-0.2mm}\lfloor\hspace{-0.1mm}\frac{n}{2}\hspace{-0.1mm}\rfloor)$\!\!\!} 
\topstrut \\ \cdashline{4-11}
&&&&&&&&&\\[-0.3cm]
&&& \makecell{$\begin{tikzpicture}[scale=0.8,baseline=-5mm] 
\dynkin[name=upper]{C}{*...*o*......*o*...**o}
\node (current) at ($(upper root 1)+(0,-10mm)$) {};
\dynkin[at=(current),name=middle,labels={,,d,,,(r-1)d,,,,m=rd}]{C}{*...*o*......*o*...**o}
\begin{scope}[on background layer]
\foreach \x in {1,...,10} 
\draw[/Dynkin diagram/fold style]
($(upper root \x)$) -- ($(upper root \x)+(0,-3mm)$);
\foreach \x in {1,...,10} 
\draw[/Dynkin diagram/fold style]
($(middle root \x)+(0,3mm)$) -- ($(middle root \x)$);
\foreach \x in {1,...,10} 
\draw[dotted,black!40,line width=.5mm]
($(upper root \x)+(0,-3mm)$) -- ($(middle
root \x)+(0,3mm)$);
\end{scope}
\end{tikzpicture}$ \\ $m=rd$, $d \geq 1$, $d=2^a$} & $A_{d-1}^{rn_m}$ & $A_{dn_m-1}^{r}$  & \makecell{$\dynkin[scale=0.8,labels={,,c,,,,,,ic}, rotate=0]{C}{*...*o*......*o*...*o}$ \\ where $c|dn_m$,\\ $n=ic$}  & \makecell{$\checkmark$ \\ iff \\ $\!c\!=\!d\!=\!1\!$} & \makecell{$\checkmark$ \\ iff $d\!=\!1$ \\ \!and $n_m\hspace{-0.5mm}\!=\!2$\!} & \makecell{$\checkmark$ \\ \!\!if $c\!=\!d\!=\!1$,\!\! \\ \!\!$n_m$ is prime,\hspace{-0.5mm}\!\!\\ \!\!$\xmark$ if \hspace{-0.2mm}$c$ \hspace{-0.2mm}or\hspace{-0.3mm} $d$ \hspace{-0.5mm}$>\!2$\!\!\!}
\topstrut \\ \hline

\makecell{$A_{n-1}$ \\ $p \neq 2$} & {$\Z_2$} & {$\Z_2$} &
\makecell{$\begin{tikzpicture}[scale=0.8,baseline=-5mm] 
\node (up) at (0,2) {};
\dynkin[at=(up),name=upper,labels={,,d,,,rd,,,,{n-rd},,,,,,n-1}, ply=2,fold radius=2mm]{A}{*.*o*.*o*.*.*.*o*.*o*.*};
\node (up1) at ($(up)+(4.5,1.5)$) {};
\node (up2) at ($(up1)+(0,-0.75)$) {};
\node (up3) at ($(up2)+(0,-1)$) {};
\node (up4) at ($(up3)+(0,-0.75)$) {};
\dynkin[at=(up1),name=upper_right, labels={,},ply=2,fold radius=2mm]{A}{.**.};
\draw[densely dotted] ($(upper_right root 1)+(-0.4,0)$) -- ($(upper_right root 1)+(-0.2,0)$);
\draw ($(upper_right root 1)+(-0.2,0)$) -- ($(upper_right root 1)+(0.01,0)$);
\draw[densely dotted] ($(upper_right root 2)+(-0.4,0)$) -- ($(upper_right root 2)+(-0.2,0)$);
\draw ($(upper_right root 2)+(-0.2,0)$) -- ($(upper_right root 2)+(0.01,0)$); 
\dynkin[at=(up2),name=middle1_right, labels={,}, ply=2,fold radius=2mm]{A}{.oo.};
\draw[densely dotted] ($(middle1_right root 1)+(-0.4,0)$) -- ($(middle1_right root 1)+(-0.2,0)$);
\draw ($(middle1_right root 1)+(-0.2,0)$) -- ($(middle1_right root 1)+(0.01,0)$);
\draw[densely dotted] ($(middle1_right root 2)+(-0.4,0)$) -- ($(middle1_right root 2)+(-0.2,0)$);
\draw ($(middle1_right root 2)+(-0.2,0)$) -- ($(middle1_right root 2)+(0.01,0)$);
\node[align=left] at ($(up2) + (0,-5mm)$) {\tiny{if $m=0$}};
\dynkin[at=(up3),name=middle_right, labels={,,}, ply=2,fold radius=2mm]{A}{.***.};
\draw[densely dotted] ($(middle_right root 1)+(-0.4,0)$) -- ($(middle_right root 1)+(-0.2,0)$);
\draw ($(middle_right root 1)+(-0.2,0)$) -- ($(middle_right root 1)+(0.01,0)$);
\draw[densely dotted] ($(middle_right root 3)+(-0.4,0)$) -- ($(middle_right root 3)+(-0.2,0)$);
\draw ($(middle_right root 3)+(-0.2,0)$) -- ($(middle_right root 3)+(0.01,0)$);
\dynkin[at=(up4),name=lower_right, labels={,,},ply=2,fold radius=2mm]{A}{.*o*.};
\draw[densely dotted] ($(lower_right root 1)+(-0.4,0)$) -- ($(lower_right root 1)+(-0.2,0)$);
\draw ($(lower_right root 1)+(-0.2,0)$) -- ($(lower_right root 1)+(0.01,0)$);
\draw[densely dotted] ($(lower_right root 3)+(-0.4,0)$) -- ($(lower_right root 3)+(-0.2,0)$);
\draw ($(lower_right root 3)+(-0.2,0)$) -- ($(lower_right root 3)+(0.01,0)$);
\draw [decorate, decoration = {calligraphic brace}] ($(lower_right root 3)+(-7mm,-5mm)$) --  ($(upper_right root 1)+(-7mm,2mm)$);
\node[align=left] at ($(up4) + (0,-5mm)$) {\tiny{if $m=-1$}}; \end{tikzpicture}$ \vspace{-12mm} \\ $\times$  ${\mathcal{T}_0}^1$  \\ $m=n-1-2rd \geq -1$, \\ $d,r \geq 1$, $d|n$} & \makecell{$A_{d-1}^{2r}A_{m}$ \\ ($A_{-1}:=\varnothing$)} & $A_{2d-1}^{r}C_{m+1}$ & \makecell{$\dynkin[scale=0.8,labels={,,c,,,ic,,,}, rotate=0]{C}{*...*o*......*o*...**}$ \\ where $c|2d$,\\ $2rd \leq ic \leq n$, \\ if $c=1$ then $i=n$} & \makecell{$\checkmark$ \\ iff $m\!\leq\!0$ \\ and \\ $c\!=\!d\!=\!1$} & \makecell{$\checkmark$ \\ iff $d\!=\!1$ \\\!\!and either\!\! \\ $c\!=\!1$ or \\\!\!$(c,i)\!=\!(2,r)$\!\!} & \makecell{$\checkmark$ \\ \!\!if $c\!=\!d\!=\!1$\!\! \\ and \\ \!\!$m\!\in\!\{0,-1\}$,\!\! \\ $\xmark$ \\ \!\!if $d\!\neq\!1$\!\! \\ \!\!or $c\!>\!2$\!\!} 
\topstrut \\ \hline
&&&&&&&&&\\[-0.4cm]
\multirow{12}{*}{\makecell{$\widetilde{D_n}$ \\ $p=2$}} & \multirow{12}{*}{$\Z_2$} & \multirow{4}{*}{$1$} & \makecell{$\dynkin[scale=1,labels={,,d,,,rd,,,,n}, rotate=0]{D}{*...*o*......*o*...***}$\\  $j=n-rd \geq 2$, \\ $d,r \geq 1$, $d=2^a|2n$} & $A_{d-1}^rD_j$ & $A_{d-1}^rC_j$ & \makecell{$\dynkin[scale=0.8,labels={,,d,,,id,,,}, rotate=0]{C}{*...*o*......*o*...**}$ \\ $rd \leq id \leq n$, \\ if $d=1$ then $i=n$} & $\xmark$ & $\xmark$ & $\xmark$ \topstrut \\ 
\cdashline{4-11} 
&&&&&&&&&\\[-0.4cm]
& & & \makecell{$\dynkin[scale=1,labels={,,d,,,(r-1)d,,,,rd},rotate=0]{D}{*...*o*......*o*...**o}$  \\ $n=rd$, $d=2^a \geq 1$} & $A_{d-1}^{r}$ & $A_{d-1}^{r}$ & \makecell{$\dynkin[scale=0.8,labels={,,d,,,,,,rd}, rotate=0]{C}{*...*o*......*o*...*o}$} & \makecell{$\checkmark$ \\ iff $d\!=\!1$} & $\xmark$ & $\xmark$ \topstrut \\ \cline{3-11} 
& & \multirow{5}{*}{$\Z_2$} & \makecell{$\dynkin[scale=1,labels={,,d,,,rd,,,,n}, ply=2,rotate=0,fold radius = 2mm]{D}{*...*o*......*o*...***}$ \\ $j=n-rd \geq 2$,\\ $d,r \geq 1$, $d=2^a|2n$} & $A_{d-1}^rD_j$ & $A_{d-1}^rC_j$  & \makecell{$\dynkin[scale=0.8,labels={,,d,,,id,,,}, rotate=0]{C}{*...*o*......*o*...**}$ \\ $rd \leq id \leq n$, \\ if $d=1$ then $i=n$}  & $\xmark$ & $\xmark$ & $\xmark$ \topstrut \\ \cdashline{4-11}
&&& \makecell{$\!\dynkin[scale=1,labels={,,d,,,(r-1)d,,,rd,n}, ply=2,rotate=0,fold radius = 2mm]{D}{*...*o*......*o*...*oo}\!$ \\ $n-rd=1$, $d \in \{1, 2\}$} & $A_{d-1}^r$ & $A_{d-1}^{r}A_1$ & \makecell{$\dynkin[scale=0.8,labels={,2,,4,,,2r,}, rotate=0]{C}{*o*o*...*o*}$ \\ only if $d=2$, \\ $\dynkin[scale=0.8, rotate=0]{C}{ooo...oo}$ \\ only if $d=1$} & \makecell{$\checkmark$ \\ iff $d\!=\!1$} & $\xmark$ & $\xmark$ \topstrut \\ \hline
\end{longtable}\restoregeometry
\end{small}
\end{landscape}
\newpage \restoregeometry
\subsection{\texorpdfstring{ \hspace*{-0.65cm}. $\Phi=D_n=\dynkin[labels={\alpha_1, \alpha_2,,\alpha_{n-1},\alpha_n},scale=1,
edge length=.5cm] D{**...***}$, $n \geq 4 $}{Dn}} \label{Dn}
\subsubsection{\texorpdfstring{$\Phi=D_4$}{D4}}
Let $\Phi=D_4$. We use the following standard parametrisation in $\mathbb{R}_4$ of $\Phi$. Let $\alpha_0=e_1+e_2$, $\alpha_1=e_1-e_2$, $\alpha_2=e_2-e_3$, $\alpha_3=e_3-e_4$, and $\alpha_4=e_3+e_4$ where $e_i$ denotes the $i$'th standard basis vector in $\R^4$.\\
Following \cite[Table II]{T2}, the isotropic indices of type $\Phi$ are given as follows:
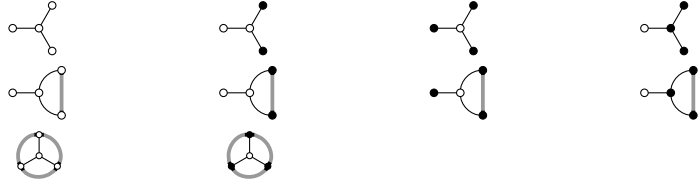
\begin{center}
\begin{tabular}{cp{1cm}cp{1cm}cp{1cm}c}
$\dynkin{D}{oooo}$ & & $\dynkin{D}{oo**}$ & & $\dynkin{D}{*o**}$ & & $\dynkin{D}{o***}$\\
$\dynkin[fold]{D}{oooo}$ & & $\dynkin[fold]D{oo**}$ & & $\dynkin[fold]{D}{*o**}$ & & $\dynkin[fold]{D}{o***}$\\
 \begin{tikzpicture}[scale=0.8,baseline=-2.2ex] \dynkin[ply=3,rotate=270,fold radius = 2mm]{D}{oooo} \end{tikzpicture} & & \begin{tikzpicture}[scale=0.8,baseline=-2.2ex] \dynkin[ply=3,rotate=270,fold radius = 2mm]{D}{*o**} \end{tikzpicture} & &   & &
\end{tabular}
\end{center}
Recall that $\Delta$ is one of the following types:
\begin{enumerate}
    \item $A_1^4$
    \item $A_3$ (of which there are 3 $W$-conjugacy classes)
    \item $A_2$.
\end{enumerate}

We start by considering $\Delta=A_1^{4}$ and choose a parametrisation of this subsystem by taking $\alpha_0, \alpha_1, \alpha_3,$ and $\alpha_4$. We have that $\Stab_{\Iso(\Phi)}(\Delta)=S_{4}$ by Table \ref{classlist} acting on the set of roots $\{\alpha_0, \alpha_1,\alpha_3,\alpha_4\}$. We denote the elements of $S_4$ by cycles in $\{1,2,3,4\}$ where each integer corresponds to a node. Here, we use $1$ to denote $\alpha_0$, $2$ to denote $\alpha_1$, $3$ to denote $\alpha_3$, and $4$ to denote $\alpha_4$. Up to permuting the nodes, there are four possibilities for the index $\mathcal{H}$ given in Table \ref{D_4} below. 
We start with the case of $\Pi=1$. Then $E_a=E_{\Delta_{0}}$, so the sets $\Psi_0$ and $\Phi_a$ are equal. \\
Let $\Pi=\Z_2^{(1)}=\langle (1,2)(3,4)\rangle$. Then we have the indices $\mathcal{H}=\begin{tikzpicture}[baseline=0.5ex] 
\dynkin[name=first]{A}{o}
\node (current) at ($(first root 1)+(-0,3.4mm)$) {};
\dynkin[at=(current),name=second]{A}{o}
\node (currenta) at ($(first root 1)+(-3.4mm,0)$) {};
\dynkin[at=(currenta),name=third]{A}{o}
\node (currenta) at ($(third root 1)+(-0,3.4mm)$) {};
\dynkin[at=(currenta),name=fourth]{A}{o}
\begin{scope}[on background layer]
\draw[/Dynkin diagram/fold style]
($(first root 1)$) -- ($(third root 1)$);%
\draw[/Dynkin diagram/fold style]
($(second root 1)$) -- ($(fourth root 1)$);%
\end{scope}
\end{tikzpicture}$ and  
$\mathcal{H}=\begin{tikzpicture}[baseline=0.5ex] 
\dynkin[name=first]{A}{o}
\node (current) at ($(first root 1)+(-0,3.4mm)$) {};
\dynkin[at=(current),name=second]{A}{*}
\node (currenta) at ($(first root 1)+(-3.4mm,0)$) {};
\dynkin[at=(currenta),name=third]{A}{o}
\node (currenta) at ($(third root 1)+(-0,3.4mm)$) {};
\dynkin[at=(currenta),name=fourth]{A}{*}
\begin{scope}[on background layer]
\draw[/Dynkin diagram/fold style]
($(first root 1)$) -- ($(third root 1)$);%
\draw[/Dynkin diagram/fold style]
($(second root 1)$) -- ($(fourth root 1)$);%
\end{scope}
\end{tikzpicture}$.\\
For $\mathcal{H}=\begin{tikzpicture}[baseline=0.5ex] 
\dynkin[name=first]{A}{o}
\node (current) at ($(first root 1)+(-0,3.4mm)$) {};
\dynkin[at=(current),name=second]{A}{o}
\node (currenta) at ($(first root 1)+(-3.4mm,0)$) {};
\dynkin[at=(currenta),name=third]{A}{o}
\node (currenta) at ($(third root 1)+(-0,3.4mm)$) {};
\dynkin[at=(currenta),name=fourth]{A}{o}
\begin{scope}[on background layer]
\draw[/Dynkin diagram/fold style]
($(first root 1)$) -- ($(third root 1)$);%
\draw[/Dynkin diagram/fold style]
($(second root 1)$) -- ($(fourth root 1)$);%
\end{scope}
\end{tikzpicture}$ we have $E_s=(E_{\Delta_0}^{\perp})^{\Pi}=\langle \alpha_0+\alpha_1, \alpha_3+\alpha_4\rangle=\langle e_1, e_3\rangle$, so $E_s^{\perp}=\langle e_2, e_4\rangle$, that is $\Phi_a=A_1^2$. In the case of $\mathcal{H}=\begin{tikzpicture}[baseline=0.5ex] 
\dynkin[name=first]{A}{o}
\node (current) at ($(first root 1)+(-0,3.4mm)$) {};
\dynkin[at=(current),name=second]{A}{*}
\node (currenta) at ($(first root 1)+(-3.4mm,0)$) {};
\dynkin[at=(currenta),name=third]{A}{o}
\node (currenta) at ($(third root 1)+(-0,3.4mm)$) {};
\dynkin[at=(currenta),name=fourth]{A}{*}
\begin{scope}[on background layer]
\draw[/Dynkin diagram/fold style]
($(first root 1)$) -- ($(third root 1)$);%
\draw[/Dynkin diagram/fold style]
($(second root 1)$) -- ($(fourth root 1)$);%
\end{scope}
\end{tikzpicture}$, we have $E_s=(E_{\Delta_0}^{\perp})^{\Pi}=\langle \alpha_3+\alpha_4\rangle=\langle e_3\rangle$, so $E_s^{\perp}=\langle e_1, e_2, e_4\rangle$, that is $\Phi_a=A_3$.\\

If $\Pi=\Z_2^{(2)}=\langle (1,2)\rangle$, we compute the following for each possible index $\mathcal{H}$:

\begin{table}[!htb]\begin{center}
\begin{tabular}{| c | c | c | K{3.5cm} | c | c |}    \hline
 $\mathcal{H}$ & $\Psi_0$ & $E_{\Delta_0}^{\perp}$ & $E_s$ & $E_a$ & $\Phi_a$ \\ \hline \hline 
$\!\!\!\begin{tikzpicture}[baseline=1ex] 
\dynkin[name=first]{A}{o}
\node (current) at ($(first root 1)+(-0,3.4mm)$) {};
\dynkin[at=(current),name=second]{A}{o}
\node (currenta) at ($(first root 1)+(-3.4mm,0)$) {};
\dynkin[at=(currenta),name=third]{A}{o}
\node (currenta) at ($(third root 1)+(-0,3.4mm)$) {};
\dynkin[at=(currenta),name=fourth]{A}{o}
\begin{scope}[on background layer]
\draw[/Dynkin diagram/fold style]
($(second root 1)$) -- ($(fourth root 1)$);%
\end{scope}
\end{tikzpicture}\!\!\!$ & $\varnothing$ & $\R^4$ & $\langle \alpha_0+\alpha_1, \alpha_3, \alpha_4\rangle =\langle e_2, e_3, e_4\rangle$ & $\langle e_1 \rangle$ & $\varnothing$ \\ \hline
$\!\!\!\begin{tikzpicture}[baseline=1ex] 
\dynkin[name=first]{A}{*}
\node (current) at ($(first root 1)+(-0,3.4mm)$) {};
\dynkin[at=(current),name=second]{A}{o}
\node (currenta) at ($(first root 1)+(-3.4mm,0)$) {};
\dynkin[at=(currenta),name=third]{A}{o}
\node (currenta) at ($(third root 1)+(-0,3.4mm)$) {};
\dynkin[at=(currenta),name=fourth]{A}{o}
\begin{scope}[on background layer]
\draw[/Dynkin diagram/fold style]
($(second root 1)$) -- ($(fourth root 1)$);%
\end{scope}
\end{tikzpicture}\!\!\!$ & $A_1$ & $\langle \alpha_0,\alpha_1,\alpha_3\rangle$ & $\langle \alpha_0+\alpha_1, \alpha_3,\rangle =\langle e_1, e_3-e_4\rangle$ & $\langle e_2,e_3+e_4 \rangle$ & $A_1$ \\ \hline
$\!\!\!\begin{tikzpicture}[baseline=1ex] 
\dynkin[name=first]{A}{o}
\node (current) at ($(first root 1)+(-0,3.4mm)$) {};
\dynkin[at=(current),name=second]{A}{*}
\node (currenta) at ($(first root 1)+(-3.4mm,0)$) {};
\dynkin[at=(currenta),name=third]{A}{o}
\node (currenta) at ($(third root 1)+(-0,3.4mm)$) {};
\dynkin[at=(currenta),name=fourth]{A}{*}
\begin{scope}[on background layer]
\draw[/Dynkin diagram/fold style]
($(second root 1)$) -- ($(fourth root 1)$);%
\end{scope}
\end{tikzpicture}\!\!\!$ & $A_1^2$ & $\langle \alpha_3, \alpha_4 \rangle = \langle e_3, e_4 \rangle$ & $\langle \alpha_3, \alpha_4 \rangle = \langle e_3, e_4 \rangle$ & $\langle \alpha_0, \alpha_1 \rangle = \langle e_1, e_2 \rangle$ & $A_1^2$ \\ \hline
$\!\!\!\begin{tikzpicture}[baseline=1ex] 
\dynkin[name=first]{A}{o}
\node (current) at ($(first root 1)+(-0,3.4mm)$) {};
\dynkin[at=(current),name=second]{A}{*}
\node (currenta) at ($(first root 1)+(-3.4mm,0)$) {};
\dynkin[at=(currenta),name=third]{A}{o}
\node (currenta) at ($(third root 1)+(-0,3.4mm)$) {};
\dynkin[at=(currenta),name=fourth]{A}{*}
\begin{scope}[on background layer]
\draw[/Dynkin diagram/fold style]
($(first root 1)$) -- ($(third root 1)$);%
\end{scope}
\end{tikzpicture}\!\!\!$ & $A_1^2$ & $\langle \alpha_3, \alpha_4 \rangle = \langle e_3, e_4 \rangle$ & $\langle \alpha_3 + \alpha_4 \rangle = \langle e_3\rangle$ & $\langle e_1, e_2, e_4 \rangle$ & $A_3$ \\ \hline
$\!\!\!\begin{tikzpicture}[baseline=1ex] 
\dynkin[name=first]{A}{*}
\node (current) at ($(first root 1)+(-0,3.4mm)$) {};
\dynkin[at=(current),name=second]{A}{*}
\node (currenta) at ($(first root 1)+(-3.4mm,0)$) {};
\dynkin[at=(currenta),name=third]{A}{o}
\node (currenta) at ($(third root 1)+(-0,3.4mm)$) {};
\dynkin[at=(currenta),name=fourth]{A}{*}
\begin{scope}[on background layer]
\draw[/Dynkin diagram/fold style]
($(second root 1)$) -- ($(fourth root 1)$);%
\end{scope}
\end{tikzpicture}\!\!\!$ & $A_1^3$ & $\langle \alpha_3 \rangle $ & $\langle \alpha_3\rangle$ & $\langle \alpha_0, \alpha_1, \alpha_4 \rangle$ & $A_1^3$ \\ \hline
\end{tabular}\end{center}\end{table}
  \newpage
We have two conjugacy classes of $\Z_2^2$ in $S_4$, where $(\Z_2^2)^{(1)}$ acts on the simple roots of $A_1^4$ in two orbits and $(\Z_2^2)^{(2)}$ gives rise to one orbit. Observe that the subgroups $\Z_4$ and $\Di_4$ have the same orbit as $(\Z_2^2)^{(2)}$ acting on $A_1^4$ and on $D_4$, but the subgroups $\Alt_4$ and $S_4$ only act on $A_1^4$ with the same orbits as $(\Z_2^2)^{(2)}$. For each set of subgroups we again take a closer look at the computations of $\Phi_a$. First, let $\Pi=(\Z_2^2)^{(1)}$.
  
\begin{table}[H] \begin{center}
\begin{tabular}{| c | c | c | K{3.2cm} | c | c |}    \hline
 $\mathcal{H}$ & $\Psi_0$ & $E_{\Delta_0}^{\perp}$ & $E_s$ & $E_a$ & $\Phi_a$ \\ \hline \hline 
$\!\!\!\begin{tikzpicture}[baseline=1ex] 
\dynkin[name=first]{A}{o}
\node (current) at ($(first root 1)+(-0,3.4mm)$) {};
\dynkin[at=(current),name=second]{A}{o}
\node (currenta) at ($(first root 1)+(-3.4mm,0)$) {};
\dynkin[at=(currenta),name=third]{A}{o}
\node (currenta) at ($(third root 1)+(-0,3.4mm)$) {};
\dynkin[at=(currenta),name=fourth]{A}{o}
\begin{scope}[on background layer]
\draw[/Dynkin diagram/fold style]
($(second root 1)$) -- ($(fourth root 1)$);%
\draw[/Dynkin diagram/fold style]
($(first root 1)$) -- ($(third root 1)$);%
\end{scope}
\end{tikzpicture}\!\!\!$ & $\varnothing$ & $\R^4$ & $\langle \alpha_0+\alpha_1, \alpha_3+\alpha_4\rangle =\langle e_2, e_3\rangle$ & $\langle e_1, e_4 \rangle$ & $A_1^2$ \\ \hline
$\!\!\!\begin{tikzpicture}[baseline=1ex] 
\dynkin[name=first]{A}{o}
\node (current) at ($(first root 1)+(-0,3.4mm)$) {};
\dynkin[at=(current),name=second]{A}{*}
\node (currenta) at ($(first root 1)+(-3.4mm,0)$) {};
\dynkin[at=(currenta),name=third]{A}{o}
\node (currenta) at ($(third root 1)+(-0,3.4mm)$) {};
\dynkin[at=(currenta),name=fourth]{A}{*}
\begin{scope}[on background layer]
\draw[/Dynkin diagram/fold style]
($(second root 1)$) -- ($(fourth root 1)$);%
\draw[/Dynkin diagram/fold style]
($(first root 1)$) -- ($(third root 1)$);%
\end{scope}
\end{tikzpicture}\!\!\!$ & $A_1^2$ & $\langle \alpha_3, \alpha_4 \rangle = \langle e_3, e_4 \rangle$ & $\langle \alpha_3+ \alpha_4 \rangle = \langle e_3 \rangle$ & $\langle \alpha_0, \alpha_1, \alpha_2+\alpha_3 \rangle = \langle e_1, e_2, e_4 \rangle$ & $A_3$ \\ \hline
\end{tabular} 
\end{center}
\end{table}
  
The action of the subgroups $\Pi=(\Z_2^2)^{(2)}, \Z_4, \Di_8, \Alt_4,$ and $S_4$ on $A_1^4$ all have the same orbits, so we describe them using the index $\mathcal{H}= \begin{tikzpicture}[baseline=0.5ex] 
\dynkin[name=first]{A}{o}
\node (current) at ($(first root 1)+(-0,3.4mm)$) {};
\dynkin[at=(current),name=second]{A}{o}
\node (currenta) at ($(first root 1)+(-3.4mm,0)$) {};
\dynkin[at=(currenta),name=third]{A}{o}
\node (currenta) at ($(third root 1)+(-0,3.4mm)$) {};
\dynkin[at=(currenta),name=fourth]{A}{o}
\begin{scope}[on background layer]
\draw[/Dynkin diagram/fold style]
($(first root 1)$) -- ($(second root 1)$);%
\draw[/Dynkin diagram/fold style]
($(second root 1)$) -- ($(fourth root 1)$);%
\draw[/Dynkin diagram/fold style]
($(third root 1)$) -- ($(fourth root 1)$);%
\draw[/Dynkin diagram/fold style]
($(third root 1)$) -- ($(first root 1)$);%
\end{scope}
\end{tikzpicture}$. Clearly, $\Psi_0=\varnothing$, and so $E_{\Delta_0}^{\perp}=\R^4$. Then $E_s=\langle \alpha_0+\alpha_1+\alpha_3+\alpha_4\rangle = \langle e_1+e_3 \rangle$ which gives us $E_a = \langle e_2, e_4, e_1-e_3\rangle$ and therefore $\Phi_a = A_1^3$.\\
Finally, consider $\Pi=\Z_3$ or $S_3$. Here we get the following results by embedding into $\R^4$:

  \begin{table}[!htb]\begin{center}
\begin{tabular}{| c | c | c | K{3.5cm} | c | c |}    \hline
 $\mathcal{H}$ & $\Psi_0$ & $E_{\Delta_0}^{\perp}$ & $E_s$ & $E_a$ & $\Phi_a$ \\ \hline \hline 
 
$\!\!\!\begin{tikzpicture}[baseline=1ex] 
\dynkin[name=first]{A}{o}

\node (current) at ($(first root 1)+(-0,3.4mm)$) {};
\dynkin[at=(current),name=second]{A}{o}

\node (currenta) at ($(first root 1)+(-3.4mm,0)$) {};
\dynkin[at=(currenta),name=third]{A}{o}
\node (currenta) at ($(third root 1)+(-0,3.4mm)$) {};
\dynkin[at=(currenta),name=fourth]{A}{o}

\begin{scope}[on background layer]

\draw[/Dynkin diagram/fold style]
($(third root 1)$) -- ($(second root 1)$);%
\draw[/Dynkin diagram/fold style]
($(second root 1)$) -- ($(fourth root 1)$);%
\draw[/Dynkin diagram/fold style]
($(third root 1)$) -- ($(fourth root 1)$);%

\end{scope}
\end{tikzpicture}\!\!\!$ & $\varnothing$ & $\R^4$ & $\langle \alpha_0+\alpha_1+\alpha_3, \alpha_4\rangle =\langle e_1+e_3, e_3+e_4\rangle$ & $\langle e_2, e_1-e_3+e_4 \rangle$ & $\varnothing$ \\ \hline
$\!\!\!\begin{tikzpicture}[baseline=1ex] 
\dynkin[name=first]{A}{o}

\node (current) at ($(first root 1)+(-0,3.4mm)$) {};
\dynkin[at=(current),name=second]{A}{*}

\node (currenta) at ($(first root 1)+(-3.4mm,0)$) {};
\dynkin[at=(currenta),name=third]{A}{*}
\node (currenta) at ($(third root 1)+(-0,3.4mm)$) {};
\dynkin[at=(currenta),name=fourth]{A}{*}

\begin{scope}[on background layer]

\draw[/Dynkin diagram/fold style]
($(third root 1)$) -- ($(second root 1)$);%
\draw[/Dynkin diagram/fold style]
($(second root 1)$) -- ($(fourth root 1)$);%
\draw[/Dynkin diagram/fold style]
($(third root 1)$) -- ($(fourth root 1)$);%

\end{scope}
\end{tikzpicture}\!\!\!$ & $A_1^3$ & $\langle \alpha_4 \rangle $ & $\langle \alpha_4 \rangle = \langle e_3+e_4 \rangle$ & $\langle e_1, e_2, e_3 - e_4\rangle = \langle \alpha_0, \alpha_1, \alpha_4\rangle$ & $A_1^3$ \\ \hline
$\!\!\!\begin{tikzpicture}[baseline=1ex] 
\dynkin[name=first]{A}{*}

\node (current) at ($(first root 1)+(-0,3.4mm)$) {};
\dynkin[at=(current),name=second]{A}{o}

\node (currenta) at ($(first root 1)+(-3.4mm,0)$) {};
\dynkin[at=(currenta),name=third]{A}{o}
\node (currenta) at ($(third root 1)+(-0,3.4mm)$) {};
\dynkin[at=(currenta),name=fourth]{A}{o}

\begin{scope}[on background layer]

\draw[/Dynkin diagram/fold style]
($(third root 1)$) -- ($(second root 1)$);%
\draw[/Dynkin diagram/fold style]
($(second root 1)$) -- ($(fourth root 1)$);%
\draw[/Dynkin diagram/fold style]
($(third root 1)$) -- ($(fourth root 1)$);%

\end{scope}
\end{tikzpicture}\!\!\!$ & $A_1$ & $\langle \alpha_0, \alpha_1, \alpha_3 \rangle$ & $\langle \alpha_0+\alpha_1+\alpha_3 \rangle= \langle 2e_1+e_3-e_4\rangle$ & $\langle e_2, e_1+2e_4, e_3+e_4 \rangle$ & $A_1$ \\ \hline
\end{tabular}\end{center}\end{table}
  
In the case of $\Delta=A_3$, we have that $\Stab_{\Iso(\Phi)}(\Delta) = \Z_2^2$. The subsystem $A_3$ can be embedded into $D_4$ by parametrising it as $\langle e_2-e_3, e_3-e_4,e_3+e_4\rangle= \langle \alpha_2, \alpha_3, \alpha_4\rangle$. In this case we can write the generators of $\Z_2^2$ as the map
$$\varphi: e_1\mapsto -e_1, e_2\mapsto e_2, e_3\mapsto e_3, e_4\mapsto -e_4$$ and the map 
$$\psi: e_1\mapsto -e_1, e_2\mapsto e_2, e_3\mapsto e_3, e_4\mapsto e_4.$$

We may ignore the cases where $\Pi$ acts trivially on the 1-dimensional space $\Delta^{\perp}=\langle e_1 \rangle$, as any associated index $\mathcal{H}$ is not maximal in $\Phi$. The remaining possibilities are when $\Pi= \Stab_{\Iso(\Phi)}(\Delta) \cong \Z_2^2$, and the order 2 subgroups $\langle \varphi \rangle=:(\Z_2)^{(1)}$ and $\langle \psi \rangle=:(\Z_2)^{(2)}$. The action of $\Z_2^2$ and $(\Z_2)^{(1)}$ on the index of $A_3$ is the same, so we combine both of these cases for the computation of $\Phi_a$. However, the action on the index of $D_4$ is different for each group.

Consider $\Pi=(\Z_2)^{(1)}$ and $\Z_2^2$. In both cases we compute $\Phi_a$ as follows.
\begin{table}[h]\begin{center}
\begin{tabular}{| c | c | c | K{3.5cm} | c | c |}    \hline
 $\mathcal{H}$ & $\Psi_0$ & $E_{\Delta_0}^{\perp}$ & $E_s$ & $E_a$ & $\Phi_a$ \\ \hline \hline 
 
$\begin{tikzpicture}[baseline=-0.3ex]
\dynkin[ply=2,rotate=90,fold radius = 2mm]{A}{ooo}\end{tikzpicture}\times\mathcal{T}^1_0$ & $\varnothing$ & $\R^4$ & $\langle e_2, e_3\rangle$ & $\langle e_1, e_4\rangle= \langle e_1-e_4, e_1+e_4 \rangle$ & $A_1^2$ \\ \hline
$\begin{tikzpicture}[baseline=-0.3ex]
\dynkin[ply=2,rotate=90,fold radius = 2mm]{A}{*o*}\end{tikzpicture}\times\mathcal{T}^1_0$  & $A_1^2$ & $\langle e_1, e_2 \rangle $ & $\langle e_2 \rangle$ & $\langle e_1, e_3, e_4 \rangle$ & $A_3$ \\ \hline
$\begin{tikzpicture}[baseline=-0.3ex]
\dynkin[ply=2,rotate=90,fold radius = 2mm]{A}{o*o}\end{tikzpicture}\times\mathcal{T}^1_0$ & $A_1$ & $\langle e_1, e_4, e_2+e_3 \rangle$ & $\langle e_2+e_3 \rangle$ & $\langle e_1, e_4, e_2-e_3 \rangle$ & $A_1^3$ \\ \hline
\end{tabular}\end{center}\end{table}
For $\Pi=(\Z_2)^{(2)}$ we have the following results.

\begin{table}[H]\begin{center}
\begin{tabular}{| c | c | c | K{4cm} | c | c |}    \hline
 $\mathcal{H}$ & $\Psi_0$ & $E_{\Delta_0}^{\perp}$ & $E_s$ & $E_a$ & $\Phi_a$ \\ \hline \hline 
 
$\begin{tikzpicture}[baseline=-0.3ex]
\dynkin{A}{ooo}\end{tikzpicture}\times\mathcal{T}^1_0$ & $\varnothing$ & $\R^4$ & $\langle e_2, e_3, e_4\rangle$ & $\langle e_1\rangle$ & $\varnothing$ \\ \hline
$\begin{tikzpicture}[baseline=-0.3ex]
\dynkin{A}{*o*}\end{tikzpicture}\times\mathcal{T}^1_0$  & $A_1^2$ & $\langle e_1, e_2 \rangle $ & $\langle e_2 \rangle$ & $\langle e_1, e_3, e_4 \rangle$ & $A_3$ \\ \hline
\end{tabular}\end{center}\end{table}
We summarise the results in the following table.
\newpage\newgeometry{top=1cm,bottom=1cm,left=0.5cm, right=0.5cm,foot=0cm}
\begin{small}
\begin{longtable}{| c | c | c | c | c | c | c | c | c | c | c |}\caption{Isotropic maximal embeddings of abstract indices of type $D_4$}\label{D_4} \\ \hline

\multicolumn{1}{|c|}{} & \multicolumn{1}{c|}{} & \multicolumn{1}{c|}{} & \multicolumn{1}{c|}{} & \multicolumn{1}{c|}{} & \multicolumn{1}{c|}{} & \multicolumn{1}{c|}{} & \multicolumn{3}{c|}{\small{Special fields}} \\

  $\Delta$ & \small{$\hspace{-0.5mm}\Stab_{\operatorname{Iso}(\Phi)}(\Delta)\hspace{-0.5mm}$} & $\Pi$ & $\textcolor{white}{iiaaa} \mathcal{H} \textcolor{white}{aaaii}$ & $\Psi_0$ & $\Phi_a$ & $\mathcal{G}$ & \small{$\hspace{-0.8mm}\!\cd 1 \!\hspace{-1mm}$} & $\!\R\!$ & $\!\!Q_{\mathfrak{p}}\!\!$ \\ \hline \hline 
\endfirsthead

\multicolumn{11}{|c|}%
{{\bfseries \tablename\ \thetable{} -- continued from previous page}} \\
\hline
\multicolumn{1}{|c|}{} & \multicolumn{1}{c|}{} & \multicolumn{1}{c|}{} & \multicolumn{1}{c|}{} & \multicolumn{1}{c|}{} & \multicolumn{1}{c|}{} & \multicolumn{1}{c|}{} & \multicolumn{3}{c|}{\small{Special fields}} \\

  $\Delta$ & \small{$\hspace{-0.5mm}\Stab_{\Iso(\Phi)}(\Delta)\hspace{-0.5mm}$} & $\Pi$ & $\mathcal{H}$ & $\Psi_0$ & $\Phi_a$ & $\mathcal{G}$ & \small{$\hspace{-0.8mm}\!\cd 1 \!\hspace{-1mm}$} & $\!\R\!$ & $\!\!Q_{\mathfrak{p}}\!\!$  \\ \hline \hline \endhead
\hline \multicolumn{11}{|c|}{{Continued on next page}} \\ \hline \endfoot 
 \endlastfoot

    \multirow{46}{*}{$A_1^4$} & \multirow{46}{*}{$S_4$} & \multirow{6}{*}{$1$} &  $\left(\dynkin{A}{o}\right)^4$ & $\varnothing$  & $\varnothing$ & $\dynkin{D}{oooo}$ & $\checkmark$ & $\checkmark$ & $\checkmark$ \topstrut  \\ \cdashline{4-11} %
&&& $\dynkin{A}{*}\hspace{-0.3mm}\times\hspace{-0.3mm} \left(\dynkin{A}{o}\right)^3$ & $A_1$ & \textcolor{white}{\Big\{}$A_1$\textcolor{white}{$\Big\{$} &&&& \topstrut \\ \cdashline{4-11}
&&& $\left(\dynkin{A}{*}\right)^2\hspace{-0.3mm}\times\hspace{-0.3mm} \left(\dynkin{A}{o}\right)^2$ & $A_1^2$ & $A_1^2$ & $\dynkin{D}{oo**}$\ & $\xmark$ & $\checkmark$ & $\checkmark$ \topstrut \\ \cdashline{4-11} 
&&& $\left(\dynkin{A}{*}\right)^3\hspace{-0.3mm}\times\hspace{-0.3mm} \dynkin{A}{o}$ & $A_1^3$ & $A_1^3$ & $\dynkin{D}{*o**}$ & $\xmark$ & $\xmark$ & $\xmark$ 
\topstrut \\ \cline{3-11} 
&& \multirow{9}{*}{\makecell{$\Z_2^{(1)}=\langle (1,2)(3,4)\rangle$ \\ 3 $\Stab_W(\Delta)$-c.c.}} &  \multirow{2.6}{*}{$\!\!\!\begin{tikzpicture}[baseline=1ex] 
\dynkin[name=first]{A}{o}
\node (current) at ($(first root 1)+(-0,3.4mm)$) {};
\dynkin[at=(current),name=second]{A}{o}
\node (currenta) at ($(first root 1)+(-3.4mm,0)$) {};
\dynkin[at=(currenta),name=third]{A}{o}
\node (currenta) at ($(third root 1)+(-0,3.4mm)$) {};
\dynkin[at=(currenta),name=fourth]{A}{o}
\begin{scope}[on background layer]
\draw[/Dynkin diagram/fold style]
($(first root 1)$) -- ($(third root 1)$);%
\draw[/Dynkin diagram/fold style]
($(second root 1)$) -- ($(fourth root 1)$);%
\end{scope}\end{tikzpicture}\!\!\!$} & \multirow{2.6}{*}{$\varnothing$} & \multirow{2.6}{*}{$A_1^2$} &  $\dynkin{D}{oooo}$  
& $\checkmark$ & $\checkmark$ & $\checkmark$ \topstrut  \\
&& &   &  &  &   $\dynkin{D}{oo**}$ & $\xmark$ & $\checkmark$ & $?$ 
\topstrut  \\ \cdashline{4-11}
&&& 
\multirow{5}{*}{\makecell{$\!\!\!\begin{tikzpicture}[baseline=1ex] 
\dynkin[name=first]{A}{o}
\node (current) at ($(first root 1)+(-0,3.4mm)$) {};
\dynkin[at=(current),name=second]{A}{*}
\node (currenta) at ($(first root 1)+(-3.4mm,0)$) {};
\dynkin[at=(currenta),name=third]{A}{o}
\node (currenta) at ($(third root 1)+(-0,3.4mm)$) {};
\dynkin[at=(currenta),name=fourth]{A}{*}
\begin{scope}[on background layer]
\draw[/Dynkin diagram/fold style]
($(first root 1)$) -- ($(third root 1)$);%
\draw[/Dynkin diagram/fold style]
($(second root 1)$) -- ($(fourth root 1)$);%
\end{scope}
\end{tikzpicture}\!\!\!$ 
}} & \multirow{5}{*}{$A_1^2$} & \multirow{5}{*}{$A_3$ } & $\dynkin {D}{oooo}$ 
& $\xmark$ & $\xmark$ & $?$ \topstrut \\ 
&& &   &  &  &  $\dynkin {D}{oo**}$ & $\xmark$ & $\xmark$ & $?$ \topstrut  \\ 
&& &   &  &  &  $\dynkin {D}{o***}$ & $\xmark$ & $\xmark$ & $\xmark$ \topstrut  \\ \cline{3-11}
&& \multirow{12}{*}{\makecell{$\Z_2^{(2)}=\langle (1,2)\rangle$ \\ 3 $\Stab_W(\Delta)$-c.c.}} &  $\!\!\!\begin{tikzpicture}[baseline=1ex] 
\dynkin[name=first]{A}{o}

\node (current) at ($(first root 1)+(-0,3.4mm)$) {};
\dynkin[at=(current),name=second]{A}{o}

\node (currenta) at ($(first root 1)+(-3.4mm,0)$) {};
\dynkin[at=(currenta),name=third]{A}{o}
\node (currenta) at ($(third root 1)+(-0,3.4mm)$) {};
\dynkin[at=(currenta),name=fourth]{A}{o}

\begin{scope}[on background layer]

\draw[/Dynkin diagram/fold style]
($(second root 1)$) -- ($(fourth root 1)$);%

\end{scope}
\end{tikzpicture}\!\!\!$  & $\varnothing$  & $\varnothing$ &  $\dynkin [fold]D{oooo}$ & $\checkmark$ & $\checkmark$ & $\checkmark$ \topstrut  \\ \cdashline{4-11} 

&&& $\!\!\!\begin{tikzpicture}[baseline=1ex] 
\dynkin[name=first]{A}{*}
\node (current) at ($(first root 1)+(-0,3.4mm)$) {};
\dynkin[at=(current),name=second]{A}{o}
\node (currenta) at ($(first root 1)+(-3.4mm,0)$) {};
\dynkin[at=(currenta),name=third]{A}{o}
\node (currenta) at ($(third root 1)+(-0,3.4mm)$) {};
\dynkin[at=(currenta),name=fourth]{A}{o}
\begin{scope}[on background layer]
\draw[/Dynkin diagram/fold style]
($(second root 1)$) -- ($(fourth root 1)$);%
\end{scope}
\end{tikzpicture}\!\!\!$ & $A_1$ & $A_1$ &  & &  &  \topstrut  \\ \cdashline{4-11}

&&& $\!\!\!\begin{tikzpicture}[baseline=1ex] 
\dynkin[name=first]{A}{o}
\node (current) at ($(first root 1)+(-0,3.4mm)$) {};
\dynkin[at=(current),name=second]{A}{*}
\node (currenta) at ($(first root 1)+(-3.4mm,0)$) {};
\dynkin[at=(currenta),name=third]{A}{o}
\node (currenta) at ($(third root 1)+(-0,3.4mm)$) {};
\dynkin[at=(currenta),name=fourth]{A}{*}
\begin{scope}[on background layer]
\draw[/Dynkin diagram/fold style]
($(second root 1)$) -- ($(fourth root 1)$);%
\end{scope}
\end{tikzpicture}\!\!\!$ & $A_1^2$ & $A_1^2$ & $\dynkin[fold]{D}{oo**}$ & $\xmark$ & $\xmark$ & $\xmark$ \topstrut \\ \cdashline{4-11} 

&&& \multirow{5}{*}{$\!\!\!\begin{tikzpicture}[baseline=1ex] 
\dynkin[name=first]{A}{o}
\node (current) at ($(first root 1)+(-0,3.4mm)$) {};
\dynkin[at=(current),name=second]{A}{*}
\node (currenta) at ($(first root 1)+(-3.4mm,0)$) {};
\dynkin[at=(currenta),name=third]{A}{o}
\node (currenta) at ($(third root 1)+(-0,3.4mm)$) {};
\dynkin[at=(currenta),name=fourth]{A}{*}
\begin{scope}[on background layer]
\draw[/Dynkin diagram/fold style]
($(first root 1)$) -- ($(third root 1)$);%
\end{scope}
\end{tikzpicture}\!\!\!$} & \multirow{5}{*}{$A_1^2$} & \multirow{5}{*}{$A_3$} & $\dynkin [fold]D{oooo}$ &$\xmark$&$\checkmark$&$?$ \topstrut \\
&&&&&& $\dynkin [fold]D{oo**}$ &$\xmark$&$\xmark$&$\xmark$ \topstrut \\
&&&&&& $\dynkin [fold]D{o***}$ &$\xmark$&$\checkmark$&$\xmark$ \topstrut \\ \cdashline{4-11} 
&&& $\!\!\!\begin{tikzpicture}[baseline=1ex] 
\dynkin[name=first]{A}{*}

\node (current) at ($(first root 1)+(-0,3.4mm)$) {};
\dynkin[at=(current),name=second]{A}{*}

\node (currenta) at ($(first root 1)+(-3.4mm,0)$) {};
\dynkin[at=(currenta),name=third]{A}{o}
\node (currenta) at ($(third root 1)+(-0,3.4mm)$) {};
\dynkin[at=(currenta),name=fourth]{A}{*}

\begin{scope}[on background layer]

\draw[/Dynkin diagram/fold style]
($(second root 1)$) -- ($(fourth root 1)$);%

\end{scope}
\end{tikzpicture}\!\!\!$ & $A_1^3$ & $A_1^3$ & $\dynkin [fold]D{*o**}$ & $\xmark$ & $\xmark$ & $\checkmark$ \topstrut \\ \cline{3-11} 

&& \multirow{9}{*}{\makecell{\!$(\Z_2^2)^{(1)}\hspace{-0.3mm}=\hspace{-0.3mm}\langle (1,2), (3,4)\rangle$\! \\ 3 $\Stab_W(\Delta)$-c.c.}} &  \multirow{3}{*}{$\!\!\!\begin{tikzpicture}[baseline=1ex] 
\dynkin[name=first]{A}{o}
\node (current) at ($(first root 1)+(-0,3.4mm)$) {};
\dynkin[at=(current),name=second]{A}{o}
\node (currenta) at ($(first root 1)+(-3.4mm,0)$) {};
\dynkin[at=(currenta),name=third]{A}{o}
\node (currenta) at ($(third root 1)+(-0,3.4mm)$) {};
\dynkin[at=(currenta),name=fourth]{A}{o}
\begin{scope}[on background layer]
\draw[/Dynkin diagram/fold style]
($(second root 1)$) -- ($(fourth root 1)$);%
\draw[/Dynkin diagram/fold style]
($(third root 1)$) -- ($(first root 1)$);%
\end{scope}
\end{tikzpicture}\!\!\!$}  & \multirow{3}{*}{$\varnothing$}  & \multirow{3}{*}{$A_1^2$} &  $\dynkin [fold]D{oooo}$ & $\checkmark$ & $\xmark$ & $\checkmark$ \topstrut \\
&& & & & &  $\dynkin [fold]D{oo**}$ & $\xmark$ & $\xmark$ & $\xmark$ \topstrut  \\
 \cdashline{4-11}
&& & \multirow{5}{*}{$\!\!\!\begin{tikzpicture}[baseline=1ex] 
\dynkin[name=first]{A}{o}
\node (current) at ($(first root 1)+(-0,3.4mm)$) {};
\dynkin[at=(current),name=second]{A}{*}
\node (currenta) at ($(first root 1)+(-3.4mm,0)$) {};
\dynkin[at=(currenta),name=third]{A}{o}
\node (currenta) at ($(third root 1)+(-0,3.4mm)$) {};
\dynkin[at=(currenta),name=fourth]{A}{*}
\begin{scope}[on background layer]
\draw[/Dynkin diagram/fold style]
($(second root 1)$) -- ($(fourth root 1)$);%
\draw[/Dynkin diagram/fold style]
($(third root 1)$) -- ($(first root 1)$);%
\end{scope}
\end{tikzpicture}\!\!\!$}  & \multirow{5}{*}{$A_1^2$}  & \multirow{5}{*}{$A_3$} &  $\dynkin [fold]D{oooo}$ & $\xmark$&$\xmark$&$?$  \topstrut  \\
&& & & & &  $\dynkin [fold]D{oo**}$ & $\xmark$&$\xmark$&$\xmark$ \topstrut  \\
&& & & & &  $\dynkin [fold]D{o***}$ & $\xmark$&$\xmark$&$\xmark$ \topstrut  \\ \cline{3-11}
&& \multirow{5}{*}{\makecell{~\\[-0.3cm] \!$(\Z_2^2)^{(2)}\hspace{-0.3mm}=\hspace{-0.3mm}\langle (1,2)(3,4),$\! \\$(1,4)(2,3)\rangle$}} &  \multirow{5}{*}{$\!\!\!\begin{tikzpicture}[baseline=1ex] 
\dynkin[name=first]{A}{o}
\node (current) at ($(first root 1)+(-0,3.4mm)$) {};
\dynkin[at=(current),name=second]{A}{o}
\node (currenta) at ($(first root 1)+(-3.4mm,0)$) {};
\dynkin[at=(currenta),name=third]{A}{o}
\node (currenta) at ($(third root 1)+(-0,3.4mm)$) {};
\dynkin[at=(currenta),name=fourth]{A}{o}
\begin{scope}[on background layer]
\draw[/Dynkin diagram/fold style]
($(first root 1)$) -- ($(second root 1)$);%
\draw[/Dynkin diagram/fold style]
($(second root 1)$) -- ($(fourth root 1)$);%
\draw[/Dynkin diagram/fold style]
($(third root 1)$) -- ($(fourth root 1)$);%
\draw[/Dynkin diagram/fold style]
($(third root 1)$) -- ($(first root 1)$);%
\end{scope}
\end{tikzpicture}\!\!\!$}  & \multirow{5}{*}{$\varnothing$}  & \multirow{5}{*}{$A_1^3$} &  $\dynkin D{oooo}$ & $\checkmark$ & $\xmark$ & $\checkmark$ \topstrut  \\
&& & & & &  $\dynkin D{oo**}$ & $\xmark$ & $\xmark$ & $?$ 
\topstrut  \\ 
&& & & & &  $\dynkin D{*o**}$ & $\xmark$ & $\xmark$ & $\xmark$ \topstrut \\ \cline{3-11}
&& \multirow{5}{*}{\makecell{$\Z_4$ or $\mathcal{D}i_8$ \\ 3 $\Stab_W(\Delta)$-c.c. \\ of each}} & \multirow{5}{*}{$\!\!\!\begin{tikzpicture}[baseline=1ex] 
\dynkin[name=first]{A}{o}
\node (current) at ($(first root 1)+(-0,3.4mm)$) {};
\dynkin[at=(current),name=second]{A}{o}
\node (currenta) at ($(first root 1)+(-3.4mm,0)$) {};
\dynkin[at=(currenta),name=third]{A}{o}
\node (currenta) at ($(third root 1)+(-0,3.4mm)$) {};
\dynkin[at=(currenta),name=fourth]{A}{o}
\begin{scope}[on background layer]
\draw[/Dynkin diagram/fold style]
($(first root 1)$) -- ($(second root 1)$);%
\draw[/Dynkin diagram/fold style]
($(second root 1)$) -- ($(fourth root 1)$);%
\draw[/Dynkin diagram/fold style]
($(third root 1)$) -- ($(fourth root 1)$);%
\draw[/Dynkin diagram/fold style]
($(third root 1)$) -- ($(first root 1)$);%
\end{scope}
\end{tikzpicture}\!\!\!$}  & \multirow{5}{*}{$\varnothing$}  & \multirow{5}{*}{$A_1^3$} &  $\dynkin [fold]D{oooo}$ & $\checkmark$ & $\xmark$ & $\checkmark$ \topstrut  \\
&& & & & &  $\dynkin [fold]D{oo**}$ & $\xmark$ & $\xmark$ & $\xmark$ \topstrut  \\ 
&& & & & &  $\dynkin [fold]D{*o**}$ & $\xmark$ & $\xmark$ & $?$ 
\topstrut \\ \cline{3-11}
&& \multirow{2.3}{*}{\makecell{~\\[-0.3cm]$ \operatorname{Alt}_4$ or $S_4$}} &  \multirow{3}{*}{$\!\!\!\begin{tikzpicture}[baseline=1ex] 
\dynkin[name=first]{A}{o}
\node (current) at ($(first root 1)+(-0,3.4mm)$) {};
\dynkin[at=(current),name=second]{A}{o}
\node (currenta) at ($(first root 1)+(-3.4mm,0)$) {};
\dynkin[at=(currenta),name=third]{A}{o}
\node (currenta) at ($(third root 1)+(-0,3.4mm)$) {};
\dynkin[at=(currenta),name=fourth]{A}{o}
\begin{scope}[on background layer]
\draw[/Dynkin diagram/fold style]
($(first root 1)$) -- ($(second root 1)$);%
\draw[/Dynkin diagram/fold style]
($(second root 1)$) -- ($(fourth root 1)$);%
\draw[/Dynkin diagram/fold style]
($(third root 1)$) -- ($(fourth root 1)$);%
\draw[/Dynkin diagram/fold style]
($(third root 1)$) -- ($(first root 1)$);%
\end{scope}
\end{tikzpicture}\!\!\!$} & \multirow{3}{*}{$\varnothing$}  & \multirow{3}{*}{$A_1^3$} &  $\begin{tikzpicture}[scale=0.8,baseline=-2.2ex] \dynkin[ply=3,rotate=270,fold radius = 2mm]{D}{oooo} \end{tikzpicture}$ & $\checkmark$ & $\xmark$ & $\checkmark$ \topstrut  \\ 
&& & & & &  $\begin{tikzpicture}[scale=0.8,baseline=-2.2ex] \dynkin[ply=3,rotate=270,fold radius = 2mm]{D}{*o**} \end{tikzpicture}$ & $\xmark$ & $\xmark$ & $\xmark$ \topstrut  \\\cline{3-11}
\pagebreak
\multirow{4.5}{*}{$A_1^4$}&\multirow{4.5}{*}{$S_4$}& \multirow{4.5}{*}{$\Z_3$ or $S_3$} &  $\!\!\!\begin{tikzpicture}[baseline=1ex] 
\dynkin[name=first]{A}{o}

\node (current) at ($(first root 1)+(-0,3.4mm)$) {};
\dynkin[at=(current),name=second]{A}{o}

\node (currenta) at ($(first root 1)+(-3.4mm,0)$) {};
\dynkin[at=(currenta),name=third]{A}{o}
\node (currenta) at ($(third root 1)+(-0,3.4mm)$) {};
\dynkin[at=(currenta),name=fourth]{A}{o}

\begin{scope}[on background layer]

\draw[/Dynkin diagram/fold style]
($(third root 1)$) -- ($(second root 1)$);%
\draw[/Dynkin diagram/fold style]
($(second root 1)$) -- ($(fourth root 1)$);%
\draw[/Dynkin diagram/fold style]
($(third root 1)$) -- ($(fourth root 1)$);%

\end{scope}
\end{tikzpicture}\!\!\!$  & $\varnothing$  & $\varnothing$ &  $\begin{tikzpicture}[scale=0.8,baseline=-2.2ex] \dynkin[ply=3,rotate=270,fold radius = 2mm]{D}{oooo} \end{tikzpicture}$ & $\checkmark$&$\xmark$&$\checkmark$ \topstrut  \\ \cdashline{4-11}
&&& $\!\!\!\begin{tikzpicture}[baseline=1ex] 
\dynkin[name=first]{A}{o}

\node (current) at ($(first root 1)+(-0,3.4mm)$) {};
\dynkin[at=(current),name=second]{A}{*}

\node (currenta) at ($(first root 1)+(-3.4mm,0)$) {};
\dynkin[at=(currenta),name=third]{A}{*}
\node (currenta) at ($(third root 1)+(-0,3.4mm)$) {};
\dynkin[at=(currenta),name=fourth]{A}{*}

\begin{scope}[on background layer]

\draw[/Dynkin diagram/fold style]
($(third root 1)$) -- ($(second root 1)$);%
\draw[/Dynkin diagram/fold style]
($(second root 1)$) -- ($(fourth root 1)$);%
\draw[/Dynkin diagram/fold style]
($(third root 1)$) -- ($(fourth root 1)$);%
\end{scope}
\end{tikzpicture}\!\!\!$ & $A_1^3$ & $A_1^3$ & $\begin{tikzpicture}[scale=0.8,baseline=-2.2ex] \dynkin[ply=3,rotate=270,fold radius = 2mm]{D}{*o**} \end{tikzpicture}$ & $\xmark$&$\xmark$&$\xmark$ \topstrut  \\ \cdashline{4-11}
&&& $\!\!\!\begin{tikzpicture}[baseline=1ex] 
\dynkin[name=first]{A}{*}

\node (current) at ($(first root 1)+(-0,3.4mm)$) {};
\dynkin[at=(current),name=second]{A}{o}

\node (currenta) at ($(first root 1)+(-3.4mm,0)$) {};
\dynkin[at=(currenta),name=third]{A}{o}
\node (currenta) at ($(third root 1)+(-0,3.4mm)$) {};
\dynkin[at=(currenta),name=fourth]{A}{o}

\begin{scope}[on background layer]

\draw[/Dynkin diagram/fold style]
($(third root 1)$) -- ($(second root 1)$);%
\draw[/Dynkin diagram/fold style]
($(second root 1)$) -- ($(fourth root 1)$);%
\draw[/Dynkin diagram/fold style]
($(third root 1)$) -- ($(fourth root 1)$);%
\end{scope}
\end{tikzpicture}\!\!\!$ & $A_1$ & $A_1$  &&&& \topstrut \\ \hline 

\multirow{34}{*}{\makecell{$A_3$ \\ \!3 c.c.\!}} & \multirow{34}{*}{$\Z_2^2$} & \multirow{12}{*}{$\Z_2^{(1)}$} & \multirow{3}{*}{$\begin{tikzpicture}[baseline=-0.3ex]
\dynkin[ply=2,rotate=90,fold radius = 2mm]{A}{ooo}\end{tikzpicture}\times\mathcal{T}^1_0$} & \multirow{3}{*}{$\varnothing$} & \multirow{3}{*}{$A_1^2$} & $\dynkin{D}{oooo}$ & $\checkmark$ & $\checkmark$ & $\checkmark$ \\
&&&&&& $\dynkin{D}{oo**}$ & $\xmark$ & $\checkmark$ & $?$ 
\\ \cdashline{4-11}
&&& \multirow{5}{*}{$\begin{tikzpicture}[baseline=-0.3ex]
\dynkin[ply=2,rotate=90,fold radius = 2mm]{A}{*o*}\end{tikzpicture}\times\mathcal{T}^1_0$} & \multirow{5}{*}{$A_1^2$} & \multirow{5}{*}{$A_3$}  & $\dynkin{D}{oooo}$ & $\xmark$ & $\xmark$ & $\xmark$ \\
&&&&&& $\dynkin{D}{oo**}$ & $\xmark$ & $\xmark$ & $\xmark$ \\ 
&&&&&& $\dynkin{D}{o***}$ & $\xmark$ & $\xmark$ & $\xmark$ \\
\cdashline{4-11} 
&&& \multirow{2.3}{*}{$\begin{tikzpicture}[baseline=-0.3ex]
\dynkin[ply=2,rotate=90,fold radius = 2mm]{A}{o*o}\end{tikzpicture}\times\mathcal{T}^1_0$} & \multirow{2.3}{*}{$A_1$} & \multirow{2.3}{*}{$A_1^3$} & $\dynkin{D}{oo**}$ & $\xmark$ & $\checkmark$ & $?$ 
\\
&&&&&& $\dynkin{D}{*o**}$ & $\xmark$ & $\xmark$ & $\xmark$ \\
\cline{3-11}

 &  & \multirow{7}{*}{$\Z_2^{(2)}$} & $
\dynkin{A}{ooo}\times\mathcal{T}^1_0$ & $\varnothing$ & $\varnothing$ & $\dynkin[fold]{D}{oooo}$ & $\checkmark$ & $\checkmark$ & $\checkmark$ \\ \cdashline{4-11} 
&&&\multirow{5}{*}{$\dynkin{A}{*o*}\times\mathcal{T}^1_0$} & \multirow{5}{*}{$A_1^2$} & \multirow{5}{*}{$A_3$} & $\dynkin[fold]{D}{oooo}$ & $\xmark$ & $\checkmark$ & $?$ \\
&&&&&& $\dynkin[fold]{D}{oo**}$ & $\xmark$ & $\xmark$ & $\xmark$ \\ 
&&&&&& $\dynkin[fold]{D}{o***}$ & $\xmark$ & $\checkmark$ & $\xmark$ \\\cline{3-11}

&& \multirow{11}{*}{$\Z_2^2$} & \multirow{3}{*}{$\begin{tikzpicture}[baseline=-0.3ex]
\dynkin[ply=2,rotate=90,fold radius = 2mm]{A}{ooo}\end{tikzpicture}\times\mathcal{T}^1_0$} & \multirow{3}{*}{$\varnothing$} & \multirow{3}{*}{$A_1^2$}  & $\dynkin[fold]{D}{oooo}$ & $\checkmark$ & $\xmark$ & $\checkmark$ \\
&&&&&& $\dynkin[fold]{D}{oo**}$ & $\xmark$ & $\xmark$ & $\xmark$ \\ \cdashline{4-11} 

 &&& \multirow{5}{*}{$\begin{tikzpicture}[baseline=-0.3ex]
\dynkin[ply=2,rotate=90,fold radius = 2mm]{A}{*o*}\end{tikzpicture}\times\mathcal{T}^1_0$} & \multirow{5}{*}{$A_1^2$} & \multirow{5}{*}{$A_3$} & $\dynkin[fold]{D}{oooo}$ & $\xmark$ & $\xmark$ & $\xmark$ \\ 
&&&&&& $\dynkin[fold]{D}{oo**}$ & $\xmark$ & $\xmark$ & $\xmark$ \\ 
 &&&&&& $\dynkin[fold]{D}{o***}$ & $\xmark$ & $\xmark$ & $\xmark$ \\ \cdashline{4-11}

&&& $\begin{tikzpicture}[baseline=-0.3ex]
\dynkin[ply=2,rotate=90,fold radius = 2mm]{A}{o*o}\end{tikzpicture}\times\mathcal{T}^1_0$ & $A_1$ & $A_1^3$ & $\dynkin[fold]{D}{*o**}$ & $\xmark$ & $\xmark$ & $?$ 
\\\hline

\multirow{4.6}{*}{$A_2$} & \multirow{4.6}{*}{$\mathcal{D}i_{12}$} &$\Z_3$, $S_3^{(1)}$ & $\dynkin{A}{oo} \times\mathcal{T}^2_0$ & $\varnothing$ & $\varnothing$ & \begin{tikzpicture}[scale=0.8,baseline=-2.2ex] \dynkin[ply=3,rotate=270,fold radius = 2mm]{D}{oooo} \end{tikzpicture} & $\checkmark$ & $\xmark$ & $\checkmark$ \\ \cline{3-11}
&& \multirow{2.7}{*}{$\Z_6$, $S_3^{(2)}$ or $\mathcal{D}i_{12}$} & \multirow{2.7}{*}{$\begin{tikzpicture}[baseline=-0.3ex]
\dynkin[ply=2,rotate=90,fold radius = 2mm]{A}{oo}\end{tikzpicture}\times\mathcal{T}^2_0$} & \multirow{2.7}{*}{$\varnothing$} & \multirow{2.7}{*}{$A_1^3$} & \begin{tikzpicture}[scale=0.8,baseline=-2.2ex] \dynkin[ply=3,rotate=270,fold radius = 2mm]{D}{oooo} \end{tikzpicture} & $\checkmark$ & $\xmark$ & $\checkmark$ \\
&&&&&& \begin{tikzpicture}[scale=0.8,baseline=-2.2ex] \dynkin[ply=3,rotate=270,fold radius = 2mm]{D}{*o**} \end{tikzpicture} & $\xmark$ & $\xmark$ & $\xmark$ \\ \hline
\end{longtable}\restoregeometry
\end{small}
\newpage\restoregeometry
\subsubsection{\texorpdfstring{$\Phi=D_n, n>4$}{Dn}}
\noindent Let $\Phi=D_n , n > 4$. We use the following \textit{standard} parametrisation in $\R^n$ of $\Phi$. Let $\alpha_i = e_i-e_{i+1}$ for $1\leq i \leq n-1$ and $\alpha_n = e_{n-1}+e_n$. 

By Table II of \cite{T2}, the isotropic indices of type $\Phi$ consist of the following diagrams:\\[0.2cm]
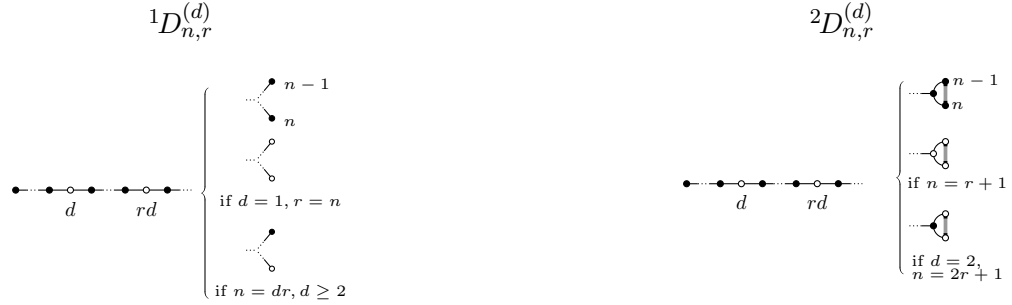
\begin{tabular}{cp{1.5cm}c}
$^{1\!}D^{(d)}_{n,r}$ & & $^{2\!}D^{(d)}_{n,r}$ \\[0.2cm]
\makecell{\begin{tikzpicture}[scale=0.8,baseline=-5mm] 
\node (up) at (0,2) {};
\dynkin[at=(up),name=upper,labels={,,d,,,rd,}]{A}{*.*o*.*o*}
\draw ($(upper root 7)$) -- ($(upper root 7)+(0.2,0)$);
\draw[densely dotted] ($(upper root 7)+(0.2,0)$) -- ($(upper root 7)+(0.4,0)$);
\node (up1) at ($(up)+(4,1.5)$) {};
\node (up2) at ($(up1)+(0,-1)$) {};
\node (up3) at ($(up2)+(0,-1.5)$) {};
\node[inner sep=0pt] (d1) at (up1) {};
\node[label=right:{\tiny $n-1$}, inner sep=0pt] (d2) at ($(d1)+(50:0.4)$) {};
\node[label=right:{\tiny $n$}, inner sep=0pt] (d3) at ($(d1)+(310:0.4)$) {};
\draw[fill=black] (d2) circle [radius=1.3pt];
\draw[fill=black] (d3) circle [radius=1.3pt];
\draw[densely dotted] (d1) -- ($(d1)+(-2mm,0)$);
\draw[densely dotted] (d1) -- ($(d1)+(50:0.2)$);
\draw[densely dotted] (d1) -- ($(d1)+(310:0.2)$);
\draw ($(d1)+(50:0.2)$)--(d2);
\draw ($(d1)+(310:0.2)$)--(d3);
\node[inner sep=0pt] (d4) at (up2) {};
\node[inner sep=0pt] (d5) at ($(d4)+(50:0.4)$) {};
\node[inner sep=0pt] (d6) at ($(d4)+(310:0.4)$) {};
\draw (d5) circle [radius=1.2pt];
\draw (d6) circle [radius=1.2pt];
\draw[densely dotted] (d4) -- ($(d4)+(-2mm,0)$);
\draw[densely dotted] (d4) -- ($(d4)+(50:0.2)$);
\draw[densely dotted] (d4) -- ($(d4)+(310:0.2)$);
\draw ($(d4)+(50:0.2)$)--(d5);
\draw ($(d4)+(310:0.2)$)--(d6);
\node[inner sep=0pt] (d7) at (up3) {};
\node[inner sep=0pt] (d8) at ($(d7)+(50:0.4)$) {};
\node[inner sep=0pt] (d9) at ($(d7)+(310:0.4)$) {};
\draw[fill=black] (d8) circle [radius=1.2pt];
\draw (d9) circle [radius=1.2pt];
\draw[densely dotted] (d7) -- ($(d7)+(-2mm,0)$);
\draw[densely dotted] (d7) -- ($(d7)+(50:0.2)$);
\draw[densely dotted] (d7) -- ($(d7)+(310:0.2)$);
\draw ($(d7)+(50:0.2)$)--(d8);
\draw ($(d7)+(310:0.2)$)--(d9);
\draw [decorate, decoration = {calligraphic brace}] ($(d7)+(-8mm,-8mm)$) --  ($(d1)+(-8mm,2mm)$);
\node[align=left] at ($(up2) + (4mm,-7mm)$) {\tiny{if $d=1, r= n$}};
\node[align=left] at ($(up3) + (4mm,-7mm)$) {\tiny{if $n=dr, d\geq 2$}};
\end{tikzpicture}} & &
\makecell{\begin{tikzpicture}[scale=0.8,baseline=-5mm] 
\node (up) at (0,2) {};
\dynkin[at=(up),name=upper,labels={,,d,,,rd,}]{A}{*.*o*.*o*}
\draw ($(upper root 7)$) -- ($(upper root 7)+(0.2,0)$);
\draw[densely dotted] ($(upper root 7)+(0.2,0)$) -- ($(upper root 7)+(0.4,0)$);
\node (up1) at ($(up)+(4.3,1.5)$) {};
\node (up2) at ($(up1)+(0,-1)$) {};
\node (up3) at ($(up2)+(0,-1.2)$) {};
\begin{scope}
\dynkin[at =(up1), name=upper1,labels={,,},ply=2,fold radius=2mm, rotate = 180]{A}{***}
\draw[densely dotted] ($(upper1 root 2)+(0.4,0)$) -- ($(upper1 root 2)+(0.2,0)$);
\draw ($(upper1 root 2)+(0.2,0)$) -- ($(upper1 root 2)$);
\node at ($(upper1 root 3) + (-0.5,0)$) {\tiny{$n-1$}};
\node[align=left] at ($(upper1 root 1) + (-0.2,0)$) {\tiny{$n$}};
\end{scope}
\begin{scope}
\dynkin[at=(up2),name=middle1_right, labels={,,},ply=2,fold radius=2mm, rotate = 180]{A}{ooo}
\draw[densely dotted] ($(middle1_right root 2)+(0.4,0)$) -- ($(middle1_right root 2)+(0.2,0)$);
\draw ($(middle1_right root 2)+(0.2,0)$) -- ($(middle1_right root 2)+(0.05,0)$);
\end{scope}
\begin{scope}
\dynkin[at=(up3),name=lower1_right, labels={,,},ply=2,fold radius=2mm, rotate = 180]{A}{o*o}
\draw[densely dotted] ($(lower1_right root 2)+(0.4,0)$) -- ($(lower1_right root 2)+(0.2,0)$);
\draw ($(lower1_right root 2)+(0.2,0)$) -- ($(lower1_right root 2)$);
\end{scope}
\draw [decorate, decoration = {calligraphic brace}] ($(up3)+(-7mm,-8mm)$) --  ($(up1)+(-7mm,2mm)$);
\node[align=left] at ($(up2) + (2mm,-5mm)$) {\tiny{if $n=r+1$}};
\node[align=left] at ($(up3) + (2mm,-7mm)$) {\tiny{if $d=2,$}\\[-0.01cm] \tiny{$n=2r+1$}};
\end{tikzpicture}}\\
$d=2^a|2n$ for $a \in \mathbb{N} \cup \{0\}$, $rd\leq n, n \neq rd+1$ & & \hspace{-9mm}$d=2^a|2n$ for $a \in \mathbb{N}\cup \{0\}$, $rd\leq n-1$
\end{tabular}
~\\[0.2cm] We follow the convention that anisotropic indices of type $D_n$ satisfy $r=0$, and $d$ may be any positive integer. 

Recall from Table \ref{classlist} that $\Delta$ is one of the following types: 
\begin{enumerate}
\item $A_{n-1}$ (of which there are 2 $W$-conjugacy classes if $n$ is even, 1 otherwise).
\item $D_mD_{n-m}$ for $1 \leq m \leq \frac{n}{2}$.
\item $D_m^{n_m}$ for $mn_m=n$ and $n_m\geq 2$.
\end{enumerate}

We start by considering $\Delta= A_{n-1}$. We have that $\Stab_{\Iso(\Phi)}(\Delta)=\Z_2$ by Table \ref{classlist}. We may ignore the case $\Pi = 1$, as any associated index $\mathcal{H}$ is not maximal in $\Phi$. So we need only consider the case $\Pi = \Z_2$. Then $\mathcal{H}$ is given as in Table \ref{D_n}. Note that here $\Delta_0=\left\{ \alpha_i, \alpha_{n-i} \mid 1\leq i \leq \lfloor\frac{n}{2}\rfloor, d\nmid i \right\}$, so 
  \begin{align*}
      E_{\Delta_0} = &\spn[\R]{e_{i}-e_{i+1}\mid
    1\leq i \leq rd, d\nmid i}\\
  \oplus &\spn[\R]{e_{i}-e_{i+1} \mid rd+1\leq i \leq m-rd}\\
  \oplus &\spn[\R]{e_{i}-e_{i+1} \mid m-rd+1\leq i \leq m, d\nmid i }
    \end{align*} as for type $A_n$. We can use the same computations as in Appendix \ref{An} for $\Delta = A_m^{n_m}$ and $\Pi \nsubseteq W$ where $\rho(\Pi)$ is a primitive subgroup of $S_{n_m}$ and $\rho: \Z_2\times S_{n_m}\rightarrow S_{n_m}$ is the projection. We get a basis $B$ for $E_a$ given by 
    \begin{align*}
    B_{1,1} & := \left\{e_{i}-e_{i+1}| 1\leq i \leq rd, d \nmid i \right\}\\
B_{1,2} & := \left\{e_{n-i}-e_{n-i+1} |
    1\leq i \leq rd, d \nmid i \right\}\\
B_2 & := \left\{e_{i} |
    rd +1 \leq i \leq n-rd \right\}\\
B_3 & := \left\{v_{i,0}^{(d)}+v_{n/d-i+1,0}^{(d)} \Big|1\leq i\leq r \right\}\\
B & := B_{1,1}\cup B_{1,2}\cup B_{2} \cup B_{3}.
\end{align*}
Note that we have 
\begin{align*}
    e_{(i-1)d+1}-e_{n-(i-1)d} & = \frac{1}{d}\left(v_{i,0}^{(d)}+v_{\frac{n}{d}-i+1,0}^{(d)}+\sum_{j=1}^{d-1} j(e_{id-j}-e_{id-j+1}) + j(e_{n-id+j+1}-e_{n-id+j})\right)
\end{align*}
so for $B_3':=\{e_{(i-1)d+1}-e_{n-(i-1)d}\mid 1\leq i\leq r\}$ we get another basis $B':=B_{1,1}\cup B_{1,2}\cup B_{2}\cup B'_{3}$. Then $B_{1,1}\cup B_{1,2}\cup B'_{3}$ gives rise to a root system of type $A_{2d-1}^r$ and $B_2$ gives rise to a root system of type $D_{n-2rd}$, that is $\Phi_a = A_{2d-1}^r\+  D_{n-2rd}$. \\

Next, let $\Delta = D_mD_{n-m}$ where $1\leq m \leq \frac{n}{2}$. We have that $S=\Z_2^2$ by Table \ref{classlist}. We start with $\Pi = 1$. In this case, $$\Delta_0 = \left\{\alpha_i \middle| \begin{array}{l}
   0\leq i\leq m-1,  \\
     m-i+1\neq sd, 1\leq s \leq r
\end{array} \right\} \cup \left\{\alpha_i \middle| \begin{array}{l}
   m+1\leq i\leq n,  \\
     i-m\neq sd', 1\leq s\leq r'
\end{array} \right\},$$ where $\alpha_0$ denotes the longest root in $\Phi = D_n$.\\
We have $\dim E_s = r+r'$. Note that $$\Delta_s:=\{\alpha_{m-(sd-1)}\mid 1\leq s\leq r\}\cup\{\alpha_{m+sd'}\mid 1\leq s\leq r'\}\subseteq E_s$$ is a set of $r+r'$ linear independent vectors and therefore spans $E_s$.  In particular, they also span $\Phi_s = \Phi\cap E_s.$ As in Appendix \ref{Cn} for $\Delta = C_mC_{n-m}$, it follows that $$\Phi_a = (\Phi_s)^{\perp}\simeq A_{d-1}^r\+A_{d'-1}^{r'}\+ D_{n-(rd+r'd')}.$$
If $\Pi= \Z_2^{(i)}$ for $i=2,3$, then either $\Delta_s$ is as before (this occurs when $d'\neq 1,2$), if $d'= 2$ and $n-j = 2r'+1$ we get $$\Delta_s:=\{\alpha_{j-(sd-1)}\mid 1\leq s\leq r\}\cup\{\alpha_{j+sd'}\mid 1\leq s\leq r'-1\}\cup\{\alpha_{n-1}+\alpha_{n}\}.$$ In the case of $d'=1$, $n-j=r'+1$ we have $$\Delta_s:=\{\alpha_{j-(sd-1)}\mid 1\leq s\leq r\}\cup \{\alpha_i\mid j+1\leq i\leq n-2\}\cup\{\alpha_{n-1}+\alpha_n\}.$$ Again, we note that in each case, $\Delta_s$ spans $E_s$ and therefore is a base for $\Phi_s$ as well, and it follows that $$\Phi_a = (\Phi_s)^{\perp}\simeq A_{d-1}^r\+A_{d'-1}^{r'}\+ D_{n-(rd+r'd')}.$$
Finally, we get a similar situation for $\Pi=\Z_2^{(1)}, \Z_2^2$. First, consider $d,d'\neq 1,2$. Then, as before 
$$\Delta_s:=\{\alpha_{j-(sd-1)}\mid 1\leq s\leq r\}\cup\{\alpha_{j+sd'}\mid 1\leq s\leq r'\}.$$
If $d = 1$ and $j = r+1$ the first set is just $\{\alpha_i \mid 2\leq i \leq j-1\}\cup \{\alpha_0+\alpha_1\}$ and for $d = 2$, $j=2r+1$ we get $\{\alpha_{j-(sd-1)}\mid 1\leq s\leq r-1\}\cup \{\alpha_0+\alpha_1\}.$ For $d'$ the same as the previous case holds, i.e. if $d'= 2$ and $n-j = 2r'+1$ the second set is $\{\alpha_{j+sd'}\mid 1\leq s\leq r'-1\}\cup\{\alpha_{n-1}+\alpha_{n}\}$ and if $d'=1$, $n-j=r'+1$ we have $\{\alpha_i\mid j+1\leq i\leq n-2\}\cup\{\alpha_{n-1}+\alpha_n\}.$ 
As before, we have $$\Phi_a = (\Phi_s)^{\perp}\simeq A_{d-1}^r\+A_{d'-1}^{r'}\+ D_{n-(rd+r'd')}.$$\\


We next consider $\Delta= D_m^{n_m}$. In this case, $\Stab_{\Iso{\Phi}}(\Delta)=N\rtimes Q$, where $N:=\Z_2^{n_m}$ and $Q:=S_{n_m}$. Let $\rho: N\rtimes Q \rightarrow Q$ be the projection.\\
First, consider the case for $\Pi \cap N = 1$, and the index is given as in Table \ref{D_n}. Let $$\Delta_0= \Delta - \cup_{i=0}^{n_m-1}\{\alpha_{im+d},...\alpha_{im+rd} \} \text{ if } m \neq rd$$
$$\Delta_0= \Delta - \cup_{i=0}^{n_m-1}\{\alpha_{im+d},...,\alpha_{im+(r-1)d}, \alpha_{(i+1)m-1}^*\} \text{ if } m=rd,$$
where $\alpha_{k}^* = \alpha_{k}+\alpha_{n-1}+\alpha_{n}+2\sum_{j=0}^{n-k-3} \alpha_{k+1+j}=e_k+e_{k+1}$.
The subspace $E_{\Delta_0}$ has dimension $n-n_mr$ and therefore $E_{\Delta_0}^\perp$ has dimension $n_mr$ and we compute that:
$$E_{\Delta_0}^\perp = \oplus_{c=0}^{r-1} \oplus_{b=0}^{n_m-1} \left\langle \sum_{f=1}^{d} e_{bm+cd+f} \right\rangle.$$

The action of $\Pi$ on the basis vectors $e_i$ is given here: $\Pi(e_i)=e_{i+m}$ for all $1 \leq i \leq n-m$, $\Pi(e_i)=e_{i-n+m}$ for all $n-m < i < n$ and $\Pi(e_n)=-e_m$. 
Hence the fixed point space $E_s$ of dimension $r$ is computed:
$$E_s= \oplus_{c=0}^{r-1} \left\langle \sum_{f=1}^{d}  (\sum_{b=0}^{n_m-1} e_{bm+cd+f}) \right\rangle.$$

Then it follows that $E_a$ of dimension $n-r$ is given by $E_a=B_1+B_2$ where:
$$B_1=\oplus_{c=0}^{r-1} (\oplus_{f=2}^{d} \left\langle e_{cd+1}-e_{m+cd+f} \right\rangle \oplus_{b=1}^{n_m-1} \oplus_{f=1}^{d} \left\langle e_{cd+1}-e_{bm+cd+f}\right\rangle),$$
$$ B_2 = \oplus_{f=1}^{m-rd} \oplus_{b=0}^{n_m-1} \langle e_{bm+rd+f} \rangle.$$
$B_1$ gives rise to a root system of type $A_{dn_m-1}^{r}$ and analogously $B_2$ gives rise to a root system $D_{n_m(m-rd)}$. Then it follows that: $$\Phi_a=E_a \cap \Phi \cong A_{dn_m-1}^{r}D_{n_m(m-rd)}.$$

If we now consider the case where $\Pi\cap N \neq 1$, then the only cases that differ from the previous case is when $m=rd+1$ and $d=1$ or $2$. Assume that $m=dr+1$ and $d=1$ or $2$ then $\Delta_0= \Delta - \cup_{i=0}^{n_m-1}\{\alpha_i+d,...\alpha_{i+rd},\alpha_{i+rd}^*\}$;
where $\alpha_{k}^* := \alpha_{k}+\alpha_{n-1}+\alpha_{n}+2\sum_{j=0}^{n-k-3} \alpha_{k+1+j}=e_k+e_{k+1}$.

The subspace $E_{\Delta_0}$ has dimension $n-n_m(r+1)=rn_m(d-1)$ and therefore $E_{\Delta_0}^\perp$ has dimension $n_m(r+1)$ and we compute that:
$$E_{\Delta_0}^\perp = \oplus_{c=0}^{r-1} \oplus_{b=0}^{n_m-1} \left\langle \sum_{f=1}^{d} e_{bm+cd+f} \right\rangle \oplus_{b=0}^{n_m-1} \left\langle e_{(b+1)m} \right\rangle .$$

The action of $\Pi$ on the basis vectors $e_i$ is the same as the previous case so $\Pi(e_i)=e_{i+m}$ for all $1 \leq i \leq n-m$, $\Pi(e_i)=e_{i-n+m}$ for all $n-m < i < n$ and $\Pi(e_n)=-e_m$. 
Hence the fixed point space $E_s$ of dimension $r$ is computed:
$$E_s= \oplus_{c=0}^{r-1} \left\langle \sum_{f=1}^{d}  (\sum_{b=0}^{n_m-1} e_{bm+cd+f}) \right\rangle.$$

Then it follows that $E_a$ is the same as the previous case and given by $E_a=B_1+B_2$ where:
$$B_1=\oplus_{c=0}^{r-1} (\oplus_{f=2}^{d} \left\langle e_{cd+1}-e_{m+cd+f} \right\rangle \oplus_{b=1}^{n_m-1} \oplus_{f=1}^{d} \left\langle e_{cd+1}-e_{bm+cd+f}\right\rangle),$$
$$ B_2 = \oplus_{f=1}^{m-rd} \oplus_{b=0}^{n_m-1} \left\langle e_{bm+rd+f} \right\rangle.$$
Again $B_1$ gives rise to a root system of type $A_{dn_m-1}^{r}$ and analogously $B_2$ gives rise to a root system $D_{n_m(m-rd)}$. Then it follows again that: $$\Phi_a=E_a \cap \Phi \cong A_{dn_m-1}^{r}D_{n_m(m-rd)}.$$

\newpage\newgeometry{top=0.5cm,bottom=0.4cm,left=3cm, right=0.5cm,foot=0cm}
\begin{landscape}
\begin{small}
\begingroup
\begin{longtable}{| c | c | c | c | c | c | c | c | c | c |}\caption{Isotropic maximal embeddings of abstract indices of type $D_n$, $n>4$}\label{D_n} \\ \hline

\multicolumn{1}{|c|}{} & \multicolumn{1}{c|}{} & \multicolumn{1}{c|}{} & \multicolumn{1}{c|}{} & \multicolumn{1}{c|}{} & \multicolumn{1}{c|}{} & \multicolumn{1}{c|}{} & \multicolumn{3}{c|}{\small{Special fields}} \\

  $\Delta$ & \small{$\hspace{-0.5mm}\Stab_{\operatorname{Iso}(\Phi)}(\Delta)\hspace{-0.5mm}$} & $\Pi$ & $\textcolor{white}{iiaaa} \mathcal{H} \textcolor{white}{aaaii}$ & $\Psi_0$ & $\Phi_a$ & $\mathcal{G}$ & \small{$\hspace{-0.8mm}\!\cd 1 \!\hspace{-1mm}$} & $\!\R\!$ & $\!\!Q_{\mathfrak{p}}\!\!$ \\ \hline \hline 
\endfirsthead

\multicolumn{10}{|c|}%
{{\bfseries \tablename\ \thetable{} -- continued from previous page}} \\
\hline
\multicolumn{1}{|c|}{} & \multicolumn{1}{c|}{} & \multicolumn{1}{c|}{} & \multicolumn{1}{c|}{} & \multicolumn{1}{c|}{} & \multicolumn{1}{c|}{} & \multicolumn{1}{c|}{} & \multicolumn{3}{c|}{\small{Special fields}} \\

  $\Delta$ & \small{$\hspace{-0.5mm}\Stab_{\Iso(\Phi)}(\Delta)\hspace{-0.5mm}$} & $\Pi$ & $\mathcal{H}$ & $\Psi_0$ & $\Phi_a$ & $\mathcal{G}$ & \small{$\hspace{-0.8mm}\!\cd 1 \!\hspace{-1mm}$} & $\!\R\!$ & $\!\!Q_{\mathfrak{p}}\!\!$  \\ \hline \hline \endhead
\hline \multicolumn{10}{|c|}{{Continued on next page}} \\ \hline \endfoot 
 \endlastfoot
\hline
\makecell{$A_{n-1}$ \\\!\!\hspace{-0.2mm}$(n,2)$\hspace{-0.2mm} c.c.\!\!} & \makecell{$\Z_2$} & \makecell{$\Z_2$} &%
\makecell{\begin{tikzpicture}[scale=0.8,baseline=-5mm] 
\node at (0,1.5) {};
\node at (0,-2) {};
\node (up) at (0,0) {};
\dynkin[at=(up),name=upper,labels={,,d,,,rd,,,,{n-rd},,,,,,n-1}, ply=2,fold radius=2mm]{A}{*.*o*.*o*.*.*.*o*.*o*.*};
\node (up1) at ($(up)+(4.5,1.05)$) {};
\node (up2) at ($(up1)+(0,-0.75)$) {};
\node (up3) at ($(up2)+(0,-0.75)$) {};
\dynkin[at=(up1),name=upper_right, labels={,},ply=2,fold radius=2mm]{A}{.**.};
\draw[densely dotted] ($(upper_right root 1)+(-0.4,0)$) -- ($(upper_right root 1)+(-0.2,0)$);
\draw ($(upper_right root 1)+(-0.2,0)$) -- ($(upper_right root 1)+(0.01,0)$);
\draw[densely dotted] ($(upper_right root 2)+(-0.4,0)$) -- ($(upper_right root 2)+(-0.2,0)$);
\draw ($(upper_right root 2)+(-0.2,0)$) -- ($(upper_right root 2)+(0.01,0)$);
\dynkin[at=(up2),name=middle_right, labels={,,}, ply=2,fold radius=2mm]{A}{.***.};
\draw[densely dotted] ($(middle_right root 1)+(-0.4,0)$) -- ($(middle_right root 1)+(-0.2,0)$);
\draw ($(middle_right root 1)+(-0.2,0)$) -- ($(middle_right root 1)+(0.01,0)$);
\draw[densely dotted] ($(middle_right root 3)+(-0.4,0)$) -- ($(middle_right root 3)+(-0.2,0)$);
\draw ($(middle_right root 3)+(-0.2,0)$) -- ($(middle_right root 3)+(0.01,0)$);
\dynkin[at=(up3),name=lower_right, labels={,,},ply=2,fold radius=2mm]{A}{.*o*.};
\draw[densely dotted] ($(lower_right root 1)+(-0.4,0)$) -- ($(lower_right root 1)+(-0.2,0)$);
\draw ($(lower_right root 1)+(-0.2,0)$) -- ($(lower_right root 1)+(0.01,0)$);
\draw[densely dotted] ($(lower_right root 3)+(-0.4,0)$) -- ($(lower_right root 3)+(-0.2,0)$);
\draw ($(lower_right root 3)+(-0.2,0)$) -- ($(lower_right root 3)+(0.01,0)$);
\draw [decorate, decoration = {calligraphic brace}] ($(lower_right root 3)+(-7mm,-0.8)$) --  ($(upper_right root 1)+(-7mm,0.1)$);
\node[align=left] at ($(up3) + (0,-5.5mm)$) {if $n$};
\node[align=left] at ($(up3) + (0.2,-0.9)$) {$ = 2rd$};
\node at ($(middle_right root 3)+(1cm,0)$) {$\+ \cT_0^1$};
\end{tikzpicture} \\
$d,r \geq 1$, $d|n$, $2rd \leq n$ \vspace{1mm}}

%
& \makecell{$A_{d-1}^{2r}$\\$\times$\\$A_{n-1-2rd}$ \\ $(A_{-1}:=\varnothing)$} & \makecell{$A_{2d-1}^{r}$\\$\times$\\$D_{n-2rd}$} & %
\makecell{\begin{tikzpicture}[scale=0.8,baseline=-5mm] \dynkin[name=dyn,labels = {,,c,,,ic,,,n-1,n}]{D}{*.*o*..*o*..***}
\end{tikzpicture} \\ only if $n$ is even \\ \begin{tikzpicture}[scale=0.8,baseline=-5mm] \dynkin[name=dyn,ply=2,labels = {,,c,,,ic,,,n-1,n}]{D}{*.*o*..*o*..***}
\end{tikzpicture} \\ only if $n$ is odd; \\ $c|2d$, $2rd \leq ic \leq 2\lfloor\frac{n}{2}\rfloor$} & \makecell{$\checkmark$ \\ iff \\ \!\!$r\!=\!\lfloor\frac{n}{2}\rfloor$\!\! \\ and \\\!\!$i\!=\!2\lfloor\frac{n}{2}\rfloor$\!\!} & \makecell{$\checkmark$ \\ iff $d\!=\!1$, and \\ $(c,i)$ is either \\ $(1,2r)$ or \\ $(2,\lfloor\frac{n}{2}\rfloor)$\!\!} & \makecell{$\checkmark$ \\ \!\!if $c\!=\!d\!=\!1$,\!\! \\ \!\!\hspace{-0.3mm}$r\!=\!\lfloor\frac{n}{2}\rfloor$, \hspace{-0.3mm}and\hspace{-0.3mm}\!\! \\ \!\!$i\!=\!2\lfloor\frac{n}{2}\rfloor$,\!\! \\ $\xmark$ \\ if $d\!\neq \!1$ \\ or $c\!>\!2$}
\\
\hline
\end{longtable}%
\endgroup
\end{small}
\end{landscape}
\restoregeometry

\bigskip\noindent
{\textbf {Acknowledgements}}:
The three authors VG, DS, and LV, are respectively affiliated with the University of Warwick, the University of Freiburg, and the Ruhr University Bochum. DS was partially supported by a Humboldt Postdoctoral Research Fellowship hosted at the Ruhr University Bochum.

\end{document}

\end{proof}

